\documentclass[12pt]{article}

\usepackage{a4wide}
\usepackage{amssymb}
\usepackage{amsfonts}
\usepackage{amsmath}
\input xy
\xyoption{arrow} \xyoption{matrix}

\date{}

\newtheorem{proposition}{Proposition}[section]
\newtheorem{theorem}[proposition]{Theorem}
\newtheorem{lemma}[proposition]{Lemma}

\newtheorem{corollary}[proposition]{Corollary}

\def\der{\partial }

\def\nFM0{{\nu }_{F,M_0}}
\def\nFN0{{\nu }_{F,N_0}}
\def\nGN0{{\nu }_{G,N_0}}

\def\N0{ {\bf N}_0 }

\def\g{\gamma}

\def\ra{\rightarrow}

\def\lra{\leftrightarrow}
\def\Xpm{X^{\pm }}

\def\s{\sigma}
\def\Z{\mathbb{Z}}

\def\l1{{\lambda}_1}

\def\a{\alpha}
\def\a0{ {\alpha }_0}
\def\a1{ {\alpha }_1}

\def\l{\lambda}
\def\o{\omega}

\def\nFGM0{{\nu }_{F,G,M_0}}

\def\nFN0{{\nu}_{F,N_0}}

\def\sm{{\sigma}^m}

\def\sm1{{\sigma}^{-1}}

\def\smtp1{{\sigma}^{-t+1}}

\def\o{\omega }
\def\S1{S^{-1}}

\def\Xpm1{X^{\pm 1}_1}

\def\sPM1{{\sigma }^{\pm 1}}
\def\sMP1{{\sigma }^{\mp 1 }}

\def\d{\delta}

\def\di{{\rm d.ind}}

\def\L{\Lambda}
\def\O{\Omega}

\def\G{\Gamma}

\def\OO{{\cal O}}

\def\Ytm1{Y^{t-1}}
\def\Yim1{Y^{i-1}}

\def\CN{{\cal N}}

\def\CF{{\cal F}}

\def\Aut{{\rm Aut}}

\def\dim{{\rm dim }}

\def\ker{ {\rm ker } }

\def\CJ{ {\cal J}}

\def\D{ \Delta }

\def\SL2Z{ {\rm SL}_2({\bf Z}) }

\def\th{ \theta }

\def\Gp1{ G^{1 , 1 } }
\def\P11{ P^{-1 , 1 } }
\def\Pp1{ P^{1 , 1 } }

\def\Supp{{\rm Supp}}

\def\th{\theta}
\def\CE{{\cal E}}

\def\nCLsr{{}^\nu\kern-2pt {\cal L}^{\sigma , \rho  }}
\def\nP{{}^\nu \kern-2pt P}
\def\nL{{}^\nu\kern-2pt L}
\def\nLL{{}^\nu\kern-2pt \Lambda}
\def\nPsr{{}^\nu\kern-2pt P^{\sigma , \rho  }}
\def\nLsr{{}^\nu\kern-2pt L^{\sigma , \rho  }}
\def\nuCL{{}^\nu\kern-2pt  {\cal L}}
\def\nCLsr{{}^\nu\kern-2pt {\cal L}^{\sigma , \rho  }}
\def\nCL1m{{}^\nu\kern-2pt {\cal L}^{-1 , 1  }}
\def\x1nu{x^\frac{1}{\nu}}
\def\xm1nu{x^{-\frac{1}{\nu}}}

\def\CN{{\cal N}}
\def\ra{\rightarrow }

\def\CB{{\cal B}}

\def\CE{ {\cal E} }

\def\nAM0{{\nu }_{{\cal A},M_0}}
\def\nAN0{{\nu }_{{\cal A},N_0}}

\def\End{ {\rm End }}

\def\CJ{ {\cal J }}

\def\det{ {\rm det }}

\def\tr{{\rm tr}}

\def\gm{\mathfrak{m}}

\def\hi{\widehat{i}}

\def\GL{{\rm GL}}
\def\SL{{\rm SL}}

\def\di!{\frac{\der^i}{i!}}
\def\dik!{\frac{\der^k_i}{k!}}

\def\N{\mathbb{N}}

\def\0{\overline{0}}
\def\1{\overline{1}}
\def\Lnz{\L_{n,\overline{0}}}
\def\Ln1{\L_{n,\overline{1}}}

\def\oa{\overline{a}}
\def\Lnod{\L_n^{od}}
\def\Lnev{\L_n^{ev}}

\def\a1{a_{\overline{1}}}

\def\S{\Sigma}
\def\tS{\widetilde{\Sigma}}
\def\grad{{\rm grad}}
\def\bCJ{\overline{\CJ}}

\def\od{{\rm od}}
\def\Od{{\rm Od}}

\begin{document}

\author{V. V. \  Bavula 
}

\title{The Jacobian map, the Jacobian group and the group of automorphisms of the Grassmann algebra}

\maketitle
\begin{abstract}
There are nontrivial dualities and parallels between polynomial
algebras and the Grassmann algebras. This paper is an attempt to
look at the Grassmann algebras at the angle of the Jacobian
conjecture for  polynomial algebras (which is the
question/conjecture about the $ $ {\em Jacobian set} -- the set of
all algebra endomorphisms of a polynomial algebra with the
Jacobian $1$ -- the Jacobian conjecture claims that the Jacobian
set is a {\em group}). In this paper, we study in detail  the
Jacobian set for the Grassmann algebra which turns out to be a
{\em group} -- the {\em Jacobian group} $\Sigma$ -- a
sophisticated (and large) part of the group of automorphisms of
the Grassmann algebra $\L_n$. It is proved that the Jacobian group
$\Sigma$ is a rational unipotent algebraic group. A (minimal) set
of generators for the algebraic group $\Sigma$, its dimension and
coordinates are found explicitly. In particular, for $n\geq 4$,
$$\dim (\S )=\begin{cases}
(n-1)2^{n-1} -n^2+2& \text{if $n$ is even},\\
(n-1)2^{n-1} -n^2+1& \text{if $n$ is odd}.\\
\end{cases}
 $$
The same is done for the Jacobian ascents - some natural algebraic
overgroups of $\Sigma$. It is proved that the Jacobian map $\s
\mapsto \det (\frac{\der \s (x_i)}{\der x_j})$ is surjective for
odd $n$, and is {\em not} for even $n$ though, in this case,  the
image of the Jacobian map is an algebraic subvariety of
codimension 1 given by a single equation.



$${\bf Contents}$$
\begin{enumerate}
\item Introduction. \item The group of automorphisms of the
Grassmann ring. \item The group $\G$, its subgroups, and the
Inversion Formula.\item The Jacobian group $\S$ and the equality
$\S = \S'\S''$. \item The algebraic group $\S'$ and its dimension.
\item A (minimal) set of generators for the Jacobian group $\S$
and its dimension. \item The image of the Jacobian map, the
dimensions of the Jacobian ascents and of $\G / \S$.  \item
Analogues of the Poincar\'{e} Lemma. \item The unique presentation
$\s = \o_{1+a} \g_b\s_A $ for $\s\in \Aut_K(\L_n)$.
\end{enumerate}
\end{abstract}


\section{Introduction}
Throughout, ring means an associative ring with $1$. Let $K$ be an
 arbitrary ring (not necessarily commutative). The {\em Grassmann
algebra} (the {\em exterior algebra}) $\L_n = \L_n (K)= K\lfloor
x_1, \ldots , x_n\rfloor$ is generated freely over $K$ by elements
$x_1, \ldots , x_n$ that satisfy the defining relations:
$$ x_1^2=\cdots = x_n^2=0 \;\; {\rm and}\;\; x_ix_j=-x_jx_i\;\;
{\rm for \; all} \;\; i\neq j.$$

{\bf What is the paper about? Motivation}.  Briefly, for the
Grassmann algebra $\L_n $ over a commutative ring $K$ we study in
detail the {\em Jacobian map} $$\CJ (\s ) := \det (\frac{\der \s
(x_i)}{\der x_j})$$ which is a `straightforward'
 generalization of the usual Jacobian map $\CJ (\s ) := \det
(\frac{\der \s (x_i)}{\der x_j})$ for a polynomial algebra  $P_n=
K[x_1, \ldots , x_n]$, $\s \in {\rm End}_{K-alg}(P_n)$. The
polynomial Jacobian map is  not yet a well-understood map, one of
the open questions about this map is the Jacobian conjecture (JC)
which claims {\em that $\CJ (\s ) =1$ implies} $\s \in
\Aut_K(P_n)$ (where $K$ is a field of characteristic zero).
Obviously, one can reformulate the Jacobian conjecture as the
question of whether the Jacobian monoid $ \S (P_n) : = \{ \s \in
{\rm End}_{K-alg}(P_n)\, | \, \CJ (\s ) =1\}$ is a {\em group}?
The analogous Jacobian monoid $\S = \S (\L_n)$ for the Grassmann
algebra $\L_n$ {\em is}, by a trivial reason,  a {\em group}, it
is a subgroup of the group $\Aut_K(\L_n)$ of automorphisms of the
Grassman algebra $\L_n$. It turns out that  properties of the
Jacobian map $\CJ$ are closely related to properties of the
Jacobian group $\S$ which should be treated as the  `kernel' of
the Jacobian map $\CJ$ despite the fact that $\CJ$ is {\em not} a
 homomorphism.

It turns out that the Jacobian group $\S$ is a large subgroup of
$\Aut_K(\L_n)$, so we start the paper considering the structure of
the group $\Aut_K(\L_n)$ and its subgroups. The Jacobian map and
the Jacobian group are not transparent objects to deal with.
Therefore, several (important) subgroups of $\Aut_K(\L_n)$ are
studied first. Some of them are given by explicit generators,
another are defined via certain `geometric' properties. That is
why we study these subgroups in detail. They are building blocks
in understanding the structure of the Jacobian map and the
Jacobian group. Let us describe main results of the paper.

In the Introduction,  $K$ is a {\em reduced commutative} ring with
$\frac{1}{2}\in K$,   $n\geq 2$ (though many results of the paper
are true under milder assumptions, see in the text),  $\L_n =
\L_n(K)= K\lfloor x_1, \ldots , x_n \rfloor$ be the  Grassmann
$K$-algebra and $\gm := (x_1, \ldots , x_n)$ be its augmentation
ideal. The algebra $\L_n$ is  endowed with the $\Z$-grading $\L_n=
\oplus_{i=0}^n\L_{n,i}$ and $\Z_2$-grading $\L_n= \Lnev \oplus
\Lnod$, and so each element $a\in \L_n$ is a unique sum $a=
a^{ev}+a^{od}$ where $a^{ev}\in \Lnev$ and $a^{od} \in \Lnod$. For
each $s\geq 2$, the algebra $\L_n$ is also a $\Z_s$-graded algebra
($\Z_s := \Z / s\Z $).

{\bf The structure of the group of automorphisms of the Grassmann
algebra and its subgroups}.  In Sections \ref{GOAG} and
\ref{PREXS}, we study the group $G:= \Aut_K(\L_n (K))$ of
$K$-algebra automorphisms of $\L_n$ and various its subgroups (and
their relations):
\begin{itemize}
\item $G_{gr}$, the subgroup of $G$ elements of which respect
$\Z$-grading, \item $G_{\Z_2-gr}$, the subgroup of $G$ elements of
which respect $\Z_2$-grading, \item $G_{\Z_s-gr}$, the subgroup of
$G$ elements of which respect $\Z_s$-grading,\item $U:= \{ \s \in
G\, | \, \s (x_i)= x_i+\cdots\; {\rm for \; all} \; i \}$ where
the three dots mean bigger terms with respect to the $\Z$-grading,
\item $G^{od}:= \{ \s \in G\, | \, \s (x_i) \in \L_{n,1}+\Lnod
\;\; {\rm for \; all}\;\; i\}$ and $G^{ev}:= \{ \s \in G\, | \, \s
(x_i) \in \L_{n,1}+\Lnev  \;\; {\rm for \; all}\;\; i\}$, \item
${\rm Inn} (\L_n):= \{ \o_u :x\mapsto uxu^{-1}\}$ and ${\rm
Out}(\L_n):=G/ {\rm Inn}(\L_n)$, the groups of inner and  outer
automorphisms, \item $\O := \{ \o_{1+a}\, | \, a\in \Lnod \}$,
\item For each odd number $s$ such that $1\leq s \leq n$, $\O
(s):= \{ \o_{1+a}\, | \, a\in \sum_{1\leq j\; {\rm is\, odd}}
\L_{n,js} \}$, \item $\G := \{ \g_b \, | \, \g_b(x_i) = x_i+b_i,
\; b_i \in \Lnod \cap \gm^3, i=1, \ldots , n\}$, $ b=(b_1, \ldots
, b_n)$, \item For each even number $s$ such that $3\leq s \leq
n$, $\G (s):= \{ \g_b \, | \, {\rm all} \; b_i\in \sum_{j\geq 1}
\L_{n, 1+js} \}$, \item $U^n:= \{ \tau_\l \, | \, \tau_\l (x_i) =
x_i+\l_ix_1\cdots x_n, \; \l = (\l_1, \ldots , \l_n)\in K^n\}
\simeq K^n$, $\tau_\l \lra \l$, \item $\GL_n(K)^{op}:= \{ \s_A\, |
\, \s_A(x_i)= \sum_{j=1}^n a_{ij}x_j, \; A=(a_{ij})\in
\GL_n(K)\}$, \item $\Phi := \{ \s : x_i\mapsto x_i(1+ a_i)\, | \,
a_i\in \Lnev \cap \gm^2, \, i=1, \ldots , n\}$.
\end{itemize}

If $K= \mathbb{C}$ the group $\G \GL_n(\mathbb{C})^{op}$ was
considered in \cite{Berezin67}. If $K= k$ is a field of
characteristic $\neq 2$ it was proved in \cite{Djokovic78} that
$G$ is a semidirect product ${\rm Inn}(\L_k(k))\rtimes
G_{\Z_2-gr}$. One can find a lot of information about the
Grassmann algebra (i.e. the exterior algebra) in
\cite{BourbakiAlgCh1-3}.

\begin{itemize}
\item (Lemma \ref{m28Sep06}.(5)) $\O$ {\em is an abelian group
canonically isomorphic to the additive group} $\Lnod/\Lnod\cap
Kx_1\cdots x_n$ via $\o_{1+a}\mapsto a$. \item (Lemma
\ref{i30Sep06}, Corollary \ref{G30Sep06}.(3)) ${\rm Inn}(\L_n)=
\O$ {\em and} ${\rm Out} (\L_n) \simeq G_{\Z_2-gr}$. \item
(Theorem \ref{29Sep06}) $U= \O \rtimes \G$. \item (Theorem
\ref{28Sep06}) $\O$ {\em is a maximal abelian subgroup of $U$ if
$n$ is even ($\O \supseteq U^n$); and $\O U^n=\O \times U^n$ is a
maximal abelian subgroup of $U$ if $n$ is odd} $(\O \cap U^n = \{
e\}$). \item (Theorem \ref{29Sep06}, Corollary \ref{G30Sep06},
Lemma \ref{e30Sep06})
\begin{enumerate}
\item $G= U\rtimes \GL_n(K)^{op} = (\O \rtimes \G ) \rtimes
\GL_n(K)^{op}$, \item $G= \O \rtimes G_{\Z_2-gr}$, {\em and} \item
$G= G^{ev} G^{od} = G^{od} G^{ev}$.
\end{enumerate}
\item (Lemma \ref{e30Sep06}.(1)) $G^{od} = G_{\Z_2-gr} = \G
\rtimes\GL_n(K)^{op}$. \item (Lemma \ref{4Oct06})  {\em Let} $s=2,
\ldots , n$. {\em Then}
$$
G_{\Z_s-gr} = \begin{cases}
\G (s) \rtimes \GL_n(K)^{op},& \text{if  $s$ is even},\\
\O (s) \rtimes \GL_n(K)^{op},& \text{if $s$ is odd}.\\
\end{cases}$$
\end{itemize}

{\bf The Jacobian matrix and an analogue of the Jacobian
Conjecture  for $\L_n$}. The even subalgebra $\Lnev$ of $\L_n$
belongs to the centre of the algebra $\L_n$. A $K$-linear map $\d
: \L_n \ra \L_n$ is called a {\em left skew derivation} if $\d
(a_ia_j) = \d (a_i) a_j +(-1)^ia_i\d (a_j)$ for all homogeneous
elements $a_i$ and $a_j$ of graded degree $i$ and $j$
respectively. Surprisingly, skew derivations rather than ordinary
derivations are more important in study of the Grassmann algebras.
This and the forthcoming paper \cite{skedgras} illustrate this
phenomenon.

The `partial derivatives' $\der_1:= \frac{\der}{\der x_1}, \ldots
,\der_n:= \frac{\der}{\der x_n}$ are  {\em left  skew
$K$-derivations} of $\L_n(K)$ ($\der_i(x_j)=\d_{ij}$, the
Kronecker delta; $\der_k (a_ia_j) = \der_k
(a_i)a_j+(-1)^ia_i\der_k(a_j)$). Let $\CE^{od}:={\rm
End}_{K-alg}(\L_n)^{od}:= \{ \s \in {\rm End}_{K-alg}(\L_n) \, |$
all $\s (x_i) \in \Lnod\}$. For each endomorphism $\s \in
\CE^{od}$, its {\em Jacobian matrix} $\frac{\der \s }{\der x}:=
(\frac{\der \s (x_i) }{\der x_j})\in M_n(\Lnev )$ has even (hence
central) entries, and so, the {\em Jacobian} of $\s$,
$$ \CJ (\s ) := \det (\frac{\der \s (x_i) }{\der x_j}),  $$
is well-defined element of the even subalgebra $\Lnev$. For $\s ,
\tau \in \CE^{od}$, the `chain rule' holds
$$\frac{\der (\s \tau )}{\der x}=\s (\frac{\der \tau }{\der x})
\cdot \frac{\der \s }{\der x},$$ which implies, $\CJ (\s \tau ) =
\CJ (\s )\cdot  \s (\CJ (\tau ))$, i.e. the Jacobian map is almost
a homomorphism of monoids (with zeros). It follows that the sets
$\S \subset \tS \subset \OO$ are monoids where
$$ \OO := \{ \s \in \CE^{od}\, | \, \CJ (\s ) \;\; {\rm is \; a \;
unit}\}, \;  \tS := \{ \s \in \CE^{od}\, | \, \CJ (\s )=1\}, \; \S
:= \{ \s \in \CE^{od}\cap \G \, | \, \CJ (\s )=1\}.$$

If $\s \in \CE^{od}\cap \Aut_K(\L_n)=G_{\Z_2-gr}= \G
\GL_n(K)^{op}$ (Lemma \ref{e30Sep06}.(1)) then $\s^{-1} \in
\CE^{od}\cap \Aut_K(\L_n)$ and
$$ 1= \CJ ({\rm id}_{\L_n})= \CJ (\s \s^{-1})= \CJ (\s ) \s (\CJ
(\s^{-1})),$$ and so $\CJ (\s )$ is a unit in $\L_n$.

An analogue of the Jacobian Conjecture  for the Grassmann algebra
$\L_n$, i.e. $\s \in \OO$ {\em implies $\s$ is an automorphism
(i.e. $\OO\subseteq G:=\Aut_K(\L_n)$)}, is trivially true.

\begin{itemize}
\item $\OO\subseteq G$.
\end{itemize}
{\it Proof}. Let $\s \in \OO$. Then $\s (x) = Ax+\cdots$ where
$x:= (x_1, \ldots , x_n)^t$, $A\in M_n(K)$, and the three dots
mean  higher terms with respect to the $\Z$-grading. Since
$$\CJ (\s ) \equiv \det (A)\mod \gm$$ and $\CJ (\s )$ is a unit,
the determinant must be a unit in $K$. Changing $\s$ for
$\s\s_{A^{-1}}$ where $\s_{A^{-1}}(x) = A^{-1} x$
($\s\s_{A^{-1}}(x) = A^{-1} Ax+\cdots = x+\cdots $) one can assume
that, for each $i=1, \ldots , n$, $\s (x_i) = x_i-a_i$ where
$a_i\in \gm^3\cap \Lnod$. Let us denote $\s (x_i)$ by $x_i'$, then
$x_i= x_i'+a_i(x)$. After repeating several times ($\leq n$ times)
these substitutions simultaneously  in the tail of each element
$x_i$, i.e. elements of degree $\geq 3$, it is easy to see that
$x_i= x_i'+b_i(x')$ for some element $b_i(x')\in \gm^3\cap \Lnod$.
This gives the {\em inverse} map for $\s$,  i.e. $\s \in G$. The
elements $b_i$ can be found even explicitly using the inversion
formula (Theorem \ref{i9Sep06}). $\Box$

The next two facts follow directly from the inclusion
$\OO\subseteq G$ and the formula $\CJ(\s^{-1})= \s^{-1} (\CJ (\s
)^{-1})$, for all $\s \in G$.

\begin{itemize}
\item $\S$, $\tS$, {\em and $ \OO$ are groups.} \item  $\OO
=G_{\Z_2-gr} = \G \rtimes \GL_n(K)^{op}$ {\em and} $\tS = \S
\rtimes \SL_n(K)^{op}$.
\end{itemize}

{\bf The Jacobian map and the Jacobian group}. The set $E_n:=
K^*+\sum_{m\geq 1} \L_{n,2m}$ is the group of units of the even
subalgebra $\Lnev := \oplus_{m\geq 0} \L_{n,2m}$ of $\L_n$ where
$K^*$ is the group of units of the ring $K$; and $E_n':=
1+\sum_{m\geq 1} \L_{n,2m}$ is the subgroup of $E_n$.  Due to the
equality $\CJ (\s \tau )= \CJ (\s ) \s (\CJ (\tau ))$, to study
the {\em Jacobian map}
$$ \CJ : \G \GL_n(K)^{op}\ra E_n, \;\; \s \mapsto \CJ (\s ),$$
is the same as to study its restriction to $\G$:
$$ \CJ : \G \ra E_n', \;\; \s \mapsto \CJ (\s ).$$
When we mention the Jacobian map it means as a rule this map. The
Jacobian group $\S = \{ \s \in \G \, | \, \CJ (\s ) =1\}$ is
trivial iff $n\leq 3$. So, we always assume that $n\geq 4$ in the
results on the Jacobian group $\S$ and its subgroups.

An algebraic group $A$ over $K$ is called {\em affine} if its
algebra of regular functions is a polynomial algebra $K[t_1,
\ldots , t_d]$ with coefficients in $K$ where $d:= \dim (A)$ is
called the {\em dimension} of $A$ (i.e. $A$ is an {\em affine
space}). If $K$ is a field then $\dim (A)$ is the usual dimension
of the algebraic group $A$ over the field $K$.

\begin{itemize}
\item (Theorem \ref{19Nov06}) {\em The Jacobian group $\S $ is an
affine group over $K$ of dimension}
$$\dim (\S )=\begin{cases}
(n-1)2^{n-1} -n^2+2& \text{{\em if $n$ is even}},\\
(n-1)2^{n-1} -n^2+1& \text{{\em if $n$ is odd}}.\\
\end{cases}
 $$
\item  {\em The coordinate functions on $\S$ are given explicitly
by (\ref{scooS}) and Corollary \ref{a8Oct06}.(5).} \end{itemize}

A subgroup of an algebraic group $A$ over $K$ is called a $1$-{\em
parameter subgroup} if it is isomorphic to the algebraic group
$(K,+)$. {\em A minimal set of generators} for an affine algebraic
group $A$ over $K$ is a set of $1$-parameter subgroups that
generate the group $A$ as an abstract group but each smaller
subset does not generate $A$.

\begin{itemize}
\item (Theorem \ref{14Nov06}) {\em A (minimal) set of generators
for $\S$ is given explicitly. } \item (Corollary \ref{n8Oct06})
{\em The Jacobian group $\S$ is not a normal subgroup of $\G$ iff}
$n\geq 5$. \item (Theorem \ref{8Oct06}) {\em The Jacobian map $\CJ
: \G \ra E_n'$, $ \s\mapsto \CJ (\s )$, is surjective if $n$ is
odd, and it is not surjective if $n$ is even but in this case its
image is a closed affine subvariety of $E_n'$ of codimension $1$
which is given by a single equation. }
\end{itemize}

{\bf The subgroups $\S'$ and $\S''$ of the Jacobian group $\S$}.
To prove the (above) results about the Jacobian group $\S$, we,
first, study in detail two of its subgroups:
$$\S':= \S \cap \Phi =\{ \s : x_i\mapsto x_i(1+a_i)\, | \, \CJ (\s ) =1,
\; a_i\in \Lnev\cap \gm^2, \; 1\leq i \leq n\}$$ and the subgroup
$\S''$ is generated by the explicit automorphisms of $\S$ (see
(\ref{ksibi})):
$$ \xi_{i,b_i}: x_i\mapsto x_i+b_i, \;\; x_j\mapsto x_j, \;\;
j\neq i,$$ where $b_i\in K\lfloor x_1, \ldots , \widehat{x}_i,
\ldots , x_n\rfloor_{\geq 3}^{od}$ and $i=1, \ldots , n$.

The importance of these subgroups is demonstrated by the following
two facts.

\begin{itemize}
\item (Corollary \ref{a8Oct06}.(1)) $\S = \S'\S''$. \item (Theorem
\ref{13Oct06}.(1)) $\G = \Phi \S''$.
\end{itemize}
Note that each element $x_i$ is a {\em normal} element of $\L_n$:
$x_i\L_n= \L_nx_i$. Therefore, the ideal $(x_i)$ of $\L_n$
generated by the element $x_i$ determines a coordinate
`hyperplane.' The groups $\S'$ and $\S''$ have the following
geometric interpretation: the group $\S'$ preserves the coordinate
`hyperplanes' and elements of the group $\S''$ can be seen as
`rotations.'

By the definition, the group $\S'$ is a closed subgroup of $\S$,
it is not a normal subgroup of $\S$ unless $n\leq 5$. It is not
obvious from the outset whether the subgroup $\S''$ is closed or
normal. In fact, it is.

\begin{itemize}
\item (Theorem  \ref{20Nov06}.(2)) $\S''$ {\em is the closed
normal subgroup of $\S$, $\S''$ is an affine group of dimension}
$$\dim (\S'') = \begin{cases}
(n-1)2^{n-1}-n^2+2-(n-3){n\choose 2}& \text{{\em if $n$ is even}},\\
(n-1)2^{n-1}-n^2+1-(n-3){n\choose 2}& \text{{\em if $n$ is odd}},\\
\end{cases}
$$
{\em and the factor group $\S/\S''\simeq \S'/\S'^5$ is an abelian
affine group of dimension} $\dim (\S/\S'') = n{n-1\choose 2}-
{n\choose 2}= (n-3){n\choose 2}$.
 \item (Corollary \ref{d25Oct06}) {\em The group $\S'$ is an affine
 group over $K$ of dimension}
$$
 \dim (\S') =\begin{cases}
(n-2)2^{n-2}-n+2& \text{{\em if $n$ is even}},\\
(n-2)2^{n-2}-n+1& \text{{\em if $n$ is odd}}.\\
\end{cases}
$$
\item (Lemma \ref{i14Nov06}) {\em The intersection $\S'\cap \S''$
is a closed subgroup of $\S$, it is an affine  group over $K$ of
dimension}
$$
 \dim (\S'\cap \S'') =\begin{cases}
(n-2)2^{n-2}-n+2-(n-3){n\choose 2}& \text{{\em if $n$ is even}},\\
(n-2)2^{n-2}-n+1-(n-3){n\choose 2}& \text{{\em if $n$ is odd}}.\\
\end{cases}
$$
\item{\em The coordinates on $\S'$ and $\S''$ are given explicitly
by (\ref{scooS}) and (\ref{1scooS}).}
\end{itemize}
To find coordinates for the groups $\S$, $\S'$, and  $\S''$
explicitly, we introduce avoidance functions and a series of
subgroups $ \{ \Phi'^{2s+1}\}$ , $s=1, 2, \ldots ,
[\frac{n-1}{2}]$, of $\Phi$ that are given explicitly (see Section
\ref{GSID}). They are too technical to explain in the
introduction.

\begin{itemize}
\item (Theorem  \ref{g23Oct06}) {\em This theorem is a key result
in finding coordinates for the groups $\S$, $\S'$, $\S''$, etc.}
\end{itemize}

{\bf The Jacobian ascents $\G_{2s}$}. In order to study the image
of the Jacobian map $\CJ : \G \ra E_n'$, $\s \mapsto \CJ (\s )$,
certain overgroups of the Jacobian group $\S$ are introduced. They
are called the {\em Jacobian ascents}. The problem of finding the
image ${\rm im }(\CJ )$ is equal to the problem of finding
generators for these groups. Let us give some details. The
Grassmann algebra $\L_n$ has the $\gm$-adic filtration $\{
\gm^i\}$. Therefore, the group $E_n'$ has the induced $\gm$-adic
filtration:
$$ E_n'=E_{n,2}'\supset E_{n,4}'\supset \cdots \supset
E_{n,2m}'\supset \cdots \supset E_{n,2[\frac{n}{2}]}'\supset
E_{n,2[\frac{n}{2}]+2}'=\{ 1\},$$ where $E_{n,2m}':= E_n'\cap
(1+\gm^{2m})$. Correspondingly, the group $\G$ has the {\em
Jacobian filtration}:
$$ \G=\G_2\supseteq \G_4\supseteq \cdots \supseteq
\G_{2m}\supseteq \cdots \supseteq \G_{2[\frac{n}{2}]}\supseteq
\G_{2[\frac{n}{2}]+2}=\S ,$$ where $\G_{2m}:= \G_{n,2m}:=
\CJ^{-1}(E_{n,2m}')= \{ \s \in \G \, | \, \CJ (\s ) \in
E_{n,2m}'\}$. It follows from the equality $\CJ (\s \tau ) = \CJ
(\s ) \s (\CJ (\tau ))$ that all $\G_{2m}$ are subgroups of $\G$,
they are called, the {\em Jacobian ascents} of the Jacobian group
$\S$.

The Jacobian ascents are {\em distinct} groups with a {\em single}
exception when two groups coincide. This is  a subtle fact, it
explains (partly) why formulae for various dimensions  differ by 1
in odd and even cases.

\begin{itemize}
\item (Corollary \ref{a20Nov06})  {\em Let $K$ be a commutative
ring and} $n\geq 4$.
\begin{enumerate}
\item {\em If $n$ is an odd number then the Jacobian ascents}
$$
\G = \G_2\supset \G_4\supset \cdots \supset \G_{2s} \supset \cdots
\supset \G_{2[\frac{n}{2}]}\supset \G_{2[\frac{n}{2}]+2}=\S$$ {\em
are distinct groups.}

\item {\em If $n$ is an even  number then the Jacobian ascents}
$$ \G = \G_2\supset \G_4\supset \cdots \supset \G_{2s} \supset
\cdots \supset\G_{2[\frac{n}{2}]-2}\supset\G_{2[\frac{n}{2}]}
=\G_{2[\frac{n}{2}]+2}=\S$$ {\em are distinct groups except the
last two groups, i.e.} $\G_{2[\frac{n}{2}]}
=\G_{2[\frac{n}{2}]+2}$.
\end{enumerate}
\end{itemize}

The subgroups $\{ \G^{2s+1}\}$ of $\G$ are given explicitly,
$$\G^{2s+1}:= \{ \s : x_i\mapsto x_i+a_i\, | \, a_i\in \Lnod \cap
\gm^{2s+1}, \; 1\leq i \leq n \},\;\; s\geq 1,  $$
 they
have clear structure. The next result explains that the Jacobian
ascents $\{\G_{2s}\}$ have clear structure too, $\G_{2s}=
\G^{2s+1} \S$, and so  the structure of the Jacobian ascents is
completely determined by the structure of the Jacobian group $\S$.

\begin{itemize}
\item (Theorem \ref{G20Nov06}) {\em Let $K$ be a commutative ring
and $n\geq 4$. Then}
\begin{enumerate}
\item $\G_{2s}= \G^{2s+1}\S = \Phi^{2s+1}\S = \Phi'^{2s+1} \S$
{\em for
 each }  $s=1,2, \ldots , [\frac{n-1}{2}]$.  \item {\em If $n$ is an even
number then $\G_n= \S $, i.e. } $\G_n = \G_{n+2} =\S $.
\end{enumerate}
\end{itemize}

The next theorem introduces an isomorphic affine structure on the
algebraic group $\G$.

\begin{itemize}
\item (Theorem \ref{27Nov06}) {\em Let $K$ be a commutative ring,
$n\geq 4$, and $s=1, \ldots , [\frac{n-1}{2}]$. Then each
automorphism $\s \in \G$ is a unique product $\s =
\phi'_{a(2)}\phi'_{a(4)}\cdots \phi'_{a(2[\frac{n-1}{2}])}\g$ for
unique elements $a(2s)\in \L_{n,2s}$ and $\g \in
\G_{2[\frac{n-1}{2}]+2}=\S$ (by (\ref{G2n2})). Moreover,}
\begin{eqnarray*}
 a(2)&\equiv & \CJ (\s )-1 \mod E_{n,4}', \\
  a(2t)&\equiv & \CJ ({\phi'}_{a(2t-2)}^{-1} \cdots {\phi'}_{a(2)}^{-1}\s )-1
  \mod E_{n,2t+2}', \; t=2, \ldots , [\frac{n-1}{2}], \\
  \g &= & (\phi'_{a(2)}\phi'_{a(4)}\cdots
\phi'_{a(2[\frac{n-1}{2}])})^{-1} \s.
\end{eqnarray*}
\end{itemize}

The automorphisms $\phi'_{a(2s)}$ are give explicitly (see Section
\ref{IJMJA} for details), they are too technical to explain here.
The theorem above is a key result in proving that various quotient
spaces, like $\G_{2s}/\G_{2t}$ $(s<t)$, are {\em affine}, and in
finding their dimensions. An algebraic variety $V$ over $K$ is
called {\em affine} (i.e. an {\em affine space} over $K$) if its
algebra of regular functions $\OO (V)$ is a polynomial algebra $K[
v_1, \ldots , v_d]$ over $K$ where $d:= \dim (V)$ is called the
{\em dimension} of $V$ over $K$. In this paper, all  algebraic
groups and varieties will turn out to be affine (i.e. affine
spaces), and so the word `affine' is used only in this sense. This
fact strengthen relations between the Grassmann algebras and
polynomial algebras even more.

\begin{itemize}
\item (Corollary \ref{a28Nov06}) {\em Let $K$ be a commutative
ring, $n\geq 4$. Then all the Jacobian ascents are affine groups
over $K$ and  closed subgroups of $\G$, and}
$$\dim (\G_{2s})= \dim (\S ) +
\sum_{i=s}^{[\frac{n-1}{2}]}{n\choose 2i}, \;\; s=1, \ldots
,[\frac{n-1}{2}].$$ \item (Corollary \ref{y27Nov06}) {\em Let $K$
be a commutative ring and $n\geq 4$. Then the quotient space $\G /
\S : = \{ \s \S \, | \, \s \in \G \} $ is an affine variety.}
\begin{enumerate}
\item {\em If $n$ is odd then the Jacobian map $\G / \S \ra E_n'$,
$\s\S\mapsto \CJ (\s )$, is an isomorphism of the affine varieties
over $K$, and} $\dim (\G / \S ) = 2^{n-1} -1$. \item {\em If $n$
is even then the Jacobian map $\G / \S \ra E_n'/E_{n,n}'$,
$\s\S\mapsto \CJ (\s )E_{n,n}'$, is an isomorphism of the affine
varieties over $K$ (where $E_{n,n}'= 1+Kx_1\cdots x_n$), and}
$\dim (\G / \S ) = 2^{n-1} -2$.
\end{enumerate}
\end{itemize}

\begin{itemize}
\item (Theorem \ref{8Oct06}) {\em  Let $K$ be a commutative ring,
$n\geq 4$,  $\CJ :\G \ra E_n'$, $\s \mapsto \CJ (\s )$, be the
Jacobian map, and $s=1,2, \ldots , [\frac{n-1}{2}]$. Then,}
\begin{enumerate}
\item {\em for an odd number $n$, the Jacobian map $\CJ$ is
surjective, and}
 \item {\em for an even  number $n$, the Jacobian map $\CJ$ is not surjective. In more detail,
 the image ${\rm im } (\CJ )$ is a closed algebraic variety of
$E_n'$ of codimension 1.}
\end{enumerate}
\end{itemize}

{\bf The unique presentation $\s = \o_{1+a} \g_b\s_A$ for $\s \in
{\rm Aut}_K(\L_n)$}.  Each automorphism $\s \in G= \O \G
\GL_n(K)^{op}$ is a {\em unique} product (Theorem \ref{29Sep06})
$$\s = \o_{1+a} \g_b\s_A$$
where $\o_{1+a} \in \O$ ($a\in \L_n'^{od}$), $ \g_b\in \G$, and $
\s_A\in \GL_n(K)^{op}$ where $\L_n'^{od}:= \oplus_i \L_{n,i}$ and
$i$ runs through {\em odd} natural numbers such $1\leq i \leq
n-1$. The next theorem determines explicitly the elements $a$,
$b$, and $A$ via the vector-column $\s (x) := (\s (x_1) , \ldots ,
\s (x_n))^t$ (for, only one needs to know explicitly the inverse
$\g_b^{-1}$ for each $\g_b\in \G$ which is given by the inversion
formula below, Theorem \ref{i9Sep06}).

\begin{itemize}
\item  (Theorem \ref{M30Sep06})  {\em  Each element $\s \in G$ is
a unique product $\s = \o_{1+a} \g_b\s_A$ (Theorem
\ref{29Sep06}.(3)) where $a\in \L_n'^{od}$ and}
\begin{enumerate}
\item $\s (x) = Ax +\cdots $ {\em (i.e. $\s (x) \equiv Ax\mod
\gm$) for some} $A\in \GL_n(K)$, \item $b= A^{-1} \s(x)^{od}-x$,
{\em and } \item $a= -\frac{1}{2}\g_b(\sum_{i=1}^{n-1} x_1\cdots
x_i\der_i\cdots \der_1\der_{i+1} (a_{i+1}') +\der_1(a_1'))$ {\em
where $a_i':= (A^{-1} \g_b^{-1} (\s (x)^{ev}))_i$, the $i$'th
component of the column-vector} $A^{-1} \g_b^{-1} (\s (x)^{ev})$,
\end{enumerate}
{\em where $\der_1:= \frac{\der}{\der x_1}, \ldots ,\der_n:=
\frac{\der}{\der x_n}$ are left skew $K$-derivations of $\L_n(K)$
($\der_i(x_j)=\d_{ij}$, the Kronecker delta).}
\end{itemize}

{\bf The inversion formula for automorphisms}. The formula
(\ref{mul1}) for multiplication of elements of $G$ shows that the
most non-trivial (difficult) part of the group $G$ is the group
$\G$. Elements of the group $\G$ should be seen as $n$-tuples of
noncommutative polynomials in anti-commuting variables $x_1,
\ldots , x_n$, and the multiplication of two $n$-tuples is the
composition of functions. The group $\G$ (and $\G \GL_n(K)^{op}$)
is the part of the group $G$ that `behaves' in a similar fashion
as polynomial automorphisms. This very observation we will explore
in the paper. The analogy between the group $\G$ and the group of
polynomial automorphisms is far more reaching than one may expect.

\begin{itemize}
\item  (Theorem \ref{i9Sep06}) {\rm (The Inversion Formula)} {\em
Let $K$ be a commutative ring,
 $\s \in \G \GL_n(K)^{op}$
 and $a\in \L_n(K)$. Then} $\s^{-1}(a)=\sum_{\alpha \in \CB_n} \l_\alpha
 x^\alpha$ {\em where}
\begin{eqnarray*}
 \l_\alpha &:=& (1-\s (x_n)\der_n')
(1-\s (x_{n-1})\der_{n-1}') \cdots (1-\s (x_1)\der_1')
\der'^\alpha (a)\in K, \\
 \der'^{\alpha} &:= &\der_n'^{\alpha_n}
\der_{n-1}'^{\alpha_{n-1}} \cdots \der_1'^{\alpha_1} , \\
 &
\end{eqnarray*}
$$\der_i'  := \frac{1}{\det (\frac{\der \s (x_\nu)}{\der x_\mu})} \, \det
 \begin{pmatrix}
  \frac{\der \s (x_1)}{\der x_1} & \cdots & \frac{\der \s (x_1)}{\der x_m} \\
  \vdots & \vdots & \vdots \\
\frac{\der }{\der x_1} & \cdots & \frac{\der }{\der x_m}\\
 \vdots & \vdots & \vdots \\
\frac{\der \s (x_n)}{\der x_1} & \cdots & \frac{\der \s (x_n)}{\der x_m} \\
\end{pmatrix}, \;\;\; i=1, \ldots , n,$$
{\em where $\der_1:= \frac{\der}{\der x_1}, \ldots ,\der_n:=
\frac{\der}{\der x_n}$ are left skew $K$-derivations of
$\L_n(K)$.}
\end{itemize}

 Then, for any automorphism $\s = \o_{1+a} \g_b\s_A\in G$, one can
write explicitly the formula for the inverse $\s^{-1}$,
(\ref{invabA}). Note that $\o_{1+a}^{-1} = \o_{1-a}$ and $
\s_A^{-1}= \s_{A^{-1}}$. So, $\g_b^{-1}$ is the most difficult
part of the inverse map $\s^{-1}$. The formula for $\g_{b}^{-1}$
is written via {\em skew differential operators} (i.e. linear
combinations of products of powers of skew derivations). That is
why we start the paper with various properties of skew
derivations. Detailed study of skew derivations is continued in
\cite{skedgras}.

{\bf Analogues of the Poincar\'{e} Lemma}.  The crucial step in
finding the $b$ (in $\s = \o_{1+a} \g_b\s_A$) is an analogue of
the Poincar\'{e} Lemma for $\L_n$ (Theorem \ref{s14Sep06}) where
the solutions (as well as necessary and sufficient conditions for
existence of solutions) are given explicitly  for the following
system of equations in $\L_n$ where $a\in \L_n$ is unknown, and
$u_i\in \L_n$:
$$\begin{cases}
x_1a=u_1, \\
x_2a=u_2, \\
\;\;\;\; \;\;\;\vdots \\
x_na=u_n.
\end{cases}
$$

\begin{itemize}
\item (Theorem \ref{s14Sep06}) {\em Let $K$ be an arbitrary ring.
The system above has a solution iff $(i)$ $u_1\in (x_1),\ldots ,
u_n\in (x_n)$, and $(ii)$ $x_iu_j= - x_ju_i$ for all $i\neq j$.
Then
$$  a= x_1\cdots x_na_n+\sum_{i=1}^{n-1} x_1\cdots x_i\der_i\cdots
\der_1\der_{i+1}(u_{i+1})+\der_1(u_1), \;\; a_n \in K,
$$
are all the solutions.}
\end{itemize}
Note that the left multiplication on $x_i$ is, up to the scalar
$\frac{1}{2}$, a skew derivation in $\L_n$: $x_i(a_ja_k)=
\frac{1}{2}((x_ia_j)a_k+(-1)^ja_j(x_ia_k))$ for all homogeneous
elements $a_j$ and $a_k$ of $\Z$-graded degree $j$ and $k$
respectively.

Another version of the Poincar\'{e} Lemma for $\L_n$ is Theorem
\ref{P14Sep06} where the solutions are given explicitly for the
system (of first order partial skew differential operators):
$$\begin{cases}
\der_1(a)=u_1, \\
\der_2(a)=u_2, \\
\;\;\;\; \;\;\;\vdots \\
\der_n(a)=u_n,
\end{cases}
$$
where $ \der_1:= \frac{\der}{\der x_1}, \ldots ,\der_n:=
\frac{\der}{\der x_n}$ are the  left partial skew $K$-derivatives
of $\L_n$.

\begin{itemize}
\item (Theorem \ref{P14Sep06}) {\em Let $K$ be an arbitrary ring.
The system above has a solution iff $(i)$ $u_i\in K\langle
x_1,\ldots , x_{i-1}, x_{i+1}, \ldots  x_n\rangle$ for all $i$;
and $(ii)$ $\der_i(u_j)= - \der_j(u_i)$ for all $i\neq j$. Then
$$ a=\l +\sum_{0\neq \alpha \in \CB_n} (1-x_n\der_n)
(1-x_{n-1}\der_{n-1}) \cdots (1-x_1\der_1) (u_\alpha ) x^\alpha,
\;\; \l \in K,
$$
are all the solutions where for} $\alpha = \{ i_1<\cdots < i_k\}$,
$u_\alpha := \der_{i_k}\der_{i_{k-1}}\cdots \der_{i_2}(u_{i_1})$.
\end{itemize}

{\bf Minimal set of generators for the group $\G$ and some of its
subgroups}.  For each $i=1, \ldots , n$; $\l \in K$; and $ j
<k<l$, let us consider the automorphism $\s_{i, \l x_jx_kx_l}\in
\G$: $x_i\mapsto x_i+\l x_jx_kx_l$, $x_m\mapsto x_m$, for all
$m\neq i$. Then
$$\s_{i, \l x_jx_kx_l}\s_{i, \mu x_jx_kx_l}= \s_{i, (\l +\mu
)x_jx_kx_l}, \;\;\;\; \s_{i, \l x_jx_kx_l}^{-1} = \s_{i, -\l
x_jx_kx_l}^{-1}.$$ So, the group $\{ \s_{i, \l x_jx_kx_l} \, | \,
\l \in K\}$ is  isomorphic to the algebraic group $(K, +)$ via  $
\s_{i, \l x_jx_kx_l}\mapsto \l$.

\begin{itemize}
\item (Theorem \ref{7Oct06}.(1)) {\em The group $\G $ is generated
by all the automorphisms $\s_{i, \l x_jx_kx_l}$, i.e. $\G =
\langle \s_{i, \l x_jx_kx_l}\, | \, i=1, \ldots , n; \l \in K; j
<k<l\rangle$. The  subgroups $\{ \s_{i, \l x_jx_kx_l}\}_{\l \in
K}$  form a minimal set of generators for $\G$.} \item (Theorem
\ref{7Oct06}.(3)) {\em The group $U$ is generated by all the
automorphisms $\s_{i, \l x_jx_kx_l}$ and all the automorphisms
$\o_{1+\l x_i}$, i.e. $U = \langle \s_{i, \l x_jx_kx_l}, \o_{1+\l
x_i}\, | \, , i=1, \ldots , n; \l \in K; j <k<l\rangle$. The
subgroups $\{ \s_{i, \l x_jx_kx_l}\}_{\l \in K}$, $\{ \o_{1+\l
x_i}\}_{\l \in K}$     form a minimal set of generators for $U$.}
\item (Corollary \ref{p15Oct06}) {\em The group $\Phi$ is
generated by all the automorphisms $\s_{i, \l x_ix_kx_l}$, i.e.
$\Phi = \langle \s_{i, \l x_ix_kx_l}\, | \, i=1, \ldots , n; \l
\in K;
 k<l; i\not\in \{ k,l\}\rangle$. The subgroups $\{ \s_{i,
\l x_ix_kx_l}\}_{\l \in  K}$  form a minimal set of generators for
$\Phi$.}
\end{itemize}


\section{The group of automorphisms of the Grassmann
ring}\label{GOAG}

For reader's convenience, at the beginning of this section some
elementary results on Grassmann rings and their left skew
derivations are collected. Later in the paper they are used in
proofs of  many explicit formulae. In the second part of this
section, all the results on the group of automorphisms of the
Grassmann ring and its subgroups (from the Introduction) are
proved.

{\bf The Grassmann algebra and its gradings}.  Let $K$ be an {\em
arbitrary} ring (not necessarily commutative). The {\em Grassmann
algebra} (the {\em exterior algebra}) $\L_n = \L_n (K)= K\lfloor
x_1, \ldots , x_n\rfloor$ is generated freely over $K$ by elements
$x_1, \ldots , x_n$ that satisfy the defining relations:
$$ x_1^2=\cdots = x_n^2=0 \;\; {\rm and}\;\; x_ix_j=-x_jx_i\;\;
{\rm for \; all} \;\; i\neq j.$$ Let $\CB_n$ be the set of all
subsets of the set of indices $\{ 1, \ldots , n\}$. We may
identify the set $\CB_n$ with the direct product $\{ 0,1\}^n$ of
$n$ copies of the two-element set $\{ 0, 1\}$ by the rule $\{ i_1,
\ldots , i_k\} \mapsto (0, \ldots, 1, \ldots , 1, \ldots , 0)$
where $1$'s are on $i_1, \ldots , i_k$ places and $0$'s elsewhere.
So, the set $\{ 0, 1\}^n$ is the set of all the characteristic
functions on the set $\{ 1, \ldots , n\}$.
$$ \L_n = \bigoplus_{\alpha \in \CB_n} Kx^\alpha = \bigoplus_{\alpha \in \CB_n} x^\alpha K,
 \;\; x^\alpha :=x_1^{\alpha_1} \cdots  x_n^{\alpha_n}, $$
where $\alpha = (\alpha_1, \ldots , \alpha_n)\in \{ 0, 1\}^n =
\CB_n$.  Note that the order in the product $x^\alpha$ is fixed.
So, $\L_n$ is a free left and right $K$-module of rank $2^n$. The
ring $\L_n(K)$ is commutative iff $K$ is commutative and either
$n=1$ or $-1=1$. Note that $(x_i):= x_i\L_n = \L_nx_i$ is an ideal
of $\L_n$. Each element $a\in \L_n$ is a unique sum $a= \sum
a_\alpha x^\alpha$, $a_\alpha \in K$. One can view each element
$a$ of $\L_n$ as a `function' $a= a(x_1, \ldots , x_n)$ in the
non-commutative variables $x_i$. The $K$-algebra epimorphism
\begin{eqnarray*}
 \L_n &\ra & \L_n / (x_{i_1}, \ldots , x_{i_k})= K\lfloor x_1, \ldots
, \widehat{x_{i_1}}, \ldots , \widehat{x_{i_k}}, \ldots ,
x_n\rfloor, \\
a&\mapsto & a|_{x_{i_1}=0, \ldots , x_{i_k}=0} :=a+ (x_{i_1},
\ldots , x_{i_k}),
\end{eqnarray*}
 may be seen as the operation of
taking value of the function $a(x_1, \ldots , x_n)$ at the point
$x_{i_1} = \cdots = x_{i_k}=0$ where here and later the hat over a
symbol means that it is missed.

For each $\alpha \in \CB_n$, let $ | \alpha | := \alpha_1 +\cdots
+ \alpha_n$. The ring $\L_n=\oplus_{i=0}^n \L_{n,i}$ is a
$\Z$-{\em graded} ring ($\L_{n,i}\L_{n,j}\subseteq \L_{n,i+j}$ for
all $i,j$) where $\L_{n,i}:= \oplus_{|\alpha | =i}Kx^\alpha$. The
ideal $\gm := \oplus_{i\geq 1} \L_{n,i}$ of $\L_n$  is called the
{\em augmentation} ideal. Clearly, $K\simeq \L_n/ \gm$, $\gm^n=
Kx_1\cdots x_n$ and $\gm^{n+1}=0$. We say that an element $\alpha$
of $\CB_n$ is {\em even} (resp. {\em odd}) if the set $\alpha$
contains even (resp. odd) number of elements. By definition, the
empty set is even. Let $\Z_2:= \Z / 2\Z = \{\0, \1 \}$. The ring
$\L_n = \Lnz \oplus \Ln1$ is a $\Z_2$-{\em graded} ring where
$\Lnz := \Lnev :=\oplus_{\alpha \; {\rm is \; even}}Kx^\alpha$ is
the subring of even elements of $\L_n$ and $\Ln1 := \Lnod
:=\oplus_{\alpha \; {\rm is \; odd}}Kx^\alpha$ is the
$\Lnev$-module of odd elements of $\L_n$. The ring $\L_n$ has the
$\gm$-{\em adic} filtration $\{ \gm^i\}_{i\geq 0}$. The even
subring $\Lnev$ has the induced $\gm$-{\em adic} filtration $\{
\L_{n, \geq i}^{ev} := \Lnev \cap \gm^i \}$. The $\Lnev$-module
$\Lnod$ has the induced $\gm$-{\em adic} filtration $\{ \L_{n,
\geq i}^{od} := \Lnod \cap \gm^i \}$.

 The
$K$-linear map $ a\mapsto \oa$ from $\L_n$ to itself which is
given by the rule
$$\oa :=
\begin{cases}
a,& \text{if $a\in \Lnz$},\\
-a,& \text{if $a\in \Ln1$},
\end{cases}$$
is a ring automorphism such that  $\overline{\oa}= a$ for all
$a\in \L_n$. For all $a\in \L_n$ and $i=1, \ldots , n$,
\begin{equation}\label{xiaib}
x_ia= \oa x_i \;\; {\rm and}\;\; ax_i= x_i\oa .
\end{equation}
So, each element $x_i$ of $\L_n$ is a {\em normal} element, i.e.
the two-sided ideal $(x_i)$ generated by the element $x_i$
coincides with both left and right ideals generated by $x_i$:
$(x_i) = \L_n x_i= x_i\L_n$.

For an arbitrary $\Z$-graded ring $A= \oplus_{i\in \Z} A_i$,  an
additive map $\d :A\ra A$ is called a {\em left skew derivation}
if 
\begin{equation}\label{dsd}
\d (a_ia_j) = \d (a_i) a_j+(-1)^ia_i\d (a_j)\;\; {\rm for \;
all}\;\; a_i\in A_i, \; a_j\in A_j.
\end{equation}
In this paper, a skew derivation means a {\em left} skew
derivation. Clearly, $1\in \ker (\d )$ ($ \d (1) = \d (1\cdot 1) =
2\d (1)$ and so $\d (1)=0$). The restriction of the left skew
derivation $\d$ to the even subring $A^{ev}:= \oplus_{i\in 2\Z}
A_i$ of $A$ is an ordinary derivation. Recall that an additive
subgroup $B$ of $A$ is called a homogeneous subgroup if $B=
\oplus_{i\in \Z } B\cap A_i$. If the kernel $\ker (\d )$ of $\d$
is a  homogeneous additive subgroup of $A$ then $\ker (\d )$ is a
 subring of $A$, by (\ref{dsd}).

{\it Definition}. For the  ring $\L_n (K)$, consider the set of
{\em left} skew $K$-derivations:
$$ \der_1:= \frac{\der}{\der x_1}, \ldots , \der_n:=
\frac{\der}{\der x_n}$$ given by the rule $\der_i (x_j)= \d_{ij}$,
the Kronecker delta. Informally, these skew $K$-derivations will
be called (left) {\em partial skew derivatives}.

{\it Example}. $\der_i (x_1 \cdots x_i \cdots x_k)= (-1)^{i-1}
x_1\cdots x_{i-1} x_{i+1} \cdots x_k$.

{\bf The Taylor formula and its generalization}.  In this paper,
$\der_1, \ldots ,\der_n$ mean {\em left} partial skew  derivatives
(if it is not stated otherwise). Note that
$$\der_1^2= \cdots =\der_n^2=0 \;\; {\rm and} \;\;
 \der_i\der_j= - \der_j\der_i\;\; {\rm for\; all}\;\; i\neq j,$$
and $K_i:= \ker (\der_i)= K\lfloor x_1, \ldots , \widehat{x_i},
\ldots , x_n\rfloor$.

\begin{lemma}\label{T9Sep06}
(The Taylor Formula) For each $a= \sum_{\alpha \in \CB_n} a_\alpha
x^\alpha \in \L_n(K)$,
 $$a= \sum_{\alpha \in \CB_n}
\der^\alpha (a) (0)\,  x^\alpha$$ where $\der^\alpha :=
\der_n^{\alpha_n} \der_{n-1}^{\alpha_{n-1}}\cdots
\der_1^{\alpha_1}$, in the \underline{reverse} order here and
\underline{everywhere}.
\end{lemma}

{\it Proof}. It is obvious since $a_\alpha \equiv \der^\alpha (a)
\mod \gm $.  $\Box $

The operation of taking value at $0$ in the Taylor Formula is
rather `annoying'. Later, we will give an `improved' (more
economical) version of the Taylor Formula without the operation of
taking value at $0$ (Theorem \ref{s9Sep06}).

\begin{lemma}\label{p9Sep06}
Recall that $K_i:= \ker (\der_i)= K\lfloor x_1, \ldots ,
\widehat{x_i}, \ldots , x_n\rfloor$ and $\L_n(K)= K_i\oplus
x_iK_i$. Then
\begin{enumerate}
\item for each $i=1, \ldots , n$, the map $\phi_i :=
1-x_i\der_i:\L_n\ra \L_n$ is the  projection onto $K_i$. \item The
composition of the maps $\phi := \phi_n\phi_{n-1}\cdots \phi_1:
\L_n\ra \L_n$ is the projection onto $K$ in $\L_n= K\oplus \gm$.
\item $\phi = (1-x_n\der_n)(1-x_{n-1}\der_{n-1})\cdots
(1-x_1\der_1) = \sum_{\alpha \in \CB_n } (-1)^{|\alpha | }x^\alpha
\der^\alpha$ where  $x^\alpha := x_1^{\alpha_1}\cdots
x_n^{\alpha_n}$ and $\der^\alpha := \der_n^{\alpha_n}
\der_{n-1}^{\alpha_{n-1}}\cdots \der_1^{\alpha_1}$, in the reverse
order.
\end{enumerate}
\end{lemma}

{\it Proof}. $1$. By the very definition, $\phi_i$ is a right
$K_i$-module endomorphism of $\L_n$ with $\phi_i(x_i)= x_i-x_i=0$,
hence $\phi_i$ is the projection onto $K_i$ since $\L_n= K_i\oplus
x_iK_i$.

$2$. This follows from statement $1$ and the decomposition $\L_n =
\oplus_{\alpha\in \CB_n}Kx^\alpha$.

$3$.\begin{eqnarray*}
 \phi &=& \sum_{i_1>\cdots >i_k} (-1)^k
x_{i_1}\der_{i_1}x_{i_2}\der_{i_2}\cdots
x_{i_k}\der_{i_k}=\sum_{i_1>\cdots >i_k} (-1)^k(-1)^{1+2+\cdots
+k-1}  x_{i_1}\cdots x_{i_k}\der_{i_1}\cdots
\der_{i_k}\\
&= &\sum_{i_1>\cdots >i_k} (-1)^k x_{i_1}\cdots
x_{i_k}\der_{i_k}\cdots \der_{i_1}=
\sum_{\alpha \in \CB_n } (-1)^{|\alpha | }x^\alpha \der^\alpha. \;\; \Box \\
\end{eqnarray*}

The next theorem gives a kind of the Taylor Formula which is more
economical then the original Taylor Formula (no evaluation at
$0$).

\begin{theorem}\label{s9Sep06}
For each $a= \sum_{\alpha \in \CB_n} a_\alpha x^\alpha \in
\L_n(K)$,
\begin{enumerate}
\item $a= \sum_{\alpha \in \CB_n} \phi (\der^\alpha (a))
x^\alpha$.\item $a= \sum_{\alpha \in \CB_n}\,  (\sum_{\beta\in
\CB_n} (-1)^{|\beta |} x^\beta \der^\beta \der^\alpha (a)) \,
x^\alpha$.
\end{enumerate}
\end{theorem}

{\it Proof}. $1$. Note that $\der^\alpha (a) \equiv a_\alpha \mod
\gm$, hence $a_\alpha = \phi (\der^\alpha (a))$, and so the
result.

$2$. This follows from statement $1$ and Lemma \ref{p9Sep06}.(3).
 $\Box $

The Grassmann ring $ \L_n (K)=\oplus_{\alpha \in \CB_n}x^\alpha K$
is a free right $K$-module of rank $2^n$. The ring $\End_K(\L_n)$
of right $K$-module endomorphisms of $\L_n$ is canonically
isomorphic to the ring $M_{2^n}(K)$ of all $2^n\times 2^n$
matrices  with entries from $K$ by taking the matrix of map with
respect to the canonical basis $\{ x^\alpha , \alpha \in \CB_n\}$
of $\L_n$  as the  right $K$-module. We often identify these two
rings. Let $ \{ E_{\alpha \beta}\, | \,  \alpha , \beta \in
\CB_n\}$ be the matrix units of $M_{2^n}(K)= \oplus_{\alpha ,
\beta \in \CB_n}KE_{\alpha \beta}$ (i.e. $E_{\alpha \beta}(x^\g )=
\d_{\beta \g}x^\alpha$). One can identify the ring $\L_n$ with its
isomorphic copy in $\End_K(\L_n)$ via the ring monomorphism $
a\mapsto (x\mapsto ax)$.
\begin{theorem}\label{E12Sep06}
Recall that $\phi := (1-x_n\der_n)(1-x_{n-1}\der_{n-1})\cdots
(1-x_1\der_1)$. Then
\begin{enumerate}
\item for each $\alpha , \beta \in \CB_n$, $E_{\alpha \beta} =
x^\alpha \phi \der^\beta$.\item $\End_K(\L_n)=\oplus_{\alpha \in
\CB_n}\L_n \der^\alpha = \oplus_{\alpha \in \CB_n} \der^\alpha
\L_n$.
\end{enumerate}
\end{theorem}

{\it Proof}. $1$. Each element $a=\sum_{\g \in \CB_n} x^\g a_\g\in
\L_n$, $a_\g \in K$, can also be written as $a= \sum_{\g \in
\CB_n} a_\g x^\g$. As we have seen in the proof of Theorem
\ref{s9Sep06}, $a_\g = \phi \der^\g (a)$, hence $x^\alpha \phi
\der^\beta (a) = x^\alpha a_\beta$ which means that $E_{\alpha
\beta} = x^\alpha \phi \der^\beta$.

$2$. Since $\der_1, \ldots , \der_n$ are skew commuting skew
derivations (i.e. $\der_i\der_j=-\der_j\der_i$), the following
equality is obvious
$$ \sum_{\alpha \in \CB_n} \L_n \der^\alpha = \sum_{\alpha \in
\CB_n} \der^\alpha \L_n.$$ Now, $\End_K(\L_n)= S:= \sum_{\alpha
\in \CB_n} \L_n \der^\alpha$ since $E_{\alpha \beta} =  x^\alpha
\phi \der^\beta\in S$. The sum $S$ is a direct sum: suppose that
$\sum u_\alpha \der^\alpha =0$ for some elements $u_\alpha \in
\L_n$ not all of which are equal to zero, we seek a contradiction.
Let $u_\beta$ be a nonzero element with $|\beta | =\beta_1+\cdots
+ \beta_n$ the least possible. Then $0=\sum u_\alpha \der^\alpha
(x^\beta )= u_\beta$, a contradiction. $\Box $

Let $A$ be a ring. For $a,b\in A$,  $\{ a,b\} := ab+ba$ is called
the {\em anti-commutator} of elements $a$ and $b$.

\begin{corollary}\label{cE12Sep06}
The ring $\End_K(\L_n )$ is generated over $K$ by elements $x_1,
\ldots , x_n, \der_1, \ldots , \der_n$ that satisfy the following
defining relations: for all $i,j$,
\begin{eqnarray*}
 x_i^2=0,& x_ix_j=-x_jx_i,\\
 \der_i^2=0,& \der_i\der_j=-\der_j\der_i,\\
 \der_ix_j+x_j\der_i=\d_{ij}, & {\rm the \; Kronecker \; delta}.
\end{eqnarray*}
\end{corollary}

{\it Proof}. It is straightforward that these relations hold.
Theorem \ref{E12Sep06}.(2) implies that these relations are
defining. $\Box $

 By Corollary
\ref{cE12Sep06}, we have the  duality $K$-automorphism of the ring
${\rm End}_K(\L_n)$:
$$\Delta : x_i\mapsto \der_i, \;\; \der_i\mapsto x_i, \;\; i=1,
\ldots , n.$$ Clearly, $\D^2= {\rm id}$.  Similarly, by Corollary
\ref{cE12Sep06}, we have the duality $K$-anti-automorphism of the
ring ${\rm End}_K(\L_n)$:
$$\nabla : x_i\mapsto \der_i, \;\; \der_i\mapsto x_i, \;\; i=1,
\ldots , n.$$ Clearly, the anti-automorphism  $\D$ is an
involution, i.e. $\nabla^2= {\rm id}$, $\nabla (ab)= \nabla (b)
\nabla (a)$.

We say that the element $2:=1+1\in K$ is {\em regular} if $2\l =0$
for some $\l \in K$  implies $\l =0$. If $K$ contains a field then
$2$ is regular iff the characteristic of $K$ is not $2$. If $K$ is
a commutative ring then the ring $\L_n(K)$ is non-commutative iff
$n\geq 2$ and $2\neq 0$ in $K$. A commutative ring is called {\em
reduced} if $0$ is the only nilpotent element of the ring.


\begin{lemma}\label{z9Sep06}
If  the ring $K$ is commutative, $2\in K$ is regular, and $n\geq
2$,  then the centre of $\L_n(K)$ is equal to
$$ Z(\L_n)=\begin{cases}
\Lnz , & \text{if $n$ is even},\\
\Lnz\oplus Kx_1\cdots x_n, & \text{if $n$ is odd }.
\end{cases}$$
\end{lemma}

{\it Proof}. Since $x_1, \ldots , x_n$ are $K$-algebra generators
for $\L_n$, an  element $u\in \L_n$ is central iff it commutes
with all $x_i$. Now, the result is obvious due to (\ref{xiaib})
and the fact that $x_i\cdot (x_1\cdots x_n)=0$, $i\geq 1$.  $\Box
$

Let $\Lnz':= \oplus Kx^\alpha$ where $\alpha$ runs through all
even subsets of $\CB_n$ distinct from $\{ 1, 2, \ldots, n\}$. So,
$\Lnz'\subseteq \Lnz\subseteq Z(\L_n)$, and $ \Lnz'= \Lnz$ iff $n$
is odd.

Let $G:= \Aut_K(\L_n(K))$ be the group of $K$-automorphisms of the
ring $\L_n(K)$. Each $K$-automorphism $\s\in G$ is a uniquely
determined by the images of the canonical generators:
$$ x_1':= \s (x_1), \ldots , x_n':= \s (x_n).$$
Note that $x_1' \ldots , x_n'$ is another set of canonical
generators for $\L_n$.

Till the end of this section, let $R$ be a {\em commutative} ring.
Consider the subgroup $G_{gr}$ of $G$, elements of which preserve
the $\Z $-grading of $\L_n$:
$$ G_{gr}:= \{ \s \in G\, | \, \s (\L_{n,i})= \L_{n,i} \;\;{\rm
for\; all}\; i\in \Z \}= \{ \s \in G\, | \, \s (\L_{n,1}) =
\L_{n,1}\}.$$ The last equality is due to the fact that the
$\L_{n,1}$ generates $\L_n$ over $K$. Clearly, $\s \in G_{gr}$ iff
$\s = \s_A$ where $\s_A(x_i)= \sum_{j=1}^na_{ij}x_j$ for some $A=
(a_{ij})\in \GL_n(K)$. This can be written in the matrix form as
$\s (x)= Ax$ where $x:=(x_1, \ldots , x_n)^t$ is the vector-column
of indeterminates. Since $\s_A\s_B= \s_{BA}$, the group $G_{gr}$
is canonically isomorphic to the group $\GL_n(K)^{op}$ {\em
opposite} to $\GL_n(K)$ via the map $\GL_n(K)^{op}\ra G_{gr}$,
$A\mapsto \s_A$. We identify the group $G_{gr}$ with
$\GL_n(K)^{op}$ via this isomorphism. Note that
$\GL_n(K)\ra\GL_n(K)^{op}$, $ A\mapsto A^{-1}$, is the isomorphism
of groups. One can write
$$G_{gr}= \{
\s_A\, | \, A\in \GL_n(K)\}.$$

From this moment and till the end of this section,  $K$ is a {\em
reduced commutative} ring (if it is not stated otherwise). For $\s
\in G$, let $x_i':= \s (x_i)$. Then $x_i'^2=\s(x_i^2)= \s (0)=0$.
If $ \l_i\equiv x_i'\mod \gm$ for some $\l_i\in K$ then
$\l_i^2=0$, hence $\l_i=0$ since $K$ is reduced. Therefore, $\s
(\gm ) = \gm$, and so 
\begin{equation}\label{smi=mi}
\s (\gm^i)=\gm^i \;\; {\rm for \; all}\;\;  i\geq 1.
\end{equation}
This proves the next lemma.

\begin{lemma}\label{s12Sep06}
If $K$ is a reduced commutative ring and $\s \in G$ then $\s (x) =
Ax+b$ for some $ A\in \GL_n(K)$ and $b:= (b_1, \ldots , b_n)^t$
where all $b_i\in \gm^2$.
\end{lemma}

Consider the following subgroup of $G$,
$$ U:= \{ \s \in G\, | \, \s (x) = x+ b\;{\rm for \;some}\; b\in (\gm^2)^{\times n}\} = \{ \s \in G\, | \, (\s
-1) (\gm ) \subseteq \gm^2\}.$$ For each $ \s \in G$ written as
$\s (x) =Ax+ b$ and each $\tau \in U$, we have $\s\s_{A^{-1}},
\s_{A^{-1}}\s \in U$, and $\s \tau \s^{-1} \in U$. These mean that
$U$ is a {\em normal} subgroup of $G$ such that 
\begin{equation}\label{G=GU}
G=G_{gr}U= UG_{gr}, \;\; G_{gr}\cap U=\{ e\}.
\end{equation}
Therefore, $G$ is a {\em skew product} of the groups $G_{gr}$ and
$U$: 
\begin{equation}\label{1G=GU}
G=G_{gr}\ltimes U= \GL_n(K)^{op}\ltimes U\simeq \GL_n(K)\ltimes U.
\end{equation}
For each $i\geq 2$, consider the subgroup $U^i$  of $U$:
\begin{equation}\label{Uidef}
U^i:= \{ \s \in U\, | \, (\s - 1) (\gm )\subseteq \gm^i\}.
\end{equation}
By the very definition,  $U=U^2\supset U^3 \supset \cdots \supset
U^n\supset U^{n+1} = \{ e\}$. Note that $\s \in U^i$ iff, for all
 $j$, $\s (x_j) = x_j+m_j$ for some $m_j\in \gm^i$. Recall that
 $\s (\gm^i)= \gm^i$ for all $\s \in G$ and  $i\geq 1$ (since $K$ is
  a reduced commutative ring). For any
 $\tau \in G$, $\s \in U^i$, and $x_j$,
 $$ \tau \s \tau^{-1} (x_j) \equiv  \tau (1+\s -1) \tau ^{-1} (x_j)
 \equiv \tau \tau^{-1} (x_j) \equiv x_j\mod \gm^i.$$
 Hence, each $U^i$ is a {\em normal} subgroup of $G$.  Each
 factor group $U^i/ U^{i+1}$ is  abelian: Let $\s , \tau \in
 U^i$, $\s (x_j)= x_j+a_j+\cdots$ and $\tau (x_j) =
 x_j+b_j+\cdots$ for some     elements $a_j, b_j\in \L_{n,i}$ and
 three dots denote elements of $\gm^{i+1}$. Then $\s \tau (x_j) =
 x_j+ a_j+b_j+ \cdots$ and $\tau \s (x_j) = x_j+b_j+a_j+\cdots $.
 Therefore, $\s \tau U^{i+1} = \tau \s U^{i+1}$, and  $U$ is
 a  nilpotent group.

Let $K^n$ be the direct sum of $n$ copies of the  additive group
$(K, +)$. Let $\th := x_1 \cdots x_n$. The map 
\begin{equation}\label{KnUn}
K^n\ra U^n, \;\; \l = (\l_1, \ldots, \l_n) \mapsto (\s_\l:
x_i\mapsto x_i+\l_i\th ),
\end{equation}
is an isomorphism of groups ($\s_{\l +\mu }= \s_\l \s_\mu $). This
follows directly from the fact that $\gm \th = \th \gm =0$. So,
$U^n= \{ \s_\l \, | \, \l \in K^n\}$. One can easily verify that
for any $\s \in G$ with $\s (x) = Ax+b$, $A= (a_{ij}) \in
\GL_n(K)$,
\begin{equation}\label{slsA}
\s^{-1} \s_\l \s  = \s_{\frac{A\l}{\det (A)}}.
\end{equation}
Indeed,
\begin{eqnarray*}
 \s^{-1} \s_\l \s (x_i) &=&\s\s_\l (\sum_{j=1}^na_{ij}x_j+b_i)= \s^{-1} (\sum_{j=1}^na_{ij}x_j +b_i +\sum_{j=1}^n a_{ij} \l_j \th ) \\
 &=& \s^{-1} (\s (x_i) +\sum_{j=1}^n a_{ij} \l_j \th ) = x_i +\frac{\sum_{j=1}^n a_{ij} \l_j}{\det
 (A)}\th  .
\end{eqnarray*}
The equality (\ref{slsA}) describes completely the group structure
of
 $\Aut_K(\L_2)$ (where  $K$ is
  a reduced commutative ring since $U= U^2$):  
\begin{equation}\label{AL2}
\Aut_K(\L_2)= \GL_2(K)^{op} U=\{ \s_A\s_\l \, | \, \s_A\in
\GL_2(K)^{op}, \s_\l \in U\}, \;\; \s_A\s_\l \cdot \s_B\s_\mu =
\s_{BA} \s_{\frac{B\l}{\det (B)}+ \mu}.
\end{equation}

An element $a(x_1, \ldots , x_n) = \l +\cdots \in \L_n$ is a unit
iff $a(0, \ldots , 0)= \l \in K$ is a unit. For each unit $a\in
\L_n$, the map $ \o_a(x):= axa^{-1}$ is called the {\em inner}
automorphism of $\L_n$. Since $\o_a(x_i)= (\l +\cdots )
x(\l+\cdots )^{-1} = \l x_i\l^{-1} +\cdots = x_i+\cdots $ we have
$\o_a\in U$, i.e. the group ${\rm Inn} (\L_n)$ of all the inner
automorphisms of $\L_n$ is a subgroup of $U$, 
\begin{equation}\label{InLU}
{\rm Inn} (\L_n) \subseteq U.
\end{equation}
Let us denote the automorphism $a\mapsto \oa$ of $\L_n$ by $s$.
Recall that $\Lnz = \{ a\in \L_n\, | \, s(a)=a\}$ and $\Ln1 = \{
a\in \L_n\, | \, s(a)=-a\}$. Consider the subgroup $G_{\Z_2-gr}$
of $G$ elements of which respect $\Z_2$-grading on $\L_n$:
$$G_{\Z_2-gr}:= \{ \s \in G\, | \, \s (\Lnz ) = \Lnz , \; \s (\Ln1
)= \Ln1\}.$$ Clearly, 
\begin{equation}\label{GZ2s}
G_{\Z_2-gr}= \{ \s \in G\, | \, \s s = s \s\}.
\end{equation}
Then 
\begin{equation}\label{GS3s}
\G := U\cap G_{\Z_2-gr}= \{ \s \in U\, | \, \s s = s\s\} = \{ \s
\in U\, | \, \s (x_i) \in \Ln1 , 1\leq i\leq n \}
\end{equation}
 is the subgroup of $U$
(the last equality follows easily from the fact that the set
$\Ln1$ generates the $K$-algebra $\L_n$ and that $\L_n = \Lnz
\oplus \Ln1$ is a $\Z_2$-graded $K$-algebra). So, $\s \in \G $
iff, for each $i=1, \ldots , n$, 
\begin{equation}\label{sGS3s}
\s (x_i) = x_i+a_{i,3}+a_{i,5}+\cdots + a_{i, 2j+1}+\cdots , \;
a_{i, 2j+1}\in \L_{n, 2j+1},
\end{equation}
all summands are odd. Note that for arbitrary commutative ring $K$
(not necessarily reduced) the set of automorphisms $\s$ from
(\ref{sGS3s}) is a group $\G$.  Clearly, $\GL_n(K)^{op}\subseteq
G_{\Z_2-gr}$ and $\GL_n(K)^{op}\cap \G =\{ e\}$ since
$\GL_n(K)^{op}\cap U= \{ e\}$ and $\G \subseteq U$.
 The group $\G $ (over an arbitrary commutative ring $K$) can be defined as
\begin{equation}\label{GS4s}
  \G = \{ \s \in G\, | \, \s (x_i) -x_i\in \Ln1 \cap \gm^3, \; i\geq
 1\}.
\end{equation}

The group $\G$ is endowed with the descending chain of its normal
subgroups
$$ \G = \G^2\supseteq \G^3\supseteq \cdots \supseteq \G^i:= \G
\cap U^i\supseteq \cdots \supseteq \G^n \supseteq \G^{n+1} = \{
e\}.$$ Since $\G^i= \{ \s \in \G \, | \, \s (x_j) - x_j\in \L_{n,
\geq i }^{od}, \; j=1, \ldots , n\}$, it is obvious that
$$ \G = \G^2=\G^3\supset \G^4=\G^5\supset \cdots \supset \G^{2m}=\G^{2m+1}\supset \cdots   .$$
Recall that $[x,y]:= xy-yx$ is  the commutator of elements $x$ and
$y$, and $\o_s: t\mapsto sts^{-1}$ is an inner automorphism.
\begin{lemma}\label{m28Sep06}
Let $K$ be a commutative ring.
\begin{enumerate}
\item For each $a\in \Lnod$, $a^2=0$. \item $[ \Lnod , \L_n ]
\subseteq \Lnev \subseteq Z(\L_n)$ and $ [ \Lnod , [ \Lnod ,
\L_n]]=0$. \item For each $a\in \Lnod$ and $x\in \L_n$, $\o_{1+a}
(x) = x+[a,x]$. \item For each $ a,b\in \Lnod$, $\o_{1+a} \o_{1+b}
= \o_{1+a+b} = \o_{1+b}\o_{1+a}$ and $\o_{1+a}^{-1} = \o_{1-a}$.
\item The map $ \o : \Lnod \ra U$, $ a\mapsto \o_{1+a}$, is a
homomorphism of groups. It is a monomorphism if $n$ is even and
has the  kernel $\ker (\o )= Kx_1\cdots x_n$ if $n$ is odd. \item
For each $a\in \L_n$, $a\oa = \oa a = a_0^2$ where $a= a_0+a_1$,
$a_0\in \Lnod$, $a_1\in \Lnev$, $\oa:= a_0-a_1$.
\end{enumerate}
\end{lemma}

{\it Proof}. $1$. For $a\in \Lnod$, $a=\sum \l_\alpha x^\alpha$
where $\alpha $ runs through non-empty {\em odd} subsets of the
set $\{ 1, \ldots , n\}$. For any two such subsets $\alpha $ and
$\beta $, $x^\alpha x^\beta = - x^\beta x^\alpha$ and $(x^\alpha
)^2=0$. Now, $a^2= \sum_{\alpha \neq \beta} \l_\alpha \l_\beta
(x^\alpha x^\beta +x^\beta x^\alpha)+\sum_{\alpha} \l_\alpha^2
(x^\alpha )^2=0$.

2. $[ \Lnod , \L_n ] = [ \Lnod , \Lnod + Z(\L_n) ]\subseteq [
\Lnod , \Lnod ]  \subseteq \Lnev \subseteq Z(\L_n)$. Then the
second equality is obvious.

3. It suffices to prove the equality for monomials $x= x^\alpha$.
If $x^\alpha$ is even, hence central, the equality is obvious. If
$x^\alpha$ is odd then
\begin{eqnarray*}
 \o_{1+a}(x^\alpha )&=& (1+a) x^\alpha (1+a)^{-1}\\
 &=& (1+a) x^\alpha (1-a) \;\;\;\;\;\;\;  ({\rm by\; statement} \; 1) \\
 &=& x^\alpha + [ a, x^\alpha ] - ax^\alpha a = x^\alpha + [ a,
 x^\alpha ] + a^2x^\alpha \\
 &=& x^\alpha + [ a, x^\alpha ] \;\;\;\;\;\;\; \;\;\;\;\;\;\;\; ({\rm by \; statement} \; 1).
\end{eqnarray*}
4. For each $a,b\in \Lnod$ and $x\in \L_n$,
\begin{eqnarray*}
 \o_{1+a}\o_{1+b}(x)&=& x+[a+b, x] +[a, [b, x]]= x+[a+b,
 x] \;\;\;  ({\rm by\; statement} \; 2)\\
 &=& \o_{1+a+b} (x) = \o_{1+b+a} (x) = \o_{1+b} \o_{1+a} (x).
\end{eqnarray*}
 $\o_{1+a}^{-1} = \o_{(1+a)^{-1}} =\o_{1-a}$ since $a^2=0$ (statement 1).

5. By statement 4, the map is a group homomorphism. By statement
3, an element $a\in \Lnod$ belongs to the kernel iff $[a, x_i]=0$
for all $i$ iff $a\in Z(\L_n)= \Lnev +Kx_1\cdots x_n$. Now, the
result is obvious since $a$ is odd.

6. $a\oa = a_0^2-a_1^2= a_0^2$ and $ \oa a = a_0-a_1^2=a_0^2$
since $a_1^2=0$, by statement 1. $\Box $

Let $ \O$ be the image of the group homomorphism $\o$ in Lemma
\ref{m28Sep06}.(5), 
\begin{equation}\label{GImo}
\O := {\rm im } (\o ) = \{ \o_{1+a} \, | \, a\in \Lnod \}.
\end{equation}

\begin{lemma}\label{i30Sep06}
Let $K$ be a commutative ring. Then $\O = {\rm Inn}(\L_n)$. In
particular, the group ${\rm Inn} (\L_n)$ of inner automorphisms of
the Grassmann algebra $\L_n(K)$ is an abelian group.
\end{lemma}

{\it Proof}. Let $u$ be a unit of $\L_n$. Then $u = \l + a+b$ for
some $0\neq \l \in K$, an odd element $a\in \gm$, and  an even
element $b\in \gm$. Note that $\l +b$ is a central element and
that the element $a':= \frac{a}{\l + b}\in \gm$ is odd. Now, $\o_u
= \o_{ (\l +b)(1+a')}= \o_{\l +b} \o_{1+a'}= \o_{1+a'}\in \O$.
Therefore, $\O = {\rm Inn}(\L_n)$.  $\Box $

So, $\O$ is a {\em normal abelian} subgroup of $G$.  Since $\O =
{\rm Inn} (\L_n)\subseteq U$, by (\ref{InLU}), the group $\O$ is
endowed with the induced filtration
$$ \O = \O^2\supseteq \O^3\supseteq \cdots \supseteq \O^i:= \O
\cap U^i \supseteq \cdots \supseteq\O^{n-1} \supseteq \O^n= \{
e\}.$$ Note that $\O = \O^2\supset  \O^3 =\O^4  \supset \O^5 =
\O^6  \supset\cdots \supset \O^{2i+1}= \O^{2i+2}\supset \cdots $.
 By Lemma \ref{m28Sep06}.(5), the group $\O$ is canonically
 isomorphic to the factor group $\Lnod / \Lnod \cap Kx_1\cdots
 x_n$. Under this isomorphism the filtration $\{ \O^i\}$ coincides
 with the filtration on $\Lnod / \Lnod \cap Kx_1\cdots
 x_n$ shifted by $-1$ that is the induced filtration from the $\gm$-adic
 filtration of the ring $\L_n$, i.e. $\O^{2m}= \{ \o_{1+a} \, | \,
 a\in \Lnod \cap \gm^{2m-1}\}, m\geq 1$. Note that (Lemma
 \ref{m28Sep06}.(3))
\begin{equation}\label{1GImo}
\o_{1+a} (x_i) - x_i = [a,x_i] \in \Lnev\;\; {\rm for \; all}\;\;
a\in \Lnod\;\;{\rm and}\;\; i.
\end{equation}

{\bf Generators for and dimension of the algebraic group $U$}.
Define the function $\d_{\cdot , ev}: \N \ra \{ 0,1\}$ by the rule
$$ \d_{n, ev}=\begin{cases}
1, & \text{if $n$ is even},\\
0, & \text{if $n$ is odd}.
\end{cases}$$
For each $n\geq 2$, $[\frac{n}{2}]-\d_{n, ev}$ (resp.
$[\frac{n}{2}]$) is the number of odd (resp. even) numbers $m$
such that $ 2\leq m \leq n$.

\begin{theorem}\label{11Sep06}
Let $K$ be a commutative ring in statement 2, and let
 $K$ be a reduced commutative ring with $\frac{1}{2}\in K$ in statements 1, 3, and
 4. Let $\L_n = \L_n(K)$, $U= U(\L_n)$,  and $n\geq 2$. Then
\begin{enumerate}
\item the associated graded group $ \prod_{i\geq 2} U^i/U^{i+1}$
is isomorphic to the direct product of $d_n$ copies of the
additive group $K$ where
 $d_n= n\sum_{m=1}^{[\frac{n}{2}]-\d_{n,ev}}{n\choose 2m+1}+\sum_{m=1}^{[\frac{n}{2}]}{n\choose
2m-1}$. In more detail, \item for each $m=1, \ldots ,
[\frac{n}{2}]-\d_{n, ev}$, the map $$ (\L_{n, 2m+1})^n\ra
U^{2m+1}/U^{2m+2}, \; a=(a_1, \ldots , a_n)\mapsto \s_aU^{2m+2},
$$
is a group isomorphism where $\s_a\in U^{2m+1}: x_i\mapsto
x_i+a_i$; \item for each $m=1, \ldots , [\frac{n}{2}]$, the map $$
\L_{n, 2m-1}\ra U^{2m}/U^{2m+1}, \; a\mapsto \o_{1+a}U^{2m+1},
$$
is a group isomorphism where $\o_{1+a}\in U^{2m}$ is the inner
automorphism of $\L_n: x\mapsto(1+a)x(1+a)^{-1}$. \item All the
elements $\s_a$ and $ \o_{1+b}$ from statements 2 and 3 are
generators for the group $U$.
\end{enumerate}
\end{theorem}

{\it Proof}. $1$.  The first statement follows directly from
statements 2 and 3 since $\L_{n,i} \simeq K^{{n\choose i}}$ for
all $i=0, 1, \ldots , n$.

$4$. This statement follows from statements 1--3.

$2$. One can easily see that $\s_a\in U^{2m+1}$ for any $a\in
(\L_{n, 2m+1})^n$; and $ \s_a\s_bU^{2m+2}= \s_{a+b}U^{2m+2}$ for
any $a$ and $b$. By the very definition, the map $a\mapsto
\s_aU^{2m+2}$ is an injection. It suffices to show that it is a
surjection. Let $\s$  be an arbitrary element of $U^{2m+1}$. Then
$\s (x_i)= x_i+a_i+\cdots$ for some $a_i\in \L_{n, 2m+1}$ where
the the dots mean bigger terms with respect to the $\Z$-grading on
$\L_n$. Then
$$ \s_{-a} \s (x_i)= \s_{-a} ( x_i+a_i+\cdots ) =
x_i-a_i+a_i+\cdots = x_i+\cdots,$$ hence $\s_{-a}\s \in U^{2m+2}$,
i.e. $\s  U^{2m+2}= \s_a U^{2m+2}$, and we are done.

$3$.  Clearly, $\o_{1+a} \in U^{2m}$ since $ \o_{1+a}(x_i) = x_i +
[a,x_i]$ (Lemma \ref{m28Sep06}.(3)). By Lemma \ref{m28Sep06}.(5),
the map $a\mapsto \o_{1+a} U^{2m+1}$ is a group homomorphism. By
the very definition, this map is injective. It suffices to show
that it is surjective. This will be done in the next lemma in the
proof of which an algorithm is given of how, for a given element
$\s$ of $U^{2m}/ U^{2m+1}$,  to find an element $a\in \L_{n,
2m-1}$ such that $\s \equiv \o_{1+a} U^{2m+1}$.

\begin{lemma}\label{a11Sep06}
We keep the assumptions of Theorem \ref{11Sep06}.(3). For each
$m=1, \ldots , [\frac{n}{2}]$, and each $\s \in U^{2m}$, there
exist elements $c_{i+1} \in K \lfloor x_{i+1} ,\ldots ,
x_n\rfloor_{2m-i}$, $i=1, \ldots, 2m$, such that $\o \s \in
U^{2m+1}$ where $ \o := \o_{1+x_1\cdots x_{2m-1}
c_{2m+1}}^{-1}\o_{1+x_1\cdots x_{2m-2} c_{2m}}^{-1}\cdots
\o_{1+x_1\cdots x_{i-1} c_{i+1}}^{-1} \cdots
\o_{1+x_1c_3}^{-1}\o_{1+c_2}^{-1}$.
\end{lemma}

{\it Remark}. If $n=2m$ then $c_{2m+1}\in K$.

{\it Proof}. Since each element $x_1 \cdots x_{i-1} c_{i+1}$ is a
homogeneous element of $\L_n$ of degree $i-1+2m -i= 2m-1$, each
automorphism $\o_{1+x_1 \cdots x_{i-1} c_{i+1}}$ belongs to the
group $U^{2m}$, and so does their product $\o$. For each $i=1,
\ldots , n$, let $ x_i':= \s (x_i) = x_i+a_i +\cdots $ for some
$a_i\in \L_{n, 2m}$. We prove the lemma in several steps.

{\it Step 1. For each $i=1, \ldots , n$, $a_i= x_i b_i$ for some
element} $b_i\in K\lfloor x_1, \ldots , \widehat{x_i}, \ldots ,
x_n\rfloor_{2m-1}$.  Note that each element $a_i \in \L_{n, 2m}
\subseteq Z(\L_n)$ is central. Then $x_i'^2= \s (x_i^2) =0$, and
so
$$ 0=x_i'^2= (x_i+a_i+\cdots )^2= 2x_ia_i+\cdots $$
hence $x_ia_i=0$ since $\frac{1}{2}\in K$; and so  $a_i= x_ib_i$
for some element $b_i\in K\lfloor x_1, \ldots , \widehat{x_i},
\ldots , x_n\rfloor_{2m-1}$.

{\it Step 2. Let us prove that for each pair} $i\neq j$,
\begin{equation}\label{st2}
b_i|_{x_j=0}= b_j|_{x_i=0}.
\end{equation}
By Step 1, $x_i'= x_i(1+b_i)+\cdots $ for all $i$. Since $b_i,
b_j\in \L_{n, 2m-1}$ and $2m-1\geq 2-1\geq 1$, we see that the
homogeneous elements $b_i-b_j$ and $b_ib_j$ have distinct degrees
(if the elements are nonzero). Computing separately both sides of
the equality $x_i' x_j'= - x_j'x_i'$ we have
\begin{eqnarray*}
x_i' x_j'&= & (x_i(1+b_i)+\cdots ) ( x_j(1+b_j)+\cdots )=  x_ix_j(1-b_i+b_j)+\cdots , \\
 -x_j' x_j'&= & -
 x_jx_i(1-b_j+b_i)+\cdots =x_ix_j(1+b_i-b_j)+\cdots .
\end{eqnarray*}
Since the degree of the elements $b_i-b_j$ is $2m-1\geq 2-1\geq 1$
and $ \frac{1}{2}\in K$, we must have 
\begin{equation}\label{1st2}
x_ix_j(b_i-b_j)=0.
\end{equation}
This equality is equivalent to the equality $b_i|_{x_i=0, x_j=0}=
b_j|_{x_i=0, x_j=0}$ which is obviously equivalent to (\ref{st2})
(since each $b_k$ does not depend on $x_k$, by Step 1).

{\it Step 3}. We are going to prove that {\it one can choose
elements $c_{i+1} \in K\lfloor x_{i+1} , \ldots ,
x_n\rfloor_{2m-i}$, $i=1, \ldots , 2m$, such that for each} $i=1,
\ldots , 2m$, 
\begin{equation}\label{w1xi}
\o_{1+x_1\cdots x_{i-1} c_{i+1}}^{-1}\cdots
\o_{1+x_1c_3}^{-1}\o_{1+c_2}^{-1}(x_j')\equiv x_j \mod\gm^{2m+1},
\;\; j=1, \ldots , i.
\end{equation}
 We use induction on $i$.
Briefly, (\ref{w1xi}) will follow from (\ref{st2}). Note that
$b_1\in K\lfloor x_2, \ldots , x_n\rfloor_{2m-1}$. Then
$c_2:=-\frac{1}{2} b_1\in K\lfloor x_2, \ldots ,
x_n\rfloor_{2m-1}$, and
$$ x_1'= x_1+x_1b_1+\cdots = x_1-b_1x_1+\cdots = x_1+[-\frac{1}{2}
b_1, x_1]+\cdots = \o_{1+c_2} (x_1) +\cdots, $$ hence
$\o_{1+c_2}^{-1} (x_1)\equiv x_1\mod \gm^{2m+1}$.

Changing $\s$ to $\o_{1+c_2}^{-1}\s$, one can assume that $x_1'=
x_1+\cdots$, i.e. $b_1=0$. Applying (\ref{st2}) in the case $j=1$
and $i\geq 2$, we have $b_i|_{x_1=0} =0$, hence $b_i=-2x_1c_{i+1}$
for unique $c_{i+1}\in K\lfloor x_2, \ldots \widehat{x_i}, \ldots
, x_n\rfloor_{2m-2}\subseteq Z(\L_n)$. Now,
$$ x_2'= x_2+x_2(-2x_1c_3)+\cdots = x_2+[x_1c_3,
x_2]+\cdots = \o_{1+x_1c_3}(x_2) +\cdots , $$ hence
$\o_{1+x_1c_3}^{-1}(x_2')= x_2+\cdots$. Note that
$\o_{1+x_1c_3}^{-1}(x_1')= x_1+\cdots$ since
$$\o_{1+x_1c_3}(x_1')=\o_{1+x_1c_3}(x_1)+\cdots = x_1+[x_1c_3,
x_1]+\cdots = x_1+\cdots.$$ This proves (\ref{w1xi}) for $i=2$.

Let $i\geq 3$, and suppose that we have found already elements
$c_{k+1} \in K\lfloor x_{k+1} , \ldots x_n\rfloor_{2m-k}$, $k=1,
\ldots , i$, that  satisfy (\ref{w1xi}). We have to find
$c_{i+2}$. Changing, if necessary, $\s $ for $\o_{1+x_1\cdots
x_{i-1} c_{i+1}}^{-1}\cdots \o_{1+x_1c_3}^{-1}\o_{1+c_2}^{-1}\s$
one can assume that $x_k'= x_k+\cdots $ for $k=1, \ldots , i$,
i.e. $b_k=0$ for $k=1, \ldots , i$. By (\ref{st2}),
$$ b_{i+1}|_{x_k=0}= b_k|_{x_{i+1}=0}=0, \;\; k=1, \ldots , i,$$
hence $b_{i+1} = -2x_1\cdots x_ic_{i+2}$ for a unique element
$c_{i+2} \in K\lfloor x_{i+2}, \ldots , x_n\rfloor_{2m-i-1}$. Now,
\begin{eqnarray*}
x_{i+1}'&=& x_{i+1} +x_{i+1} (-2x_1\cdots x_ic_{i+2})
+\cdots = x_{i+1} + [x_1\cdots x_ic_{i+2}, x_{i+1}]+\cdots \\
&=& \o_{1+x_1\cdots x_ic_{i+2}}(x_{i+1})+\cdots,
\end{eqnarray*}
hence $\o^{-1}_{1+x_1\cdots x_ic_{i+2}}(x_{i+1}')= x_{i+1}+\cdots
$. Note that $\o^{-1}_{1+x_1\cdots x_ic_{i+2}}(x_k')= x_k+\cdots
$, $k=1, \ldots , i$, since
$$\o^{-1}_{1+x_1\cdots x_ic_{i+2}}(x_k')=\o^{-1}_{1+x_1\cdots x_ic_{i+2}}(x_k+\cdots
)=x_k+ [x_1\cdots x_ic_{i+2}, x_k] +\cdots = x_k+\cdots.$$ So,
(\ref{w1xi}) holds for $i+1$. By induction, (\ref{w1xi}) holds for
all $i=1, \ldots , 2m$. In particular, it does for $i=2m$. Then
changing, if necessary, $\s$ for $\o^{-1}_{1+x_1\cdots
x_{2m-1}c_{2m+1}} \cdots \o_{1+x_1c_3}^{-1}\o_{1+c_2}^{-1}\s$ one
can assume that $x_k'= x_k+\cdots$ for $k=1, \ldots , 2m$, i.e.
$b_k=0$ for $k=1, \ldots , 2m$. Note that in order to prove Lemma
\ref{a11Sep06}, we have to show that $b_k=0$ for $k=1, \ldots ,
n$. If $n=2m$ we are done. If $n>2m$ then by (\ref{st2}) for each
$i>2m$ and $j=1, \ldots , 2m$: $b_i|_{x_j=0}= b_j|_{x_i=0}=0$.
Hence, $b_i\in (x_1\cdots x_{2m})$, but $b_i\in \L_{n, 2m-1}$,
therefore $b_i=0$. This proves Lemma \ref{a11Sep06} and Theorem
\ref{11Sep06}. $\Box $

In Step 3, it was, in fact, proved that the conditions
(\ref{w1xi}) {\em uniquely} determines the elements $c_2, \ldots ,
c_{2m+1}$ (the idea of finding the elements  $c_i$  is to kill the
`leading term', this determines uniquely $c_i$ by the expression
given in the proof above). So, Lemma \ref{a11Sep06} can be
strengthened as follows.

\begin{corollary}\label{ca11Sep06}
The elements $c_2, \ldots , c_{2m+1}$ from Lemma \ref{a11Sep06}
are unique provided (\ref{w1xi}) holds for all $1\leq j \leq i\leq
2m$.
\end{corollary}

{\it Proof}. This fact also follows at once from Lemma
\ref{m28Sep06} and Theorem \ref{14Sep06}.(1). $\Box$

\begin{corollary}\label{c11Sep06}
Let $K$ be a reduced commutative ring with $\frac{1}{2}\in K$, and
$\s \in U$. Then $\s $ is a (unique) product $\s =\cdots
\o_{1+a_5}\s_{b_5}\o_{1+a_3}\s_{b_3}\o_{1+a_1}$ for unique
elements $a_i\in \L_{n,i}$ and $b_j= (b_{j1},\ldots , b_{jn})\in
\L_{n,j}^n$.
\end{corollary}

By Corollary \ref{c11Sep06}, the coefficients of the elements
$\ldots , a_5, b_5, a_3, b_3, a_1$ are coordinate functions for
the algebraic group $U$ over $K$. Therefore, the $K$-algebra of
 (regular) functions of the algebraic group $U$ is a polynomial
algebra over $K$ in $d_n$ variables where $d_n$ is defined in
Theorem \ref{11Sep06}.

It follows from Corollary \ref{c11Sep06} that 
\begin{equation}\label{UiGi}
U^i= \O^i\rtimes \G^i, \; i\geq 2.
\end{equation}

{\bf The group structure of $G:= {\rm Aut}_K(\L_n(K))$}.

\begin{theorem}\label{29Sep06}
Let $K$ be a reduced commutative ring with $\frac{1}{2}\in K$.
Then
\begin{enumerate}
\item The group $U= \O \rtimes \G$ is the semi-direct product of
its subgroups ($\O$ is a normal subgroup of $U$, $\O \cap \G = \{
e\}$, and $U= \O \G$).
 \item  $G= U\rtimes \GL_n(K)^{op}$ ($U$ is a normal subgroup of $U$, $U \cap \GL_n(K)^{op}= \{
 e \}$, and $G= U\GL_n(K)^{op}$).\item $G= (\O \rtimes \G ) \rtimes
\GL_n(K)^{op}$.
\end{enumerate}
\end{theorem}

{\it Proof}. 1. By (\ref{GS4s}) and (\ref{1GImo}), $\O \cap \G =
\{ e\}$. For any $\g \in \G$ and $\o_{1+a} \in \O$ (resp.
$\o_{1+a} \in \O^i := \O \cap U^i$, $i\geq 2$) 
\begin{equation}\label{g1a}
\g \o_{1+a} \g^{-1} = \o_{\g (1+a)}\in \O \;\;\; ({\rm resp. }\;\;
\g \o_{1+a} \g^{-1} = \o_{\g (1+a)}\in \O^i, \;\; i\geq 2),
\end{equation}
since $\G \Lnod \subseteq \Lnod$ and $1+a\in \Lnod$ (resp. and
$\G\gm^i \subseteq \gm^i$, $i\geq 1$). So, in order to finish the
proof of statement 1 it is enough to show that $U= \G \O$. This
equality follows from Lemma \ref{c11Sep06} and (\ref{g1a}).


$2$. See (\ref{1G=GU}).

3. This follows from statements 1 and 2.  $\Box $

Let $B= B_n$ be the set of all the $n$-tuples (columns) $b= (b_1,
\ldots , b_n)^t$ where all  $b_i$ are arbitrary odd elements of
$\L_n$ of the form $b_1= x_1+\cdots, \ldots , b_n = x_n+\cdots$
where the three dots mean bigger terms. Then
 $$\G = \{ \g_b\, | \, b\in B, \; \g_b (x_i) = b_i, \; i=1, \ldots
 , n\}$$ and the map $B\ra \G$, $b\mapsto \g_b$ is a
 bijection. The product of two elements $\g_b, \g_c\in \G$ is
 given by the rule
 $$ \g_b\g_c= \g_{c\circ b}$$
 where $c\circ b$ is the  composition of functions; namely,
 the $i$'th coordinate $(c\circ b)_i$ of the $n$-tuple $c\circ b$ is equal to $c_i(b_1,
 \ldots , b_n)$ where $c_i=c_i(x_1, \ldots , x_n)$ (we have substituted
 elements $b_i$ for $x_i$ in the function $c_i=c_i(x_1, \ldots , x_n)$).

By Theorem \ref{29Sep06}, each element $\s$ of the group $G= \O \G
\GL_n(K)^{op}$ is the {\em unique} product 
\begin{equation}\label{supr}
\s = \o_{1+a} \g_b\s_A, \;\; \o_{1+a} \in \O,\;  \g_b\in \G,\;
A\in \GL_n(K)^{op}.
\end{equation}
The product of two elements of $G$ is given by the rule
\begin{equation}\label{mul1}
\o_{1+a} \g_b\s_A \cdot \o_{1+a'} \g_{b'}\s_{A'}
=\o_{1+a+\g_b\s_A(a')}\, \g_{A^{-1} \s_A(b')\circ b}\, \s_{A'A}
\end{equation}
where $\s_A(b'):= (\s_A(b_1'), \ldots ,\s_A(b_n'))^t$ and
$\s_A(b')\circ b:= (\s_A(b_1')\circ b, \ldots , \s_A(b_n')\circ
b)$. This formula shows that the most sophisticated part of the
group $G$ is the group $\G$. To prove (\ref{mul1}), note first
that $\s_A\g_{b'}\s_A^{-1} = \g_{A^{-1}\s_A(b')}$, then
\begin{eqnarray*}
 \o_{1+a} \g_b\s_A \cdot \o_{1+a'} \g_{b'}\s_{A'}
&=& \o_{1+a}\cdot  \g_b\s_A \o_{1+a'}(\g_b\s_A)^{-1}\cdot
\g_b\cdot \s_A
 \g_{b'}\s_A^{-1} \cdot \s_A\s_{A'}\\
 &=& \o_{1+a} \o_{1+\g_b\s_A(a')}\cdot \g_b\g_{A^{-1}
 \s_A(b')}\cdot \s_{A'A}\\
&=& \o_{1+a+\g_b\s_A(a')}\, \g_{A^{-1} \s_A(b')\circ b}\, \s_{A'A}. \\
\end{eqnarray*}
We know how to find inverse elements for the group $\O$
($\o_{1+a}^{-1} = \o_{1-a}$) and for the group $ \GL_n(K)^{op}$
($\s_A^{-1} = \s_{A^{-1}}$). The inversion formula for elements
$\g_b$ of the group $\G$ is given explicitly by Theorem
\ref{i9Sep06}. So, one can find explicitly an element $b'$ such
that $\g_b^{-1} = \g_{b'}$ by applying Theorem \ref{i9Sep06}:
$b_i':= \g_b^{-1} (x_i)$. Now, one can write down the explicit
formula for the inverse of any element $\s =\o_{1+a} \g_b\s_A \in
G$, 
\begin{equation}\label{invabA}
(\o_{1+a} \g_b\s_A)^{-1}=\o_{1-\s_{A^{-1}}\, \g_{b'}(a)}
\g_{A\s_{A^{-1}}(b')}\, \s_{A^{-1}}.
\end{equation}
In more detail, by (\ref{mul1}), $(\o_{1+a}
\g_b\s_A)^{-1}=\s_{A^{-1}}\g_{b'}\o_{1-a}= \o_{1-\s_{A^{-1}}\,
\g_{b'}(a)} \g_{A\s_{A^{-1}}(b')}\, \s_{A^{-1}}$.
\begin{corollary}\label{G30Sep06}
Let $K$ be a reduced commutative ring with $\frac{1}{2}\in K$.
Then
\begin{enumerate}
\item $G_{\Z_2-gr} = \G \rtimes \GL_n(K)^{op}$.\item  $G= \O
\rtimes G_{\Z_2-gr}$. \item ${\rm Out}(\L_n)\simeq G_{\Z_2-gr}$.
\end{enumerate}
\end{corollary}

{\it Proof}. By (\ref{mul1}), $G':=\G \GL_n(K)^{op}=\G \rtimes
\GL_n(K)^{op}$ is the subgroup of $G$ such that $G= \O \rtimes G'$
(Theorem \ref{29Sep06}) and $G'\subseteq G_{\Z_2-gr}$. By
(\ref{1GImo}), $G_{\Z_2-gr}\cap \O = \{ e\}$, hence $G_{\Z_2-gr}=
G'$ and $G= \O \rtimes G'= \O \rtimes G_{\Z_2-gr}$. Then ${\rm
Out}(\L_n)\simeq \O \rtimes G_{\Z_2-gr}/\O \simeq G_{\Z_2-gr}$.
$\Box $

Let $K$ be a commutative ring. Consider the sets of {\em even} and
{\em odd} automorphisms of $\L_n$:
\begin{eqnarray*}
 G^{ev}&:=&\{ \s \in G\, | \, \s (x_i) \in \L_{n,1}+\Lnev , \forall i \}, \\
G^{od}&:=&\{ \s \in G\, | \, \s (x_i) \in \L_{n,1}+\Lnod , \forall
i \}.
\end{eqnarray*}
One can easily verify that $G^{od}$ is a subgroup of $G$. It is
not obvious from the outset that $G^{ev}$ is also a subgroup of
$G$.

\begin{lemma}\label{e30Sep06}
Let $K$ be a reduced commutative ring with $\frac{1}{2}\in K$.
Then
\begin{enumerate}
\item $G^{od} = G_{\Z_2-gr} =\G\rtimes \GL_n(K)^{op}$. \item
$G^{ev}= \O \rtimes \GL_n(K)^{op}$. \item $G^{od} \cap G^{ev} =
\GL_n(K)^{op}$. \item $G= G^{od} G^{ev} = G^{ev} G^{od}$.
\end{enumerate}
\end{lemma}

{\it Proof}. 3. $G^{od} \cap G^{ev} =\{ \s \in G\, | \, \s (x_i)
\in \L_{n,1}, \forall i\}= \GL_n(K)^{op}$.

4. This follows from statements 1 and 2 since $G= \O \G
\GL_n(K)^{op}$ (Theorem \ref{29Sep06}).

1 and 2. Recall that $G_{\Z_2-gr}= \G \GL_n(K)^{op}$ (Corollary
\ref{G30Sep06}) and $ G= \O \G  \GL_n(K)^{op}=\O G_{\Z_2-gr}$
(Theorem \ref{29Sep06}). Clearly, $G_{\Z_2-gr}\subseteq G^{od}$
and $\O \cap G^{od} = \{ e \}$, it follows that  $$
G^{od}=(G^{od}\cap \O )G_{\Z_2-gr}=G_{\Z_2-gr}.$$ Similarly, $\O
\GL_n(K)^{op}\subseteq G^{ev}$ and $\G \cap G^{ev} =\{ e \}$ give
the equality $\O \GL_n(K)^{op}= G^{ev}$: $G^{ev} = \O (\G \cap
G^{ev})\GL_n(K)^{op}= \O \GL_n(K)^{op}$. $\Box$

{\it Example}. For $n=2$, $G_2= \O_2\GL_2(K)^{op}= \{ \o_{1+a}
\s_A\, | \, a= \l_1x_1+\l_2x_2, \l_i\in K, A\in \GL_n(K)\}$.
 The group $\O_2$ is canonically isomorphic to $K^2$ via $
 \o_{1+\l_1x_1+\l_2x_2} \mapsto \l : = (\l_1, \l_2)^t$. Then
 $G_2\simeq
 K^2\GL_2(K)^{op} = \{ (\l , A)\, | \, \l \in K^2, A\in
 \GL_2(K)\}$ and
 $$ (\l , A) \cdot (\l', A')= (\l +A^t\l', A'A)\;\; \; {\rm and} \;\;\; (\l ,
 A)^{-1} = (-(A^t)^{-1} \l , A^{-1} ).$$

{\it Example}. For $n=3$, $G_3= \O_3\G_3\GL_3(K)^{op}$. The group
$\O_3$ is canonically isomorphic to $K^3$ via $
 \o_{1+\l_1x_1+\l_2x_2+\l_3x_3} \mapsto \l : = (\l_1, \l_2, \l_3)^t$.
 Similarly, the group $\G_3$ is canonically isomorphic to $K^3$
 via $\g_b= \g_{(x_1+\mu_1 \th , x_2+\mu_2 \th , x_3+\mu_3 \th
 )}\mapsto \mu = (\mu_1, \mu_2, \mu_3)^t$ where $\th :=
 x_1x_2x_3$.
  Then $G_3\simeq
 K^3K^3\GL_3(K)^{op} = \{ (\l ,\mu,  A)\, | \, \l , \mu  \in K^3, A\in
 \GL_3(K)\}$ and
\begin{eqnarray*}
 (\l ,\mu,  A) \cdot (\l',\mu',  A')&=& (\l +A^t\l',\mu +\det (A) A^{-1} \mu' ,A'A),  \\
(\l
 ,\mu ,
 A)^{-1} &=& (-(A^t)^{-1} \l , -\det (A^{-1})\cdot A\mu , A^{-1} ).
\end{eqnarray*}
For $n=2,3$, the group $U:= U_n= \O_n\G_n$ is abelian: $U_2=
\O_2\simeq K^2$ and $U_3= \O_3\G_3= \O_3U_3^3\simeq K^3\times
K^3$. The next result shows that these are the only cases where
$U$ is an abelian group.

{\bf Maximal abelian subgroups of $U$}.
\begin{theorem}\label{28Sep06}
Let $K$ be a reduced commutative ring with $\frac{1}{2}\in K$.
Then
\begin{enumerate}
\item The group $\O$ is a maximal abelian subgroup of $U$ if $n$
is even ($\O \supseteq U^n$).\item The group $\O U^n=\O\times U^n$
is a maximal abelian subgroup of $U$ if $n$ is odd ($\O \cap U^n
=\{ e\}$).
\end{enumerate}
\end{theorem}

{\it Proof}. Recall that $\O $ is an abelian subgroup of $U$
(Theorem \ref{m28Sep06}.(5)) and that $U= \O \G (= \O \rtimes \G
)$ is the semidirect product of the groups $\O $ and $\G$ (Theorem
\ref{29Sep06}.(1)). Suppose that $\O'$ is an abelian subgroup of
$U$ such that $\O \subseteq \O'$.  Then $\O'= \O (\G \cap \O')$.
 Each  element $\g_b\in \G\cap \O'$  must
commute with all the elements of $\O$: $\o_{1+a} \g_b = \g_b
\o_{1+a}= \o_{1+\g_b(a)}\g_b$ for all $a$ iff $\o_{1+ \g_b(a)-a} =
e$ iff $[\g_b(a) -a, x]=0$ for all $x\in \L_n$ (since $ \o_{1+a'}
(x) = x+[a',x]$, Lemma \ref{m28Sep06}.(3)) iff $\g_b (a) -a\in
Z(\L_n)$, the centre of $\L_n$, iff $\g_b(a) = a$ for all $a$ if
$n$ is even; and $\g_b(a) -a\in Kx_1\cdots x_n$ if $n$ is odd, iff
$\g_b = e$ if $n$ is even; and $\g_b\in U^n$ if $n$ is odd.  Since
the group $U^n$ is abelian and all elements of $\O$ commute with
elements of $U^n$ if $n$ is odd, the result follows. $\Box $


\section{The group $\G$, its subgroups, and the Inversion Formula}\label{SGPG}

In this section, the {\em inversion formula} (Theorem
\ref{i9Sep06}) is given for any automorphism $\s \in
\G\GL_n(K)^{op}$; the groups $\G_{\Z_s-gr}$, $s\geq 2$, are found
(Lemma \ref{4Oct06}); minimal sets of generators are given for the
groups $\G $ and $U$ (Theorem \ref{7Oct06}) and the  commutator
series are found for them; several important subgroups are
introduced: $\Phi$, $\Phi'$, $\Phi(i)$, $\CE_{n,i}$, $\CE_{n,i}'$,
$\CE_{n,i}''$.

{\bf The Jacobian and the inversion formula for automorphism}. Let
$K$ be a commutative ring and $\der_1:= \frac{\der}{\der x_1},
\ldots , \der_n:= \frac{\der}{\der x_n}$ be the  left skew partial
derivatives for $\L_n$. For each automorphism $\s \in
\G\GL_n(K)^{op}$, the matrix of  left skew partial derivatives
$\frac{\der \s}{\der x}:= (\frac{\der \s (x_i)}{\der x_j})$ is
called the {\bf Jacobian} matrix for the automorphism $\s$. Note
that the entries of the Jacobian matrix are even elements, hence
 central elements. The determinant $\CJ (\s ):= \det (\frac{\der
\s (x_i)}{\der x_j})$ is called the {\bf Jacobian} of $\s$. One
can easily verify that the `{\em chain rule}' holds for
automorphisms $\s, \tau \in \G\GL_n(K)^{op}$: 
\begin{equation}\label{char}
\frac{\der (\s \tau )}{\der x}= \s (\frac{\der \tau}{\der x})\cdot
\frac{\der \s}{\der x}
\end{equation}
where $\s (\frac{\der \tau}{\der x}):= (\s (\frac{\der \tau
(x_i)}{\der x_j}))$. By taking the determinant of both sides we
have the equality 
\begin{equation}\label{char1}
\CJ (\s \tau ) = \s (\CJ (\tau ))\, \CJ (\s ).
\end{equation}
Then, for each $\s \in \G\GL_n(K)^{op}$, 
\begin{equation}\label{Js1i}
\CJ (\s^{-1})= \s^{-1} (\CJ (\s )^{-1}).
\end{equation}

Let $\s \in \G\GL_n(K)^{op}$ and $x_1':= \s (x_1), \ldots , x_n':=
\s (x_n)$. The elements $x_1', \ldots , x_n'$ are another set of
canonical generators for  $\L_n$: $x_i'x_j'=-x_j'x_i'$ and
$x_i'^2=0$. The corresponding left  skew derivations $\der_1':=
\frac{\der}{\der x_1'}, \ldots ,\der_n':= \frac{\der}{\der x_n'}$
are equal to 
\begin{equation}\label{skderi}
\der_i'  := \frac{1}{\CJ(\s )} \, \det
 \begin{pmatrix}
  \frac{\der \s (x_1)}{\der x_1} & \cdots & \frac{\der \s (x_1)}{\der x_m} \\
  \vdots & \vdots & \vdots \\
\frac{\der }{\der x_1} & \cdots & \frac{\der }{\der x_m}\\
 \vdots & \vdots & \vdots \\
\frac{\der \s (x_n)}{\der x_1} & \cdots & \frac{\der \s (x_n)}{\der x_m} \\
\end{pmatrix}, \;\;\; i=1, \ldots , n,
\end{equation}
where we `drop' $\s (x_i)$ in the determinant $\CJ(\s ):= \det
(\frac{\der \s (x_i)}{\der x_j})$.

For each $i=1, \ldots , n$, let 
\begin{equation}\label{shPdad3}
\phi_i':= 1-x_i'\der_i': \L_n\ra \L_n,
\end{equation}
and (the order is important)
\begin{equation}\label{shPdad4}
\phi_\s := \phi_n'\phi_{n-1}'\cdots \phi_1'= (1-x_n'\der_n')
(1-x_{n-1}'\der_{n-1}') \cdots (1-x_1'\der_1'):\L_n\ra \L_n .
\end{equation}

The next theorem gives  the  inversion formula for automorphisms
of the group $\G\GL_n(K)^{op}$.
\begin{theorem}\label{i9Sep06}
{\rm (The Inversion Formula)} Let $K$ be a commutative ring,
 $\s \in \G\GL_n(K)^{op}$
 and $a\in \L_n(K)$. Then
 $$ \s^{-1}(a)=\sum_{\alpha \in \CB_n}\phi_\s
 (\der'^\alpha (a))x^\alpha  $$
 where $x^\alpha:= x_1^{\alpha_1}\cdots x_n^{\alpha_n}$ and
 $\der'^\alpha:= \der_n'^{\alpha_n} \cdots \der_1'^{\alpha_1}$.
\end{theorem}

{\it Proof}. By Theorem \ref{s9Sep06}.(1), $a=\sum_{\alpha \in
\CB_n} \phi_\s  (\der'^\alpha(a))x'^\alpha$ where $x'^\alpha :=
x_1'^{\alpha_1} \cdots x_n'^{\alpha_n}$.  Applying $\s^{-1}$
 we have the result
 $$ \s^{-1}(a)=\sum_{\alpha \in \CB_n }\phi_\s
 (\der'^\alpha (a))\s^{-1}(x'^\alpha ) =\sum_{\alpha \in \CB_n}\phi_\s
 (\der'^\alpha (a))x^\alpha . \;\;\; \Box $$

{\bf The abelian groups of units $E_n$ and $E_n'$}.  Let $K$ be a
commutative ring and $E_n$ be the group of units of the
commutative algebra $\Lnev$. So, $E_n= K^*+\L_{n, \geq 2}^{ev}$
where $K^*$ is the group of units of the ring $K$ and $\L_{n, \geq
2}^{ev}:= \gm^2\cap \Lnev=\oplus_{m=1}^{[\frac{n}{2}]}\L_{n,2m}$.
There is the natural descending chain of subgroups of $E_n$
determined by the $\gm$-adic filtration of the Grassmann algebra
$\L_n$: 
\begin{equation}\label{Ench}
E_n= E_{n,2}\supset E_{n,4}\supset \cdots \supset
E_{n,[\frac{n}{2}]}\supset E_{n,[\frac{n}{2}]+2}=K^*
\end{equation}
where $E_{n, 2m}:= K^* +\gm^{2m} \cap \Lnev = K^*
+\sum_{i=m}^{[\frac{n}{2}]}\L_{n, 2i}$.

Each element $e\in E_n$ is a unique sum $e= \l + e^+$ for some $\l
\in K^*$ and $e^+\in \gm^2\cap \Lnev$. The map
$$ v: E_n\ra 2\Z , \;\; e\mapsto v(e) = \max \{ 2m \, | \, e^+\in
E_{n, 2m}\},$$satisfies the following  properties: for $e,f\in
E_n$,
\begin{enumerate}
\item $v(ef) \geq \min \{ v(e), v(f)\}$,
  \item $v(e^{-1} )=
v(e)$.
\end{enumerate}

The group $E_n=K^*E_n'=K^*\times E_n'$ is the direct product of
its subgroups $K^*$ and $E_n':= 1+ \L_{n, \geq 2}^{ev}$ ($E_n=
K^*E_n'$ and $K^*\cap E_n'=\{ 1\}$). The chain (\ref{Ench})
induces the chain of subgroups in $E_n'$: 
\begin{equation}\label{cEns}
E_n'= E_{n,2}'\supset E_{n,4}'\supset \cdots \supset
E_{n,[\frac{n}{2}]}'\supset E_{n,[\frac{n}{2}]+2}=\{ 1\}
\end{equation}
where $E_{n, 2m}':= E_n'\cap E_{n,2m}= 1+ \gm^{2m} \cap \Lnev =
1+\sum_{i=m}^{[\frac{n}{2}]}\L_{n, 2i}$. If $K$ is a  reduced
commutative ring then, by (\ref{smi=mi}), 
\begin{equation}\label{sEn2m}
\s (E_{n,2m})= E_{n,2m}\;\; {\rm and}\;\;\s (E_{n,2m}')=
E_{n,2m}'\;\; {\rm for\; all}\;\; \s \in G_{\Z_2-gr}, \; m\geq 1.
\end{equation}
It follows that (where $K$ is reduced) 
\begin{equation}\label{1sEn2m}
v(\s (e))= v (e)  \;\; {\rm for\; all}\;\; e\in E_n,\;\;\s \in
G_{\Z_2-gr}.
\end{equation}

\begin{lemma}\label{g4Oct06}
Let $K$ be a  commutative ring. Then
\begin{enumerate}
\item each element $\s \in \G$ is a (unique) product $\s = \cdots
\s_{b_7}\s_{b_5} \s_{b_3}$ for unique elements $b_i= (b_{i1},
\ldots , b_{in})\in \L_{n,i}^n$ (see Corollary
\ref{c11Sep06}).\item
 The elements $\{ \g_{i, \l x^\alpha}\, | \, 1\leq i\leq n , \l
 \in K, \alpha \in \CB_n, 3\leq |\alpha |\;  {\rm is \; odd}\}$ are
 generators for the group $\G$ where $\g_{i, \l x^\alpha}(x_i):= x_i+\l x^\alpha$ and $\g_{i, \l
 x^\alpha}(x_j):=x_j$ for $i\neq j$.
\end{enumerate}
\end{lemma}

{\it Proof}. 1. This follows from Theorem \ref{11Sep06}.(2).

2. This statement follows from statement 1 and Theorem
\ref{11Sep06}.(2).  $\Box $

{\bf The dimension of the algebraic group $\G$}. Let $K$ be a
commutative ring and $\L_n = \L_n (K)$. A typical element of $\G$
is an automorphism $x_1\mapsto x_1+a_1, \ldots , x_n\mapsto
x_n+a_n$ where all $a_i\in \L_{n, \geq 3}^{od}$. The group $\G$ is
a unipotent algebraic group over $K$ where the coefficients of the
elements $a_i$ are the affine coordinates for the algebraic group
$\G$ over $K$, and the algebra of (regular) functions $\OO (\G )$
of the algebraic group $\G$ is a polynomial algebra in
\begin{equation}\label{dimG}
\dim (\G ) = n(2^{n-1}-n)
\end{equation}
variables since
$$ \dim (\G )= {\rm rk}_K((\L_{n, \geq 3}^{od})^n) = n
\sum_{i=1}^{[\frac{n-1}{2}]}{n\choose 2i+1}= n (
\sum_{i=0}^{[\frac{n-1}{2}]}{n\choose 2i+1}-n)=n(2^{n-1}-n).$$ If
$K$ is a  field then $\dim (\G )$ is the usual  dimension of the
algebraic group $\G $.

 {\bf A noncommutative
analogue of the Taylor expansion}.

\begin{theorem}\label{T2Oct06}
({\rm An analogue of the Taylor expansion}) Let $K$ be a
commutative ring, $f= f(x_1, \ldots , x_n)= \sum x^\alpha
\l_\alpha \in \L_n$ where the coefficients $\l_\alpha \in K$ of
$f$ are written on the right, and $\s \in \G$. Let $x_1':= \s
(x_1)= x_1+a_1, \ldots , x_n':=\s (x_n)= x_n+a_n$ where $a_1,
\ldots , a_n$ are odd elements of $\gm^3$. Then
$$ (\s (f))(x_1, \ldots , x_n) = f(x_1+a_n, \ldots , x_n+a_n)=
\sum_{ \alpha \in \CB_n} a^\alpha \der^\alpha (f)$$ where $
a^\alpha := a_1^{\alpha_1}\cdots a_n^{\alpha_n}$ and $\der^\alpha
= \der_n^{\alpha_n}\cdots \der_1^{\alpha_1}$, $\der_i:=
\frac{\der}{\der x_i}$ are the left partial skew derivatives of
$\L_n$.
\end{theorem}

{\it Proof}. It suffices to prove the statement for $f= x_1 \cdots
x_m\l$ where $\l\in K$. Then
\begin{eqnarray*}
\s (f)& = & (x_1+a_1)\cdots (x_m+a_m)\l=\sum_{i_1<\cdots < i_s}
x_1\cdots a_{i_1} \cdots a_{i_2} \cdots a_{i_s}\cdots
x_m\l\\
&=&\sum_{i_1<\cdots < i_s} a_{i_1} \cdots a_{i_s} \der_{i_s}
\cdots \der_{i_1} (f) =\sum_{ \alpha \in \CB_n} a^\alpha
\der^\alpha (f).\;\;\; \Box
\end{eqnarray*}

{\bf The groups $\G (s)$ and $\O (s)$}.  The next Lemma introduces
subgroups determined by even subgroups of $\Z$ (the subgroups of
type $2m\Z$).

\begin{lemma}\label{n2Oct06}
Let $K$ be a commutative ring and $\L_n = \L_n(K)$.
\begin{enumerate}
\item For each even number $s=2,4, \ldots , 2[\frac{n}{2}]$, the
subset of $\G$, $\G (s) := \{ \s \in \G \, | \, \s (x_i) \in x_i
+\sum_{j\geq 1} \L_{n,1+js}$ for all $i\}$ is a subgroup of
$\G$.\item If $s|s'$ ($s$ divides $s'$) then $\G (s') \subseteq \G
(s)$.
\end{enumerate}
\end{lemma}

{\it Proof}. 1. It is easy to see that the set $\G (s)$ is closed
under multiplication and that it contains the identity. The fact
that the set $\G (s)$ is closed under the operation of  taking
inverse is a consequence of repeated use of Theorem \ref{T2Oct06}.
Let $\s \in \G (s)$ and $x_i':= \s (x_i) = x_i-a_i$ for some odd
element  $a_i\in \sum_{j\geq 1} \L_{n, 1+js}$. Now,
\begin{eqnarray*}
 x_i &=& x_i'+a_i(x_1, \ldots , x_n)= x_i'+a_i(x_1'+a_1(x), \ldots , x_n'+a_n(x)) \\
 &=& x_i'+a_i(x_1', \ldots , x_n')+ \sum_{0\neq \alpha \in \CB_n} a^{\alpha} (x) \der^{\alpha} (a_i)
 (x')\\
 &=& x_i'+a_i(x')+ \sum_{0\neq \alpha \in \CB_n} a^{\alpha} (x_1'+a_1(x), \ldots , x_n'+a_n(x)) \der^{\alpha} (a_i) (x') \\
 &=& x_i'+a_i(x')+ \sum_{0\neq \alpha \in \CB_n} a^{\alpha} (x') \der^{\alpha} (a_i)
 (x')+\sum_{0\neq \alpha \in \CB_n} (\sum_{0\neq \beta \in\CB_n}a^{\beta} (x) \der^{\beta} (a^\alpha ) (x'))\der^{\alpha} (a_i)
 (x')\\
 &=&\cdots ,
\end{eqnarray*}
keep going making substitutions $x_i= x_i'+a_i$ and then using
Theorem \ref{T2Oct06} we get the result (in no more than
$[\frac{n}{s}]+1$ steps).

2. This is obvious.  $\Box $

\begin{lemma}\label{s4Oct06}
Let $K$ be a  commutative ring, and $s=2,4, \ldots ,
2[\frac{n}{2}]$. Then
\begin{enumerate}
\item each element $\s \in \G (s)$ is a (unique) product $s
=\cdots \s_{b_{1+3s}}\s_{b_{1+2s}}\s_{b_{1+s}}$ for unique
elements $b_i= (b_{i,1}, \ldots , b_{i,n})\in \L_{n,i}^n$ (see
Lemma \ref{g4Oct06}).\item The elements $\{ \g_{i, \l x^\alpha }\,
| \, 1\leq i\leq n , \l \in K^*, \alpha \in \CB_n$ such that
$\alpha \in 1+s\N\}$ are generators for the group $\G (s)$ (see
Lemma \ref{g4Oct06}).
\end{enumerate}
\end{lemma}

{\it Proof}. These statements follow from Lemma \ref{g4Oct06}.
$\Box $

For each  odd number $s$ such that $1\leq s\leq n$, the set
\begin{equation}\label{Osdef}
\O (s):= \{ \o_{1+a} \, | \, a\in \sum_{j\geq 1} \L_{n, js}\}= \{
\o_{1+a} \, | \, a\in \sum_{1\leq j\; {\rm is\, odd}} \L_{n, js}\}
\end{equation}
is a subgroup of $\O$ (Lemma \ref{m28Sep06}.(4)). Note that $ \O
(1)=\O $ and $\O (n) = \{ e\}$.

{\bf The groups $G_{\Z_s-gr}$}. Recall that $G_{\Z_s-gr}$ is the
group of all $K$-automorphisms of the Grassmann algebra  $\L_n(K)$
that respect its $\Z_s$-grading ($\s\in G_{\Z_s-gr}$ iff $\s
(\L_{n, is})= \L_{n,is}$, $i\geq 0$). The next result describes
the groups $G_{\Z_s-gr}$.

\begin{lemma}\label{4Oct06}
Let $K$ be a reduced commutative ring with $\frac{1}{2}\in K$.
Then
\begin{enumerate}
\item if $s$ is even then $G_{\Z_s-gr} = \G (s) \GL_n(K)^{op}= \G
(s) \rtimes \GL_n(K)^{op}$ and $\G (s) = \G \cap
G_{\Z_s-gr}$;\item if $s\geq 3$ is odd then $G_{\Z_s-gr} = \O (s)
\GL_n(K)^{op}= \O (s) \rtimes \GL_n(K)^{op}$ and $\O (s) = \O \cap
G_{\Z_s-gr}= \{ \o_{1+a} \, | \, a\in \sum_{1\leq j\; {\rm is\,
odd}} \L_{n, sj}\}$.
\end{enumerate}
\end{lemma}

{\it Proof}. 1. The number $s$ is even,  hence $G_{\Z_s-gr}
\subseteq G_{\Z_2-gr} = \G \GL_n(K)^{op}$ (Lemma
\ref{G30Sep06}.(1)). Then it follows from the inclusion
$\GL_n(K)^{op}\subseteq G_{\Z_s-gr}$ that $G_{\Z_s-gr} = (\G \cap
G_{\Z_s-gr}) \GL_n(K)^{op}$. So, to finish the proof of statement
1 it suffices to show that $\G (s) = \G \cap G_{\Z_s-gr}$. The
inclusion $ \G (s) \subseteq \G \cap G_{\Z_s-gr}$ is obvious. If $
e\neq \g \in \G\cap G_{\Z_s-gr}$ then $\g = \cdots
\s_{b_{2k+3}}\s_{b_{2k+1}}$ (Lemma \ref{g4Oct06}) where $b_{2k+1}=
(b_{2k+1,1}, \ldots , b_{2k+1,n})\neq 0$, $\g (x_i) = x_i+b_{2k+1,
i}+\cdots $ for all $i$. Hence $2k+1\in 1+ s\Z$, i.e.
$\s_{b_{2k+1}}\in \G (s)$. Applying the same argument to $\g
\s_{b_{2k+1}}^{-1} = \cdots \s_{b_{2k+3}} \in \G \cap G_{\Z_s-gr}$
and using induction on $k$ we see that all $\s_{b_{2k+i}}$ in the
product for $\g$ belong to $\G (s)$. Therefore, $\G (s)= \G \cap
G_{\Z_s-gr}$.

2. By Lemma \ref{m28Sep06}.(3), $ \O \cap G_{\Z_s-gr}= \{\o_{1+a}
\, | \, a\in \sum_{1\leq j\; {\rm is\, odd}} \L_{n, sj}\}=\O (s)$.
Considering the action of automorphisms from the intersection
 $G_{\Z_s-gr}\cap \O \G$ on the generators $x_1, \ldots , x_n$
 (with help of Corollary \ref{c11Sep06}) it is easy to show that
 $G_{\Z_s-gr}\cap \O \G = G_{\Z_s-gr}\cap \O$. Recall that $G= \O
 \G \GL_n(K)^{op}$ and $\GL_n(K)^{op}\subseteq G_{\Z_s-gr}$ hence
 $$ G_{\Z_s-gr}= (G_{\Z_s-gr}\cap \O \G )\GL_n(K)^{op}= \O (s)
 \GL_n(K)^{op}= \O (s) \rtimes \GL_n(K)^{op}. \;\;\; \Box $$

{\bf The groups $\Phi$, $\Phi (i)$,  and $\Phi'$}. Let $K$ be a
commutative ring.  Clearly, 
\begin{equation}\label{Gaibi}
\G = \{ \s : x_i\mapsto x_i(1+a_i)+b_i, \; i=1, \ldots , n\}
\end{equation}
where $a_i, b_i\in K\lfloor x_1, \ldots , \widehat{x_i}, \ldots ,
x_n\rfloor$, $a_i\in \Lnev\cap \gm^2$ and $b_i\in \Lnod \cap
\gm^3$. Consider the subset $\Phi$ of $\G$ where all $b_i=0$,
\begin{equation}\label{}
\Phi := \{ \s : x_i\mapsto x_i(1+a_i), \; i=1, \ldots , n\}.
\end{equation}
The set $\Phi$ can be characterized as
$$ \Phi =\{ \s \in \G \, | \, \s (x_1) \in (x_1) , \ldots , \s
(x_n)\in (x_n)\}.$$Then it is obvious that $\Phi \Phi \subseteq
\Phi$ and ${\rm id}_{\L_n} \in \Phi$.

\begin{lemma}\label{p2Oct06}
Let $K$ be a commutative ring. Then $\Phi =\{ \s \in \G \, | \, \s
((x_1)) = (x_1) , \ldots , \s ((x_n))= (x_n)\}$ is a subgroup of
$\G $.
\end{lemma}

{\it Proof}. It remains to show that, for each $\s \in \Phi$,
$\s^{-1} (x_i) \in (x_i)$ for all $i$. The equation $x_i':= \s
(x_i) := x_i(1-a_i)$ can be written as $x_i= x_i'+x_ia_i(x_1,
\ldots , x_n)$ (we have changed the sign of the $a_i$ for
computational reason). Our aim is to show that $x_i=
x_i'(1+a_i')$, then the result will follow as $x_i'\in (x_i)$. We
use in turn, first, the substitution $x_i= x_i'+ x_ia(x)$ and then
the Taylor expansion (Theorem \ref{T2Oct06}). After repeating
these no more than $n+1$ times we will get the result (since all
elements are nilpotent and any product of $n+1$ of them is zero).

\begin{eqnarray*}
 x_i&=& x_i'+x_ia_i(x)= x_i'+ (x_i'+ x_ia_i(x))a_i(x_1'+ x_1a_1(x), \ldots , x_n'+x_na_n(x)) \\
 &=& x_i'+ (x_i'+ x_ia_i(x))(a_i(x_1', \ldots , x_n')+\sum_{0\neq
 \alpha \in \CB_n} (xa(x))^\alpha \der^\alpha (a_i) (x_1', \ldots,
 x_n'))=\cdots . \; \Box
\end{eqnarray*}
The group $\Phi$ is the solutions of the polynomial equations in
coefficients of the elements $a_i$ and $b_j$: $b_1=0, \ldots ,
b_n=0$. So, $\Phi$ is a closed subgroup of $\G$ with respect to
the Zariski topology. The group $ \Phi$ is an  algebraic group,
the algebra of functions on $\Phi$ is a polynomial algebra over
$K$ in $n\cdot {\rm rk}_K(\L_{n-1, \geq 2}^{ev})=n(2^{n-2}-1)$
variables. So, the algebraic group $\Phi$ is affine and
\begin{equation}\label{dimPhi1}
\dim (\Phi ) =n(2^{n-2}-1).
\end{equation}

Note, that, in general,  the set $\{\s \in \G  \, | \, \s (x_1) =
x_1+b_1, \ldots , \s (x_n)= x_n+b_n\}$ is not a subgroup of $\G$
where each $b_i\in K\lfloor x_1, \ldots , \widehat{x_i}, \ldots ,
x_n\rfloor\cap \gm^3$ is odd.

For each $i=1, \ldots , n$, let
$$ \Phi (i) := \{ \s \in \G \, | \, \s (x_i) \in (x_i)\} = \{ \s
\in G\, | \, \s (x_i)= x_i(1+a_i)\}$$ where $a_i\in K\lfloor x_1,
\ldots , \widehat{x_i}, \ldots , x_n\rfloor^{ev}_{\geq 2}$.
Clearly, $\Phi (i) \Phi (i) \subseteq \Phi (i)$ and the set $\Phi
(i)$ contains the identity map.

\begin{lemma}\label{i3Oct06}
Let $K$ be a commutative  ring. Then $\Phi (i)=\{ \s \in \G \, |
\, \s ((x_i))= (x_i)\}$ is a subgroup of $\G$.
\end{lemma}

{\it Proof}. It remains to prove that, given $\s \in \Phi (i)$,
$\s^{-1} \in \Phi (i)$. Repeat word for word the proof of Lemma
\ref{p2Oct06}. Note that if $K$ is a field  the result is obvious
since $\dim_K ((x_i))= \dim_K(\s ((x_i)))$: it follows from $\s
((x_i))\subseteq (x_i)$ that $\s ((x_i))= (x_i)$, hence $\s^{-1}
((x_i)) = (x_i)$, as required.  $\Box $

It is obvious that $\Phi = \cap_{i=1}^n \Phi (i)$, and
\begin{eqnarray*}
\Phi (i_1, \ldots , i_s):= \cap_{\nu =1} ^s \Phi (i_\nu )&=& \{ \s
\in \G \, | \, \s (x_{i_1})\in (x_{i_1}), \ldots , \s
(x_{i_s})\in (x_{i_s})\}\\
&=& \{ \s \in \G \, | \, \s ((x_{i_1}))= (x_{i_1}), \ldots , \s
((x_{i_s}))= (x_{i_s})\}
\end{eqnarray*}
 is a subgroup of $\G$. The set
$$ \Phi':= \{ \s \in G\, | \, \s (x_1)\in (x_1), \ldots \s
(x_n)\in (x_n)\}$$ is a group, as the next Lemma shows. Let
$\mathbb{T}^n$ be the subgroup of all the diagonal matrices of
$\GL_n(K)^{op}$.

\begin{lemma}\label{5Oct06}
Let $K$ be a reduced commutative ring with $\frac{1}{2}\in K$.
Then $\Phi'=(\O \rtimes \Phi ) \rtimes \mathbb{T}^n$ and $\Phi'=\{
\s \in G\, | \, \s ((x_1))= (x_1), \ldots , \s ((x_n))= (x_n)\}$.
\end{lemma}

{\it Proof}. It is obvious that $\O \subseteq \Phi'$,
 $\mathbb{T}^n\subseteq \Phi'$, and $\Phi'\cap \G \GL_n(K)^{op}=
\Phi\mathbb{T}^n$. Since $G= \O \G \GL_n(K)^{op}$ (Theorem
\ref{29Sep06}.(3)),
$$ \Phi'=\O (\Phi'\cap \G  \GL_n(K)^{op})= \O \Phi \mathbb{T}^n =(\O \rtimes \Phi ) \rtimes
 \mathbb{T}^n.$$ Then by Lemma \ref{p2Oct06}, $\Phi'=\{ \s \in G\, | \,
\s ((x_1))= (x_1), \ldots , \s ((x_n))= (x_n)\}$.  $\Box $

{\bf The groups $\CE_{n,i}$ and its subgroups}. Let $K$ be a
commutative ring.  For each $i=1, \ldots , n$, the set
$$\CE_{n,i} := \{ \g \in \G \, | \, \g (x_j)= x_j, \forall j\neq
i\}$$ is a subgroup of $\G$.  
\begin{equation}\label{CEni}
\CE_{n,i} =\{ \g_{1+a, b} : x_i\mapsto x_i (1+a) + b , \;
x_j\mapsto x_j, \forall j\neq i\}
\end{equation}
where $a\in  K\lfloor x_1, \ldots , \widehat{x_i}, \ldots ,
x_n\rfloor^{ev}_{\geq 2}$ and $b\in  K\lfloor x_1, \ldots ,
\widehat{x_i}, \ldots , x_n\rfloor^{od}_{\geq 3}$, and
\begin{eqnarray*}
 \g_{1+a, b}\g_{1+a', b'}&=&\g_{(1+a)(1+a'), b(1+a')+b'}, \;\; \g_{1+a, b}^{-1} = \g_{(1+a)^{-1}, -(1+a)^{-1}b},  \\
 \g_{1+a, b}\g_{1+a', b'}\g_{1+a, b}^{-1}&=&\g_{1+a', (1+a)^{-1}(ba'+b')}, \;\; \;\;\;\; \; \g_{1+a, b} = \g_{1, (1+a)^{-1}b} \g_{1+a, 0},
\end{eqnarray*}
Below, the equality (\ref{g2ab}) explains importance of these
small subgroups.  So, $\CE_{n,i}'= \{ \g_{1+a, 0}\}$ and
$\CE_{n,i}'':= \{ \g_{1, b}\}$ are abelian subgroups of
$\CE_{n,i}$ such that $\CE_{n,i}'\cap \CE_{n,i}''=\{ e\}$,
$\CE_{n,i}=\CE_{n,i}''\CE_{n,i}'$, and $\CE_{n,i}''$ is a  normal
subgroup of $\CE_{n,i}$ since 
\begin{equation}\label{g1ab}
\g_{1+a, 0} \g_{1,b} \g_{1+a, 0}^{-1} = \g_{1, (1+a)^{-1} b}.
\end{equation}
 Therefore,
$\CE_{n,i}=\CE_{n,i}''\rtimes \CE_{n,i}'$. Clearly,
$\CE_{n,i}'=\CE_{n,i}\cap \Phi$.

Let $E_{n,\hi}'$ be the group of units $E_{n-1}'$ in the case  of
the Grassmann algebra $K\lfloor x_1, \ldots , \widehat{x_i},
\ldots , x_n\rfloor$, i.e.
$$ E_{n,\hi }':= \{ 1+a\, | \, a\in  K\lfloor x_1, \ldots , \widehat{x_i}, \ldots
, x_n\rfloor^{ev}_{\geq 2}\}.$$

\begin{lemma}\label{E6Oct06}
Let $K$ be a commutative ring. Then
\begin{enumerate}
\item $\CE_{n,i} = \CE_{n,i}''\rtimes \CE_{n,i}'$.\item The map
$E_{n,\hi}'\ra \CE_{n,i}'$, $ 1+a\mapsto \g_{1+a, 0}$, is a group
isomorphism. \item The map $ K\lfloor x_1, \ldots , \widehat{x_i},
\ldots , x_n\rfloor^{od}_{\geq 3} \ra \CE_{n,i}''$, $ b\mapsto
\g_{1,b}$, is a group  isomorphism. \item $\CE_{n,i} = \langle
\g_{1+\l x^\alpha, 0}, \g _{1, \l x_jx_kx_l}\, | \, \l \in K,
|\alpha |=2 $, $j<k<l, \alpha \cup \{ j,k,l\} \subseteq \{ 1,
\ldots , \widehat{i}, $ $\ldots , n\} \rangle $.
\end{enumerate}
\end{lemma}

{\it Proof}. Statements 1--3 are obvious. Statement 4 follows from
statement 1 and the two facts: $(i)$ $\{ \g_{1+\l x^\alpha , 0}\}$
are generators for the group $\CE_{n,i}'$ since $\{ 1+\l x^\alpha
\}$ are generators for the group $E_{n,\hi}'$, and $(ii)$
\begin{equation}\label{g11ab}
\g_{1+a, 0} \g_{1,b} \g_{1+a, 0}^{-1}\g_{1,b}^{-1} = \g_{1,
(1+a)^{-1} b - b}= \g_{1, -ab +a^2b-a^3b+\cdots}= \g_{1,
-ab}\g_{1,a^2b-a^3b+\cdots} . \;\; \;\;\; \Box
\end{equation}
 For each $j\geq 2$, let $\{\CE_{n,i}^j:= \CE_{n,i} \cap U^j\}$ be
 the induced (descending) filtration on the group $\CE_{n,i}$.
 Each subgroup $\CE_{n,i}^j= \{ \s \in \CE_{n,i}\, | \, (\s -1)
 (\gm ) \subseteq \gm^j\}$ is a normal subgroup of $\CE_{n,i}$. By Lemma
 \ref{g4Oct06} and Theorem \ref{11Sep06}.(2), the group $\G$ is a
 {\em finite} product
\begin{equation}\label{g2ab}
\G = \cdots \prod_{i=1}^n \CE_{n,i}^{[2m+1]}\cdots \prod_{i=1}^n
\CE_{n,i}^{[5]}\cdot\prod_{i=1}^n \CE_{n,i}^{[3]},
\end{equation}
where $$\CE_{n,i}^{[2m+1]}:= \{ \g_{1+a, b}\in \CE_{n,i}\, | \,
a\in K\lfloor x_1, \ldots , \widehat{x_i}, \ldots ,
x_n\rfloor_{2m}, \; b\in K\lfloor x_1, \ldots , \widehat{x_i},
\ldots , x_n\rfloor_{2m+1}\}.$$ Clearly,
$\CE_{n,i}^{[2m+1]}\subseteq \CE_{n,i}^{2m+1}$.

 {\bf Minimal sets of generators for the groups $\G$,
$U$, and $\Phi$}. For each $i=1, \ldots , n$; $\l \in K$, and
$\alpha \subseteq \{ 1, \ldots , n\}$, $3\leq |\alpha | $ is odd,
let us consider the automorphism of $\G$,
$$\s_{i, \l x^\alpha} : x_i\mapsto x_i+\l x^\alpha, \; x_j\mapsto
x_j, \; \forall j\neq i.$$ Then  
\begin{equation}\label{sgco}
\s_{i, \l x_ix^\alpha}^{-1}=\s_{i, -\l x_ix^\alpha}.
\end{equation}

For two elements $a$ and $b$ of a group $A$, the {\em group
commutator} of the elements $a$ and $b$ is defined as $[a,b]:=
aba^{-1} b^{-1}\in A$. A direct (rather lengthy) calculation shows
that 
\begin{equation}\label{gcom1}
[\s_{i, \l x_ix_jx^\alpha}, \s_{j, \mu x_jx^\beta }]= \s_{i, -\l
\mu x_ix_jx^\beta x^\alpha }
\end{equation}
for all $\l ,\mu \in K$; $i\neq j$; $\alpha$ and $\beta$ are
subsets of $ \{ 1, \ldots , n\} \backslash \{ i,j\}$ such that
$\alpha \cap \beta = \emptyset$, $\alpha$ is  odd and $\beta$ is
even, $|\alpha | \geq 1$ and $|\beta | \geq 2$. Similarly,
\begin{equation}\label{gcom2}
[\s_{i, \l x_ix^\alpha}, \s_{i, \mu x^\beta }]= \s_{i, -\l \mu
x^\beta x^\alpha }
\end{equation}
for all $\l ,\mu \in K$; $i=1, \ldots , n$;  and  $\alpha , \beta
\subseteq \{ 1, \ldots , n\} \backslash \{ i\}$ such that $\alpha
\cap \beta = \emptyset$,  $\alpha$ is even and $\beta $ is  odd,
$|\alpha | \geq 2$ and $|\beta | \geq 3$.

For a group $A$, let us consider its series of commutators:
$$ A^{(0)}:= A,\; A^{(i)}:= [A, A^{(i-1)}], \; i\geq 1. $$
For each $i=1, \ldots , n$; $\l \in K$; and $ j <k<l$, let us
consider the automorphism $\s_{i, \l x_jx_kx_l}\in \G$:
$x_i\mapsto x_i+\l x_jx_kx_l$, $x_m\mapsto x_m$, for all $m\neq
i$. Then
$$\s_{i, \l x_jx_kx_l}\s_{i, \mu x_jx_kx_l}= \s_{i, (\l +\mu
)x_jx_kx_l}, \;\;\;\; \s_{i, \l x_jx_kx_l}^{-1} = \s_{i, -\l
x_jx_kx_l}^{-1}.$$ So, the set $\{ \s_{i, \l x_jx_kx_l} \, | \, \l
\in K\}$ is  isomorphic to the additive abelian group $K$, $
\s_{i, \l x_jx_kx_l}\mapsto \l$.

\begin{theorem}\label{7Oct06}
Let $K$ be a commutative ring in statements 1 and 2; and let $K$
be a reduced commutative ring with $\frac{1}{2}\in K$ in
statements 3 and 4. Then
\begin{enumerate}
\item the group $\G $ is generated by all the automorphisms
$\s_{i, \l x_jx_kx_l}$, i.e. $\G = \langle \s_{i, \l x_jx_kx_l}\,
| \, i=1, \ldots , n; \l \in K;  j <k<l\rangle$. The  subgroups
$\{ \s_{i, \l x_jx_kx_l}\}_{\l \in K}$ of $\G$ form a minimal set
of generators for $\G $.  \item $\G^{(i)}= \G^{2i+3} := \{ \s \in
\G \, | \, (\s -1)(\gm ) \subseteq \gm^{2i+3} \}$, $i\geq 0$.
\item The group $U$ is generated by all the automorphisms $\s_{i,
\l x_jx_kx_l}$ and all the automorphisms $\o_{1+\l x_i}$, i.e. $U
= \langle \s_{i, \l x_jx_kx_l}, \o_{1+\l x_i}\, | \, , i=1, \ldots
, n; \l \in K;  j <k<l\rangle$. The subgroups $\{ \s_{i, \l
x_jx_kx_l}\}_{\l \in K}$, $\{ \o_{1+\l x_i}\}_{\l \in K}$ of $U$
form a minimal set of generators for $U$. \item $U^{(i)}= U^{i+2}
:= \{ \s \in U \, | \, (\s -1)(\gm ) \subseteq \gm^{i+2} \}$,
$i\geq 0$.
\end{enumerate}
\end{theorem}

{\it Proof}. $1$. In a view of (\ref{g2ab}) and Lemma
\ref{E6Oct06}.(4), it suffices to show that each automorphism
$\g_{1+\l x^\alpha , 0}$ from Lemma \ref{E6Oct06}.(4) is a product
of the generators from statement 1. In the $\s$-notation, the
automorphism $\g_{1+\l x^\alpha , 0}$ is of the form
\begin{equation}\label{1g3ab}
\s := \s_{i_1, \l x_{i_1} x_{i_2} (x_{i_3} x_{i_4})\cdots
(x_{i_{2m-1}} x_{i_{2m}})(x_{i_{2m+1}} x_{i_{2m+2}})x_{i_{2m+3}}}
\end{equation}
for some {\em distinct} elements $i_1, i_2, \ldots ,  i_{2m+3}$, $
m\geq 0$, and $\l \in K$. The result is obvious for $m=0$. So, let
$m\geq 1$. Then, applying (\ref{gcom1}) $m$ times we have
\begin{equation}\label{g3ab}
\s = [\ldots [ \s_{i_1, \l x_{i_1} x_{i_2} x_{i_{2m+3}}}, \s_{i_2,
- x_{i_2} x_{i_{2m+1}} x_{i_{2m+2}}}], \s_{i_2, - x_{i_2}
x_{i_{2m-1}} x_{i_{2m}}}], \ldots ], \s_{i_2, - x_{i_2} x_{i_3}
x_{i_4}}].
\end{equation}
The claim that the `one dimensional' abelian subgroups $\{ \s_{i,
\l x_jx_kx_l}\}$ of $\G$ form a minimal set of generators is
obvious due to the isomorphism in Theorem \ref{11Sep06}.(2) in the
case $m=1$ there.

2. By Theorem \ref{11Sep06}.(2), $\G^{(m)}\subseteq \G^{2m+3}$ for
all $m\geq 0$. Clearly, $\G^{(n)}= \{ e\} = \G^{2n+3}$. Now, using
downward induction on $m$ (starting with $m=n$), in a view of
Theorem \ref{11Sep06}.(2) and (\ref{g3ab}), in order to prove the
equality $\G^{(m-1)}= \G^{2m+1}$, it suffices to show that each
automorphism of the type
$$\s = \s_{i, \l x_{i_1} x_{i_2}x_{i_3} (x_{i_4}x_{i_5}) \cdots (x_{i_{2m-2}}x_{i_{2m-1}})
(x_{i_{2m}}x_{i_{2m+1}})}$$ (where the elements $i, i_1, i_2,
\ldots , i_{2m+1}$ are {\em distinct},  and $\l \in K$) can be
expressed as $(m-1)$-commutator of the generators from statement 1
(i.e. $m-1$ brackets are involved). Below is such a presentation
(apply (\ref{gcom2})) 
\begin{equation}\label{g4ab}
\s = [ \s_{i, -\l x_ix_{i_{2m}} x_{i_{2m+1}}}, [\s_{i, -
x_ix_{i_{2m-2}} x_{i_{2m-1}}}[\ldots [ \s_{i, - x_ix_{i_4}
x_{i_5}}, \s_{i,
 x_{i_1}x_{i_2} x_{i_3}}]\ldots ].
\end{equation}

3. By Theorem \ref{11Sep06} and statement 1, it suffices to prove
that any inner automorphism $\o_{1+\l x_{i_1} x_{i_2} \cdots
x_{i_{2m+1}}}$, $i_1<\cdots <i_{2m+1}$,  $m\geq 0$, $\l \in K$, is
a product of the generators from statement 3. For any automorphism
$\s \in G$ and an odd element $a\in \L_n$ (Lemma
\ref{m28Sep06}.(4)): 
\begin{equation}\label{g5ab}
[\s , \o_{1+a}]=\s \o_{1+a} \s^{-1} \o_{1+a}^{-1} = \o_{1+ \s (a)
-a}.
\end{equation}
Applying this formula $m$ times we have 
\begin{equation}\label{g6ab}
\o_{1+ \l x_{i_1} x_{i_2} \cdots x_{i_{2m+1}}} = [ \s_{i_1,
x_{i_1} x_{i_2} x_{i_3}}, [\s_{i_1,  x_{i_1}x_{i_4}
x_{i_5}}[\ldots [ \s_{i_1, x_{i_1}x_{i_{2m}} x_{i_{2m+1}}}, \o_{1+
\l x_{i_1}}]\ldots ].
\end{equation}
The statement that generators in statement 3 are minimal follows
from the isomorphisms in Theorem \ref{11Sep06}.(2,3) for $m=1$
there.

$4$. By Theorem \ref{11Sep06}.(2,3), $U^{(m)}\subseteq U^{m+2}$,
$m\geq 0$. Clearly, $U^{(n)} = \{ e\} = U^{n+2}$. Now, using
downward induction on $m$ (starting with $m=n$), in a view of
Theorem \ref{11Sep06}.(2,3), statement 2 and (\ref{g6ab}), we have
$U^{(m)}= U^{m+2}$ for all $m\geq 0$.  $\Box $

\begin{corollary}\label{p15Oct06}
Let $K$ be a commutative ring. Then
\begin{enumerate}
\item the group $\Phi $ is generated by all the automorphisms
$\s_{i, \l x_ix_kx_l}$, i.e. $\Phi = \langle \s_{i, \l
x_ix_kx_l}\, | \, i=1, \ldots , n; \l \in K;  k<l; i\not\in \{
k,l\} \rangle$. The subgroups $\{ \s_{i, \l x_ix_kx_l}\}_{\l \in
K}$ of $\Phi$ form a minimal set of generators for $\Phi$. \item
$\Phi^{(i)}= \Phi^{2i+3} := \{ \s \in \Phi \, | \, (\s -1)(\gm )
\subseteq \gm^{2i+3} \}$, $i\geq 0$. \item Each element $\s \in
\Phi$ is a unique finite  product $\s = \cdots
\s_{b_7}\s_{b_5}\s_{b_3}$ for unique elements $b_i:= (b_{i1} ,
\ldots , b_{in})\in \L_{n,i}^n$ (see Corollary \ref{c11Sep06})
such that $b_{ij} \in (x_j)$ for all $j$.
\end{enumerate}
\end{corollary}

{\it Proof}. 3. Statement 3 follows from Lemma \ref{g4Oct06} and
Theorem \ref{11Sep06}.(2).

$1$. By statement 3, the elements of the type (\ref{1g3ab}) are
generators for the group $\Phi$. Then the result follows from
(\ref{g3ab}).

2. By Theorem \ref{11Sep06}.(2), $\Phi^{(i)}\subseteq \Phi^{2i+3}$
for all $i\geq 0$. The reverse inclusion follows from
 ({\ref{g3ab}). $\Box$


\section{The Jacobian group $\S$  and the equality
$\S = \S'\S''$ }\label{TJGS}

{\bf The Jacobian map}. Let $K$ be a commutative ring. Recall that
the group $\G$ consists of all automorphisms $\g_b$ : $x_1\mapsto
x_1+b_1, \ldots , x_n\mapsto x_n+b_n$, where $b:= (b_1, \ldots ,
b_n)$ is an $n$-tuple of  odd elements of $\gm^3$. Consider the
matrix $B:= \frac{\der b}{\der x}:= (\frac{\der b_i}{\der x_j})$
of the {\em skew  gradients} $\grad (b_i) : = (\frac{\der
b_i}{\der x_1}, \ldots , \frac{\der b_i}{\der x_n})$ for the
element  $b=(b_1, \ldots , b_n)$ (where $\frac{\der}{\der x_1},
\ldots , \frac{\der}{\der x_n}$ are the left partial skew
$K$-derivatives of $\L_n(K)$), and its {\em characteristic
polynomial}
$$ \det (t+B) = t^n+\sum_{i=1}^n \tr_i (B)t^{n-i}.$$
Clearly, $\tr_1(B) = \tr (B)= \sum_{i=1}^n \frac{\der b_i}{\der
x_i}$ is the  trace of the matrix $B$,   $\tr_n(B) = \det (B)$ is
its  determinant, and $\tr_i(B)= \sum_{1\leq j_1<\cdots <j_i\leq
n} \det (\frac{\der b_{j_\mu}}{\der x_{j_\nu}})_{\mu , \nu =1,
\ldots , i}$. Now, the jacobian of the automorphism $\g_b$ is
given by the rule 
\begin{equation}\label{JtriB}
\CJ (\g_b) = \det (t+B)|_{t=1} = 1+\sum_{i=1}^n \tr_i(B).
\end{equation}
Note that the sum of the  traces above is an element of $\gm^2$
since $\tr_i(B)\in \gm^{2i}$, $i\geq 1$. So, the {\em Jacobian
map} is given by the rule 
\begin{equation}\label{1JtriB}
\CJ  : \G \ra E_n', \g_b \mapsto  \CJ (\g_b)= \det (t+B)|_{t=1} =
1+\sum_{i=1}^n \tr_i(B).
\end{equation}
It is a  polynomial map in the coefficients of the elements $ b_1,
\ldots , b_n$. Recall that the abelian multiplicative group of
units $E_n'$ is equal to $E_n'= 1+ \sum_{i\geq 1} \L_{n, 2i}$.

{\bf The Jacobian group $\S$}. Let $K$ be a commutative ring.
Despite the fact that the jacobian map $\CJ : \G \ra E_n'$ is {\em
not} a group homomorphism, its  `kernel'
$$ \S := \{ \s \in \G \, | \, \CJ (\s ) =1\}$$
 is a {\em subgroup} of $\G
$ as it easily follows from (\ref{char1}) and (\ref{Js1i}). We
call $\S$ the {\bf Jacobian group}. This is a sophisticated
subgroup of $\G$. By (\ref{1JtriB}), the elements of the group
$\S$ are solutions to the system of  polynomial equations over $K$

\begin{equation}\label{Stri}
\S = \{ \g_b \in\G \, | \, \sum_{i=1}^n \tr_i(B) =0\}.
\end{equation}
So, $\S$ is a unipotent  algebraic group over the ring $K$. If $K$
is a field then $\S$ is an algebraic group over $K$ (in the usual
sense).  By Theorem \ref{s9Sep06}.(1), the system of polynomial
equations in (\ref{Stri}) can be made explicit 
\begin{equation}\label{1Stri}
\S = \{ \g_b \in\G \, | \, \phi (\der^\alpha (\sum_{i=1}^n
\tr_i(B))) =0, \;\; {\rm for \; all\; even} \;\;  \, 0\neq \alpha
\in \CB_n\}.
\end{equation}
The Jacobian group $\S$ is the solution to the system (\ref{Stri})
of skew differential operators. It looks like this is the first
example of a group of this kind.  The Jacobian group $\S$ is a
closed subgroup of $\G$ in the Zariski topology. The group $\G$
contains the descending chain of its normal subgroups
$$\G = \G^3\supset \G^5\supset \cdots \supset \G^{2m+1} : = \G
\cap U^{2m+1} \supset \cdots \supset
\G^{2[\frac{n}{2}]+1}\supseteq \G^{2[\frac{n}{2}]+3}=\{ e\},$$
with the abelian factors $\G^{2m+1} / \G^{2m+3}$ where $\G^{2m+1}
=\{ \s \in \G \, | \, (\s -1) (\gm )\subseteq \gm^{2m+1}\}$.

The group $\S$ contains the descending chain of its normal
subgroups
$$ \S = \S^3\supseteq \S^5\supseteq \cdots \supseteq \S^{2m+1} := \S
\cap \G^{2m+1} \supseteq \cdots \supseteq
\S^{2[\frac{n}{2}]+1}\supseteq \S^{2[\frac{n}{2}]+3}=\{ e\}$$ with
the abelian factors $\{ \S^{2m+1} / \S^{2m+3}\}$ since $\S^{2m+1}
/ \S^{2m+3}\subseteq \G^{2m+1} / \G^{2m+3}$, the abelian group.

{\bf The Jacobian group $\S$ and the image of the Jacobian map for
$n=3$}. For $n=3$,
$$ \G = \{\s_\l \, | \, \s (x_1)= x_1(1+\l_1x_2x_3), \s (x_2)= x_2(1+\l_2x_1x_3), \s (x_3)=
x_3(1+\l_3x_1x_2),\l \in K^3\}$$ where $\l := (\l_1, \l_2, \l_3)$,
and $\G\ra K^3$, $\s_\l \mapsto \l$, is the group isomorphism.
Since $\CJ (\s_\l ) = 1+\l_1x_2x_3+\l_2x_1x_3+\l_3x_1x_2$, the
Jacobian group $\S = \{ e\}$ is trivial, and ${\rm im }(\CJ
)=E_3'$,  i.e. the Jacobian map $\CJ : \G \ra E_3'$ is surjective.

\begin{lemma}\label{e8Oct06}
Let $K$ be a commutative ring and $\s, \tau \in \G$. Then
\begin{enumerate}
\item $\CJ (\s ) = \CJ (\tau )$ iff $\tau \in \s\S$.\item $\CJ
(\s^{-1}) = \CJ (\tau^{-1})$ iff $\tau \in \S \s$.
\end{enumerate}
\end{lemma}

{\it Proof}. $1$.  Note that $\CJ (\s ) = \s ( \CJ
(\s^{-1})^{-1})$ as it follows from the equality $ 1 = \CJ (\s
\s^{-1} ) = \CJ (\s ) \s (\CJ (\s^{-1}))$. Now, $\tau \in \s\S$
iff $\s^{-1} \tau \in \S$ iff $1= \CJ (\s^{-1} \tau ) = \CJ
(\s^{-1} ) \s^{-1} (\CJ (\tau ))$ iff $\CJ ( \tau ) = \s (\CJ
(\s^{-1})^{-1})= \CJ (\s )$.

$2$. By statement 1, $\CJ (\s^{-1}) = \CJ (\tau^{-1})$ iff
$\tau^{-1} \in \s^{-1}\S$ iff  $\tau \in \S \s$. $\Box $

{\it Remark}. Lemma \ref{e8Oct06} explains `intuitively' why, in
general,  the Jacobian group $\S$ is {\em not} a normal subgroup
of $\G$: note first that $\CJ (\s^{-1} ) = \s^{-1} (\CJ (\s
)^{-1})$ and $\CJ (\tau^{-1} ) = \tau^{-1} (\CJ (\tau )^{-1})$.
Suppose that $\S$ is a normal subgroup of $\G$, then $\s \S =
\S\s$ for all $\s \in \S$, and so the two statements of Lemma
\ref{e8Oct06} are equivalent, i.e. $\CJ (\s ) = \CJ (\tau )=:u$
iff $\CJ (\s^{-1} ) = \CJ (\tau^{-1} )$ iff $\s^{-1} (u)=
\tau^{-1} (u)$. Since the image of $\CJ$ is `big' there is no
reason to believe that the automorphisms  $\s^{-1}$ and
$\tau^{-1}$ acts always identically on $u$.

{\bf The image ${\rm im}(\CJ )$ and $\G / \S $}. By Lemma
\ref{e8Oct06}, the map
\begin{equation}\label{GSim}
 \CJ : \G / \S \ra {\rm im} (\CJ ), \;\; \s \S\mapsto \CJ (\s
),
\end{equation}
is a bijection.

{\bf The groups $\G_{2m}$, the Jacobian ascents}.  By
(\ref{cEns}), the abelian group $E_n'$ contains the descending
chain of subgroups
$$ E_n'=E_{n,2}'\supset E_{n,4}'\supset \cdots \supset E_{n,
2[\frac{n}{2}]}'\supset E_{n, 2[\frac{n}{2}]+2}'=\{ 1\}.$$
 If $K$ is a  commutative ring then, for each $m=1, 2,
 \ldots , [\frac{n}{2}]+1$, the preimage
\begin{equation}\label{G2m}
\G_{2m} := \G_{n,2m}:=\CJ^{-1} (E_{n,2m}') =\{ \s \in \G \, | \,
\CJ (\s ) \in E_{n,2m}'\}
\end{equation}
is, in fact, a  subgroup of $\G$: let $\s , \tau \in \G_{2m}$,
then
$$ \CJ (\s \tau ) = \CJ (\s ) \s (\CJ (\tau ))\subseteq
E_{n,2m}'\s (E_{n,2m}') \subseteq  E_{n,2m}'E_{n,2m}'\subseteq
E_{n,2m}',$$ i.e. $\G_{2m} \G_{2m}\subseteq \G_{2m}$; and
$$ \CJ (\s^{-1} ) = (\s^{-1} (\CJ (\s )))^{-1} \subseteq (\s^{-1}
(E_{n,2m}'))^{-1} \subseteq (E_{n,2m}')^{-1} = E_{n, 2m}',$$ i.e.
$\G^{-1}_{2m} \subseteq \G_{2m}$. $\Box$

Note that 
\begin{equation}\label{S=G2n}
\S = \G_{2[\frac{n}{2}]+2}.
\end{equation}
 We call the groups $\G_{2m}=\G_{n, 2m}$ the
{\bf Jacobian ascents}.  We have the descending chain of subgroups
in $\G$, the {\bf Jacobian filtration}: 
\begin{equation}\label{cG2m}
\G = \G_2\supseteq \G_4\supseteq \cdots \supseteq
\G_{2[\frac{n}{2}]}\supseteq \G_{2[\frac{n}{2}]+2}= \S.
\end{equation}
We will see later that all these groups are  distinct except the
last two if $n$ is even (Corollary \ref{a20Nov06}); each group
$\G_{2m+2}$ is a  normal subgroup of $\G_{2m}$ such that the
factor group $\G_{2m}/\G_{2m+2}$ is abelian (Lemma
\ref{n24Nov06}).

\begin{lemma}\label{n24Nov06}
Let $K$ be a commutative ring. For each natural number $m\geq 1$,
the group $\G_{n,2m+2}$ is a normal subgroup of $\G_{n,2m}$ such
that the factor group $\G_{n,2m}/\G_{n, 2m+2}$ is abelian.
\end{lemma}

{\it Proof}. Recall that the groups $\{ E_{n,2m}'\}$ are
$\G$-invariant. The result is an immediate consequence of the
following obvious fact: for each $m\geq 1$, $a\in E_{n,2m}'$, and
$ \s \in \G$, 
\begin{equation}\label{saa2m}
\s (a) \equiv a \mod E_{n, 2m+2}'.
\end{equation}
Indeed, to prove that $\G_{n,2m+2}$ is a normal subgroup of
$\G_{n,2m}$ we have to show that, for any $\s \in \G_{n,2m}$ and
$\tau \in \G_{n, 2m+2}$, $ \s \tau \s^{-1} \in \G_{n,2m+2}$, i.e.
$\CJ (\s \tau \s^{-1} )\in E_{n,2m+2}'$. By (\ref{char1}) and
(\ref{Js1i}),
$$ \CJ (\s \tau \s^{-1} ) = \CJ (\s )\,  \s (\CJ (\tau ))\, \s\tau (\CJ
(\s^{-1}))= \CJ (\s )\,  \s (\CJ (\tau ))\, \s\tau \s^{-1} (\CJ
(\s )^{-1}).$$ Note that $\s (\CJ ( \tau ))\in E_{n, 2m+2}'$ since
$\CJ (\tau ) \in  E_{n, 2m+2}'$; and $\s \tau \s^{-1} (\CJ (\s
)^{-1})\equiv \CJ (\s )^{-1} \mod E_{n,2m+2}'$, by (\ref{saa2m}).
Now,
$$ \CJ (\s \tau \s^{-1} ) \equiv \CJ (\s ) \CJ (\s )^{-1} \equiv 1
\mod  E_{n, 2m+2}',$$ i.e. $\CJ (\s \tau \s^{-1} ) \in  E_{n,
2m+2}'$, as required.

To prove that the factor group $\G_{n,2m}/\G_{n, 2m+2}$ is
abelian, we have to show that, for any $\s , \tau \in \G_{n,2m}$,
 $\s\tau \s^{-1}\tau^{-1}  \in
\G_{n,2m+2}$, that is $\CJ (\s\tau \s^{-1}\tau^{-1})\in
E_{n,2m+2}'$. By (\ref{char1}),  (\ref{Js1i}), and (\ref{saa2m}),
we have
\begin{eqnarray*}
 \CJ (\s\tau \s^{-1}\tau^{-1})&\equiv & \CJ (\s )\,  \s (\CJ (\tau ))\,  \s \tau (\CJ (\s^{-1} ))\, \s \tau \s^{-1} (\CJ (\tau^{-1})) \\
 &\equiv & \CJ (\s )\,  \s (\CJ (\tau ))\,  \s \tau \s^{-1} (\CJ (\s )^{-1} )\, \s \tau \s^{-1}\tau^{-1}  (\CJ (\tau )^{-1}) \\
 &\equiv & \CJ (\s )\CJ (\tau )\CJ (\s )^{-1}\CJ (\tau )^{-1} \equiv 1 \mod E_{n,2m+2}'. \;\;\; \Box\\
\end{eqnarray*}

The following result which is a part of Lemma \ref{n24Nov06} has a
more short and direct proof.

\begin{corollary}\label{n1Oct06}
The Jacobian group $\S$ is a normal subgroup of
$\G_{2[\frac{n}{2}]}$ such that the factor group
$\G_{2[\frac{n}{2}]}/\S $ is abelian.
\end{corollary}

{\it Proof}. Since each automorphism $\s \in \G$ acts trivially
(i.e. as the identity map) on $E_{n,2[\frac{n}{2}]}'$, the
Jacobian map $\CJ : \G_{2[\frac{n}{2}]}\ra E_{n,2[\frac{n}{2}]}'$,
$\tau \mapsto \CJ (\tau )$, is a group {\em homomorphism} ($\CJ
(\s \tau ) = \CJ (\s ) \s (\CJ (\tau )) = \CJ (\s ) \CJ (\tau )$)
with the kernel $\S$, hence $\S $ is a normal subgroup of
$\G_{2[\frac{n}{2}]}$ such that the factor group
$\G_{2[\frac{n}{2}]}/\S $ is abelian since the group $
E_{n,2[\frac{n}{2}]}'$ is abelian. $\Box $

The elements of the group $\G_{2m}$ are solutions to the system of
polynomial equations over $K$, 
\begin{equation}\label{Gtri}
\G_{2m} = \{ \g_b \in\G \, | \, \sum_{i=1}^n \tr_i(B)\in \gm^{2m}
\}.
\end{equation}
So, $\G_{2m}$ is an  algebraic unipotent group over the ring $K$.
By Theorem \ref{s9Sep06}.(1), the system of polynomial equations
in (\ref{Gtri}) can be made explicit 
\begin{equation}\label{1Gtri}
\G_{2m} = \{ \g_b \in\G \, | \, \phi (\der^\alpha (\sum_{i=1}^n
\tr_i(B)))=0 , \;\;\; {\rm for \; all \; even}\;  0\neq \alpha \in
\CB_n, \; 1\leq | \alpha | <2m \}.
\end{equation}

\begin{lemma}\label{m8Oct06}
Let $K$ be a  commutative ring. Then, $\G^{2m+1} \subseteq
\G_{2m}$, for each $m=1, 2, \ldots , [\frac{n}{2}]+1$.
\end{lemma}

{\it Proof}. Let $\s \in \G^{2m+1}$. Then $\s (x_1) = x_1+b_1,
\ldots , \s (x_n) = x_n+ b_n$, where all $b_i\in \L_{n , \geq
2m+1}^{od}$. Now, the result follows from 
\begin{equation}\label{Js=div}
\CJ (\s ) \equiv 1+ \sum_{i=1}^n \frac{\der b_i}{\der x_i} \equiv
1 \mod \gm^{2m},
\end{equation}
i.e. $\CJ (\s ) \in E_{2m}'$ since all $\frac{\der b_i}{\der
x_i}\in \gm^{2m}$.  $\Box $

Now, we introduce two important subgroups of $\S$, namely $\S'$
and $\S''$, and prove that $\G = \Phi \S''$ (Theorem
\ref{13Oct06}.(1)) and $\S = \S'\S''$ (Corollary
\ref{a8Oct06}.(1)).

{\bf The group $\S'$}. Consider the following subgroup of $\S$,
\begin{eqnarray*}
 \S'&:=& \S \cap \Phi = \{ \s \in \S \, | \, \s (x_1) \in (x_1),
\ldots , \s (x_n)\in (x_n)\}\\
 &=& \{ \s \in \S \, | \, \s (x_1) = x_1(1+a_1) ,\ldots , \s (x_n) =
 x_n(1+a_n)\}
\end{eqnarray*}
where each element $a_i \in K\lfloor x_1, \ldots , \widehat{x_i},
\ldots , x_n\rfloor^{ev}_{\geq 2}$. The group $\S'$ is a  closed
subgroup of $\S$ as the intersection of two closed subgroups $\S$
and $ \Phi$ of $\G$. It contains the descending chain of normal
subgroups
$$ \S' = \S'^3\supseteq \S'^5\supseteq \cdots \supseteq \S'^{2m+1} :=
\S' \cap \G^{2m+1} \supseteq \cdots \supseteq
\S'^{2[\frac{n}{2}]+1}\supseteq \S'^{2[\frac{n}{2}]+3}=\{ e\}$$
with the abelian factors $\{ \S'^{2m+1} / \S'^{2m+3}\}$.

{\it Example}. Let $a_1\in K\lfloor x_2, \ldots , x_n\rfloor_{2m}$
and $a_2\in K\lfloor x_1,x_3, \ldots , x_n\rfloor_{2m}$ be
homogeneous even elements of the same graded degree $2m\geq 2$,
and $\s \in \G$: $x_1\mapsto x_1(1+a_1), x_2\mapsto x_2(1+a_2)$,
$x_j\mapsto x_j$, $j\geq 3$. Then $\s \in \S$ iff $$1= \CJ (\s ) =
\det
\begin{pmatrix}
1+a_1 & -x_1\frac{\der a_1}{\der x_2}\\
-x_2\frac{\der a_2}{\der x_1} & 1+a_2
\end{pmatrix} = 1+a_1+a_2+a_1a_2+x_1x_2\frac{\der a_1}{\der x_2}\frac{\der a_2}{\der
x_1}$$ iff $a_1=-a_2\in K\lfloor x_3, \ldots , x_n\rfloor_{2m}$
and $a_1^2=0$. So, {\em for each even homogeneous element $a\in
K\lfloor x_3, \ldots , x_n\rfloor_{2m}$ such that $a^2=0$ the
automorphism} 
\begin{equation}\label{sama}
\s \in \G: x_1\mapsto x_1(1+a),\; x_2\mapsto x_2(1-a),\;
x_j\mapsto x_j, \; j\geq 3,
\end{equation}
 {\em belongs to the group} $\S'$.

{\bf The group $\S''$}. For each $i=1, \ldots , n$, and $b_i\in
K\lfloor x_1, \ldots , \widehat{x_i}, \ldots ,
x_n\rfloor^{od}_{\geq 3}$, consider the element of $\S$:
\begin{equation}\label{ksibi}
\xi_{i,b_i}: x_i\mapsto x_i+ b_i, \; x_j\mapsto x_j, \; \forall \,
j\neq i.
\end{equation}
Let $\S''$ be the subgroup of $\S$ generated by all the elements
$\xi_{i, b_i}$, $1\leq i \leq n$. For each $i=1, \ldots , n$, let
$$ \S''_i:= \{ \xi_{i, b_i}\, | \, b_i\in K\lfloor x_1, \ldots ,\widehat{x_i}, \ldots ,
x_n\rfloor^{od}_{\geq 3}\}.$$ Since $\xi_{i, b_i}\xi_{i, b_i'}
=\xi_{i, b_i+b_i'}$ and $ \xi_{i, b_i}^{-1} =  \xi_{i,- b_i}$, the
set $\S''_i$ is an {\em abelian} group canonically isomorphic to
the abelian additive group $K\lfloor x_1, \ldots ,\widehat{x_i},
\ldots , x_n\rfloor^{od}_{\geq 3}$ via $ \xi_{i,b_i} \mapsto b_i$.
Therefore, the group $\S''_i$ is the direct product of its
one-dimensional abelian subgroups isomorphic to $(K,+)$,
\begin{equation}\label{Siss1}
\S''_i=\prod_{\alpha}\{ \xi_{i, \l x^\alpha }\}_{\l \in K}
\end{equation}
where $\alpha$ runs through all the odd subsets of the set $\{ 1,
\ldots , \widehat{i}, \ldots , n\}$ with $ |\alpha | \geq 3$; the
map $(K,+)\ra \{ \xi_{i, \l x^\alpha } \}_{\l \in K}$, $\l \mapsto
\xi_{i, \l x^\alpha}$, is a group  isomorphism.

So, {\em the group $\S''$ is generated by its abelian subgroups} $
\S''_1, \ldots , \S''_n$.

The commutator of the elements $\xi_{i, b_i}\in \S_i''$ and
$\xi_{j, b_j}\in \S_j''$ where $i\neq j$ is given by the rule
\begin{eqnarray}\label{Siss2}
[\xi_{i, b_i}, \xi_{j, b_j}] : & & x_i\mapsto
(1-bb')x_i-b^2b'x_j-b(c'+b'c),\\
& & x_j\mapsto bb'^2x_i+(1+bb'+(bb')^2) x_j+b'(c+bc' +bb'c),\\
&&  x_k\mapsto x_k, \; k\neq i,j,
\end{eqnarray}
where $b_i=bx_j+c$ and $b_j=b'x_i+c'$ for unique elements $b,b'\in
K\lfloor x_1, \ldots ,\widehat{x_i}, \ldots , \widehat{x_j},
\ldots , x_n\rfloor^{ev}_{\geq 2}$ and $c,c'\in K\lfloor x_1,
\ldots ,\widehat{x_i}, \ldots , \widehat{x_j}, \ldots ,
x_n\rfloor^{od}_{\geq 3}$.

Given automorphisms $\xi_{1, b_1}, \ldots \xi_{n, b_n}$, and $\xi
:= \prod_{i=1}^n\xi_{i, b_i}$ is their product in an arbitrary
(fixed) order, then 
\begin{equation}\label{1ksibi}
\xi (x_i)=x_i+ b_i +\cdots ,\;\; i=1, \ldots , n.
\end{equation}
The group $\S''$ contains the descending chain of normal subgroups
$$ \S'' = \S''^3\supseteq \S''^5\supseteq \cdots \supseteq \S''^{2m+1} :=
\S'' \cap \G^{2m+1} \supseteq \cdots \supseteq
\S''^{2[\frac{n}{2}]+1}\supseteq \S''^{2[\frac{n}{2}]+3}=\{ e\}$$
with the abelian factors $\{ \S''^{2m+1} / \S''^{2m+3}\}$.

A direct calculation gives 
\begin{equation}\label{xijam}
[\xi_{i, \l x_jx^\alpha }, \xi_{j, \mu x^\beta}]= \xi_{i, -\l \mu
x^\alpha x^\beta}
\end{equation}
where $i\neq j$; $\alpha$ is an even set and $\beta$ is an  odd
set such that the sets $\{ i, j \} $ and $\alpha\cup \beta$ are
disjoint. Similarly, 
\begin{equation}\label{xijam1}
[\xi_{i, \l x^\alpha }, \xi_{j, \mu x^\beta}]= 0
\end{equation}
for all  $i\neq j$; $\l , \mu \in K$; $\alpha$ and $\beta$ are odd
sets such that $|\alpha | \geq 3$ and $|\beta |\geq 3$,
 and the sets $\{ i, j \} $ and $\alpha\cup \beta$ are disjoint.

One can verify that for each $i\neq j$, and even sets $\alpha$ and
$\beta$ such that $\{ i, j \} $ and $\alpha\cup \beta$ are
disjoint, $\l , \mu \in K$, the commutator (which is an element of
$\S''$) 
\begin{equation}\label{com1}
[\xi_{i, \l x_jx^\alpha }, \xi_{j, \mu x_ix^\beta}]: x_i\mapsto
x_i(1-\l \mu x^\alpha x^\beta ), \; x_j\mapsto x_j(1+\l \mu
x^\alpha x^\beta ),\; x_k\mapsto x_k, \;  k\neq i,j,
\end{equation}
belongs to the group $\S'$, and
$$
[\xi_{i, \l x_jx^\alpha }, \xi_{j, \mu x_ix^\beta}]^{-1}:
x_i\mapsto x_i(1+\l \mu x^\alpha x^\beta ), \; x_j\mapsto x_j(1-\l
\mu x^\alpha x^\beta ),\; x_k\mapsto x_k, \;  k\neq i,j.$$ Now,
the next corollary is obvious since $[\xi_{i, \l x_jx^\alpha},
\xi_{j, \mu x_ix^\beta} ] \in \S'\cap \S''\subseteq \Phi \cap
\S''$.

\begin{corollary}\label{x16Oct06}
If $K$ is a commutative ring and $n\geq 6$ then $\S'\cap \S''\neq
\{ e \}$ and $\Phi \cap \S''\neq \{ e\}$.
\end{corollary}

A straightforward calculation gives  
\begin{equation}\label{dvac1}
[\xi_{i, \nu x_jx^\g} , [\xi_{i, \l x_jx^\alpha }, \xi_{j, \mu
x_ix^\beta}]]= \xi_{i, -2\l \mu \nu x_jx^\alpha x^\beta x^\g},
\end{equation}
for all $\l , \mu , \nu \in K$; the sets $\alpha$, $\beta$, and
$\g$ are even and non-empty; the sets $\{ i,j\}$ and $ \alpha \cup
\beta \cup \g$ are disjoint. Similarly, 
\begin{equation}\label{dvac2}
[\xi_{i, \nu x^\g} , [\xi_{i, \l x_jx^\alpha }, \xi_{j, \mu
x_ix^\beta}]]= \xi_{i, -\l \mu \nu x^\alpha x^\beta x^\g},
\end{equation}
for all $\l , \mu , \nu \in K$; the sets $\alpha$ and  $\beta$,
are  even and non-empty; the set $\g$ is  odd and $|\g |\geq 3$;
the sets $\{ i,j\}$ and $ \alpha \cup \beta \cup \g$ are disjoint.

{\bf The group $\S''$ for $n=4,5,6$}.  Let $A$ be a group and
$A_1, \ldots , A_s$ be its subgroups. We say that the group $A$ is
an {\bf exact product} of the groups $A_i$, $$A:= A_1\cdots
A_s[exact]:= {}^{\rm ex}\prod_{i=1}^nA_i$$ if each element $a\in
A$ is a unique product $a=a_1\cdots a_s$ of elements $a_i\in A_i$.

The next lemma describes the structure of the group $\S''$ for
small values of $n=4,5,6$. These values are rather peculiar as
Theorem \ref{s15Oct06} shows.
\begin{lemma}\label{a17Oct06}
Let $K$ be a commutative ring.
\begin{enumerate}
\item If $n=4$ then $\S''= \S''_1\times \cdots \times \S''_4$ is
the abelian group.  \item If $n=5$ then $\S''= \S''_1\times \cdots
\times \S''_5$ is the abelian group. \item If $n=6$ then
\begin{enumerate}
\item $\S''= Z(\S'')\times  \prod_{i,j,k,l} \S''_{i;j,k,l}$ is the
exact product
 of the centre $Z(\S'')$ of $\S''$ and the one dimensional abelian
 subgroups $\S''_{i;j,k,l}:= \{ \xi_{i, \l x_jx_kx_l}\}_{\l \in K} \simeq K$
 where $i=1, \ldots , 6$; $j<k<l$; $i\not\in \{ j,k,l,\}$.
\item $[\S'', \S'']=\S'\cap \S''$ and the group $\S'\cap \S''$ is
the direct product $\prod_{i<j}C_{ij}$ of its subgroups $C_{ij}:=
\{ c_{ij,\l }: x_i\mapsto x_i(1-\l x^\alpha ), x_j\mapsto x_j(1+\l
x^\alpha ), x_k\mapsto x_k, k\neq i,j\}$ where $\alpha := \{ 1,
\ldots , 6\}\backslash \{ i,j\}$ and $C_{ij}\simeq (K,+)$ via
$c_{ij, \l} \mapsto \l $. \item $Z(\S'')= \S''^5=(\S'\cap \S'')
\times \prod_{i=1}^6\S''_{i,5}$ where $\S''_{i,5}:= \{
\xi_{i,b_i}\, | \, b_i\in K\lfloor x_1, \ldots , \widehat{x_i},
\ldots , x_6\rfloor_5\}$.
\end{enumerate}
\end{enumerate}
\end{lemma}

{\it Proof}. 1 and 2. If $n=4,5$ then the elements of $\S_i''$
commute with the elements of $\S_j''$, hence statements 1 and 2
are obvious.

3. Let $n=6$.  The group $\S''$ is generated by its abelian
subgroups $\S_1'', \ldots , \S_6''$, and so, by (\ref{Siss1}), the
group $\S''$ is generated by the 1-dimensional abelian subgroups
$\langle \xi_{i, \l x^\alpha }\rangle_{\l \in K}\simeq(K,+)$, $
\xi_{i, \l x^\alpha }\mapsto \l$, where $|\alpha | = 3,5$. If
$|\alpha | = 5$ then all the $\xi_{i, \l x^\alpha}\in Z:=
Z(\S'')$, the centre of the group $\S''$. Since $n=6$,
$\S''^5\subseteq Z$ and the RHS of (\ref{xijam}) is equal to zero.
By Theorem \ref{11Sep06}.(2), $[\S'', \S'']\subseteq \S''^5$, and
so $[\S'', \S'']\subseteq Z$. By (\ref{xijam}), (\ref{xijam1}),
and (\ref{com1}), the only nontrivial commutators come from
(\ref{com1}) and only in the case when $i<j$ and $\{ i,j\}\cup
\alpha \cup \beta = \{ 1, \ldots , 6\}$ there. In this case, the
commutators (\ref{com1}) is the automorphism $c_{ij, \l \mu }\in
Z$. It is obvious that the product of the groups $C_{ij}$ is the
direct product $\prod_{i<j}C_{ij}$, and $[\S'' , \S'']=
\prod_{i<j}C_{ij}\subseteq \S'\cap \S''^5$. It follows that
$\S''^5= [\S'', \S'']\times \prod_{i=1}^6\S''_{i,5}$  is the
direct product of groups since the abelian group $\S''^5$ is
generated by the subgroups $[\S'', \S'']$ and
$\prod_{i=1}^6\S''_{i,5}$, and their intersection is trivial.
Since $n=6$, $\S'' = \S''^5\times \prod_{i,j,k,l} \S''_{i; j,k,l}$
is the exact product of groups where $i,j,k,l$ are as in (a), then
$\S'\cap \S'' = \S' \cap \S''^5$. Since $[\S'' , \S'' ] \subseteq
\S'$ and $\S''^5=[\S'', \S'']\times \prod_{i=1}^6\S''_{i,5}$, we
have
$$ \S'\cap \S'' = [\S'', \S'']\times (\S'\cap \prod_{i=1}^6\S''_{i,5}) = [\S'', \S'']= \prod_{i<j} C_{ij} $$
since $\S'\cap \prod_{i=1}^6\S''_{i,5}=\{ e\}$. This proves
statement (b).

It follows from the decomposition $\S'' = \S''^5\times
\prod_{i,j,k,l} \S''_{i;j,k,l}$ and (\ref{com1}) that $Z\subseteq
\S''^5$. Since $\S''^5\subseteq Z$, we have $Z= \S''^5$. Since
 $\S'' = \S''^5\times
\prod_{i,j,k,l} \S''_{i;j,k,l}$ , (a) follows. Since
$\S''^5=[\S'', \S'']\times \prod_{i=1}^6\S''_{i,5}$, (c) follows.
  $\Box $

{\bf A minimal set of generators for the group $\S''$}. The next
result provides a (minimal) set of generators for the group
$\S''$.

\begin{theorem}\label{s15Oct06}
Let $K$ be a commutative ring and $n\geq 4$. Then
\begin{enumerate}
\item if either $n$ is odd; or $n=4$; or $n$ is even and $n\geq 8$
and $\frac{1}{2}\in K$, then
$$\S''= \langle \s_{i, \l x_jx_kx_l} \, | \, \l \in K; i=1, \ldots
, n; j<k<l; i\not\in \{ j,k,l\}\rangle,$$ and the  subgroups
$\{\s_{i, \l x_jx_kx_l}\}_{\l \in K}\simeq K$ of $\S''$ form a
minimal set of generators for $\S''$. \item If  $n$ is even, then
$$\S''= \langle \s_{i, \l x_jx_kx_l}, \; \s_{i, \l x_1\cdots \widehat{x_i}\cdots x_n} \, | \, \l \in K; i=1, \ldots
, n; j<k<l; i\not\in \{ j,k,l\}\rangle . $$
\end{enumerate}
\end{theorem}

{\it Proof}. 1.  Recall that the group $\S''$ is generated by its
abelian subgroups $\S_1'', \ldots \S_n''$. In order to prove the
claims that the elements above generate the group $\S''$ it
suffices to show that each automorphism $\xi_{i, \l x^\alpha}$
(where $\alpha$ is odd, $|\alpha |\geq 3$, and $i\not\in \alpha$)
is a product of some of them. By (\ref{xijam}), for $\l \in K$ and
 distinct indices $i,j,k_1, l_1, k_2, l_2, \ldots , k_{2m}, l_{2m},
p,q,r$ there is the equality 
\begin{eqnarray}\label{xipq1}
 \xi_{i, \l (x_{k_1}x_{l_1}) (x_{k_2}x_{l_2})\cdots
(x_{k_{2m}}x_{l_{2m}})x_px_qx_r}&=& [\xi_{i, x_jx_{k_1}x_{l_1}}, [
\xi_{j, x_ix_{k_2}x_{l_2}},[\xi_{i, x_jx_{k_3}x_{l_3}},[\xi_{j,
x_ix_{k_4}x_{l_4}},\ldots\\
 && [ \xi_{j, x_ix_{k_{2m}}x_{l_{2m}}},\xi_{i, \l x_px_qx_r}]\ldots ].
\end{eqnarray}
Similarly,  for $\l \in K$ and distinct indices $i,j,k_1, l_1,
k_2, l_2, \ldots , k_{2m+1}, l_{2m+1}, p,q,r$ there is the
equality 
\begin{eqnarray}\label{xipq2}
 \xi_{i, \l (x_{k_1}x_{l_1}) (x_{k_2}x_{l_2})\cdots
(x_{k_{2m+1}}x_{l_{2m+1}})x_px_qx_r}&=& [\xi_{i,
x_jx_{k_1}x_{l_1}}, [ \xi_{j, x_ix_{k_2}x_{l_2}},[\xi_{i,
x_jx_{k_3}x_{l_3}},[\xi_{j,
x_ix_{k_4}x_{l_4}},\ldots\\
 && [ \xi_{i, x_jx_{k_{2m+1}}x_{l_{2m+1}}},\xi_{j, -\l x_px_qx_r}]\ldots ].
\end{eqnarray}
Suppose that  $n$ is odd. The set  $\{ i\} \cup \alpha$ (for
$\xi_{i, \l x^\alpha}$)  is an even set hence it is not equal to
the set $\{ 1, \ldots , n\}$, and so one can find element $j$ such
that $j\not\in \{ i \}\cup \alpha$. It follows from (\ref{xipq1})
and (\ref{xipq2}) that the 1-dimensional subgroups $\{ \s_{i, \l
x_jx_kx_l}\}$ are generators for the group $\S''$. They are a
minimal set of generators due to the isomorphism in Theorem
\ref{11Sep06}.(2) in the case $m=1$ there.

Let $n\geq 4$ be an even number. If $n=4$ the result is obvious.
If $n=6$ the result is {\em not true} by Lemma Lemma
\ref{a17Oct06}.(3). So, let $n\geq 8$ and $\frac{1}{2}\in K$. If
$\{ i \} \cup \alpha \neq \{ 1, \ldots , n\}$ then we have already
proven that the automorphism $\xi_{i, \l x^\alpha}$ is a product
of the elements $\s_{i', \l'x_jx_kx_l}$. If $\{ i \} \cup \alpha =
\{ 1, \ldots , n\}$ then the automorphism $\xi_{i, \l x^\alpha} =
\xi_{i, \l x_1\cdots \widehat{x_i}\cdots x_n}$ is a product of the
elements $\s_{i', \l'x_jx_kx_l}$, by (\ref{dvac1}).

2. This statement has been proven already in the proof of
statement 1 (see the last two sentences above). $\Box $

The Jacobian group is a complicated group. To understand its
structure first cut it with the small group $\CE_{n,i}$.

\begin{lemma}\label{s8Oct06}
Let $K$ be a commutative ring. Then $$\CE_{n,i} \cap \S =
\CE_{n,i}'':= \{ \xi_{i,b_i}\, | \, b_i\in  K\lfloor x_1, \ldots ,
\widehat{x_i}, \ldots  x_n\rfloor_{\geq 3}^{od}\}.$$
\end{lemma}

{\it Proof}. A typical element of the group $\CE_{n,i}$ is as
follows $\g_{1+a, b} : x_i\mapsto x_i(1+a)+b$, $ x_j\mapsto x_j$,
$j\neq i$ (see (\ref{CEni})).  Then $\CJ ( \g_{1+a, b}) = 1+a$.
So, $\g_{1+a, b} \in \CE_{n,i} \cap \S$ iff $a=0$, i.e. $\CE_{n,i}
\cap \S = \CE_{n,i}''$.  $\Box $

{\bf The equality $\G =\Phi \S''$}. For a natural number $n$,
 let $\od (n)$ be the largest odd
number such that $\leq n$, and let  $\Od (n)$ be the set of all
odd natural numbers $j$ such that $3\leq j \leq n$, i.e. $\Od (n)
= \{ 3,5, \ldots , \od (n)\}$.

 For each $ t=3,5, \ldots , \od (n)$, let
$$\CF_{n,t}'':= \{ \s \in \S''\, | \, \s (x_i) =
x_i+b_i, \; b_i\in K\lfloor x_1, \ldots , \widehat{x_i}, \ldots ,
x_n\rfloor_t,\;  i=1,
 \ldots , n\}.$$ The set $\CF_{n,t}''$
is an affine  variety over $K$ of dimension $n{n-1\choose t}$,
i.e. the algebra of regular functions on $\CF_{n,t}''$ is a
polynomial algebra in $n{n-1\choose t}$ variables where the
coordinate functions are the coefficients of the polynomials
$b_i\in K\lfloor x_1, \ldots , \widehat{x_i}, \ldots ,
x_n\rfloor_{ t}$. Let 
\begin{equation}\label{dCEn}
\CF_n'':= \CF_{n, \od (n)}''\times \cdots \times \CF_{n,5}''\times
\CF_{n,3}''
\end{equation}
be the product of algebraic varieties. Then 
\begin{equation}\label{dimEns}
\dim (\CF_n'')= n\sum_{s=1}^{\od (n)}{n-1\choose 2s+1}= n
(\sum_{s=0}^{\od (n)}{n-1\choose 2s+1}-n+1)= n(2^{n-2}-n+1).
\end{equation}

 In general, $\Phi\cap
\S''\neq \{ e\}$ (Corollary \ref{x16Oct06}). The next theorem
shows that $\G = \Phi \S''$ and that any element of $\G$ is a
unique product of elements from $\Phi$ and $\S''$ when one puts
certain conditions on the choice of the multiples.

\begin{theorem}\label{13Oct06}
Let $K$ be a commutative ring. Then
\begin{enumerate}
\item $\G = \Phi \S''$. \item $\G^j = \Phi^j \S''^j$, $j=3,5,
\ldots , \od (n)$.\item Let   $j=3,5, \ldots , \od (n)$. Each
  $\s \in \G^j$ is a unique product $ \s = \phi_j
\xi_{\od (n)} \cdots \xi_{j+2} \xi_j$ where $\phi_j\in \Phi^j$ and
each  $ \xi_k$ is as in (\ref{1ksibi}) with
 $\xi_k(x_i) - x_i\in K\lfloor x_1, \ldots , \widehat{x_i}, \ldots ,
x_n\rfloor_k$, $i=1, \ldots , n$. Moreover, for all $i=1, \ldots ,
n$  the following conditions hold:
\begin{enumerate}
\item $\s^{-1} (x_i) \equiv - \xi_j (x_i) \mod ((x_i) + \L_{n,
\geq j+2}^{od})$, \item $\xi_k \cdots \xi_{j+2} \xi_j \s^{-1}
(x_i) \equiv - \xi_{k+2} (x_i) \mod ((x_i) + \L_{n, \geq
k+4}^{od})$, $k=j,j+2, \ldots , \od (n)-2$.
\end{enumerate}
\item For each odd natural number $j$ such that $3\leq j \leq n$,
let $\CF_n''^j:= \CF''_{\od (n)}\times \cdots \times
\CF_{j+2}''\times \CF_j''$. The map
$$ \Phi^j\times \CF_n''^j\ra \G^j, \;\; (\phi , \xi_{\od (n)},
\ldots , \xi_{j+2}, \xi_j) \mapsto \phi  \xi_{\od (n)} \cdots
\xi_{j+2} \xi_j, $$ is an isomorphism of the algebraic varieties
with the inverse map $\s \mapsto \phi  \xi_{\od (n)} \cdots
\xi_{j+2} \xi_j$ given by the decomposition of statement 3. \item
The map
$$ \Phi\times \CF_n''\ra \G, \;\; (\phi , \xi_{\od (n)},
\ldots , \xi_5, \xi_3) \mapsto \phi  \xi_{\od (n)} \cdots \xi_5
\xi_3, $$ is an isomorphism of the algebraic varieties with the
inverse map $\s \mapsto \phi  \xi_{\od (n)} \cdots \xi_5 \xi_3$
given by the decomposition of statement 3.
\end{enumerate}
\end{theorem}

{\it Proof}. 1. Statement 1 is a particular case of  statement 2
when $j=3$ since $ \G = \G^3$, $\Phi = \Phi^3$, and $\S'' =
\S''^3$.

2. Statement 2 follows from statement 3.

3. First, we prove that there exists a decomposition for $\s$ that
satisfies the conditions (a) and (b), then we prove the
uniqueness.

By the inversion formula (Theorem \ref{i9Sep06}), the map $\G \ra
\G $, $\s \mapsto \s^{-1}$, is an automorphism of algebraic
varieties. Let $\s \in \G^j$. Then $\s^{-1} \in \G^j$, and, for
each $i=1, \ldots , n$, $\s^{-1} (x_i) = x_i(1+ a_i) + b_i$ for
some elements $a_i \in K\lfloor x_1, \ldots , \widehat{x_i},
\ldots x_n\rfloor_{\geq j-1}^{ev}$ and $ b_i \in K\lfloor x_1,
\ldots , \widehat{x_i}, \ldots x_n\rfloor_{\geq j}^{od}$. Then
$b_i= c_i +\cdots $ for some element  $c_i \in K\lfloor x_1,
\ldots , \widehat{x_i}, \ldots x_n\rfloor_j$. Clearly,
\begin{equation}\label{1xi2m}
\s^{-1} (x_i) \equiv c_i \mod ((x_i) + \L_{n, \geq j+2}^{od}).
\end{equation}
Define the element $\xi_j \in \S''^j$ by the rule
\begin{equation}\label{xi2m1}
\xi_j (x_i) = x_i-c_i, \; \; {\rm for\; all}\;\; i.
\end{equation}
 Then the automorphism $\xi_j$
satisfies the condition $(a)$. Consider the automorphism $$\s':=
\xi_j \s^{-1} \in \G^j: x_i\mapsto x_i( 1+a_i') + b_i'$$ where
$a_i' \in K\lfloor x_1, \ldots , \widehat{x_i}, \ldots
x_n\rfloor_{\geq j-1}^{ev}$ and  $ b_i' \in K\lfloor x_1, \ldots ,
\widehat{x_i}, \ldots x_n\rfloor_{\geq j+2}^{od}$.  Then $b_i'=
c_i' +\cdots $ for some $c_i' \in K\lfloor x_1, \ldots ,
\widehat{x_i}, \ldots x_n\rfloor_{j+2}$. Clearly,
\begin{equation}\label{2xi2m}
\xi_j\s^{-1} (x_i) \equiv c_i' \mod ((x_i) + \L_{n, \geq
 j+4}^{od}).
\end{equation}
 Define the element $\xi_{j+2} \in \S''^{j+2}$
by the rule $\xi_{j+2} (x_i) = x_i-c_i'$ for all $i$. Then
$\xi_{j+2}$ satisfies the condition $(b)$  for $k=j$. Now, we can
repeat the same argument for the automorphism $\s'' := \xi_{j+2}
\xi_j \s^{-1}$. Continue in this fashion we finally come to the
inclusion
$$\xi_{\od (n)} \cdots \xi_{j+2} \xi_j\s^{-1} \in
\Phi^j,$$ and obtain a  decomposition for $\s$  that  satisfies
the properties (a) and (b) exists.

Uniqueness: Let $ \s = \phi'\xi'_{\od (n)}\cdots \xi_{j+2}'\xi_j'$
be another decomposition with $\xi'_j (x_i) = x_i-\l_i$ for some
$\l_i\in K\lfloor x_1, \ldots , \widehat{x_i}, \ldots ,
x_n\rfloor_j$, $1\leq i \leq n$. Then $\xi_j'^{-1} (x_i) =
x_i+\l_i$.  Since $\s^{-1} (x_i ) \equiv \xi_j'^{-1} (x_i) \equiv
\l_i \mod ((x_i) +\L^{od}_{n, \geq j+2})$, we must have $ \xi'_j
=\xi_j$, by (\ref{1xi2m}). Similarly, (\ref{2xi2m}) yields the
equality $\xi'_{j+2} = \xi_{j+2}$. Using the same argument again
and again (or by induction) we see that $\xi_k'=\xi_k$ for all
$k=j, j+2, \ldots , \od (n)$. These equalities imply that $\phi' =
\phi$.

4. This statement follows from statement 3 since the map $\s
\mapsto \phi \xi_{\od (n)}\cdots \xi_{j+2}\xi_j$ is a polynomial
map as the proof of statement 3 shows.

5. This statement is a particular case of statement 4 for $j=3$
since $\Phi = \Phi^3$, $\CF''= \CF''^3$, and $\G = \G^3$.
 $\Box $

In the proof of Theorem \ref{13Oct06}, the algorithm is given for
finding the automorphisms $\phi$ and $\xi_i$ in the presentation
$\s = \phi \xi_{\od (n)}\cdots \xi_{j+2}\xi_j$.

\begin{corollary}\label{c13Oct06}
Let $K$ be
a commutative ring. Then $\CJ (\G ) = \CJ (\Phi)$.
\end{corollary}

{\it Proof}. Let $\s \in \G$. Then $\s = \tau \s''$ for some $\tau
\in \Phi$ and $\s''\in \S''$ (Theorem \ref{13Oct06}.(1)). Then
$$\CJ (\s ) = \CJ (\tau \s'')= \tau (\CJ (\s''))\CJ (\tau ) = \tau
(1) \CJ (\tau ) = \CJ (\tau ).$$ Therefore, $\CJ (\G ) = \CJ
(\Phi)$. $\Box $

{\bf The equality $\S = \S'\S''$}.  The next result shows that the
Jacobian group $\S$ is the product of its subgroups $\S'$ and
$\S''$, i.e.  each element $\s \in \S$ is a product $\s = \s'\s''$
of some elements $\s'\in \S'$ an $\s''\in \S''$. This product is
not unique as,  in general, $\S'\cap \S''\neq \{ e\}$ (Corollary
\ref{x16Oct06}). Though, by putting extra conditions on the choice
of the elements $\s'$ and $\s''$ the uniqueness can be preserved.

\begin{corollary}\label{a8Oct06}
Let $K$ be a commutative ring. Then
\begin{enumerate}
\item $\S = \S' \S''$. \item $\S^j = \S'^j \S''^j$, $j=3,5, \ldots
, \od (n)$. \item Let   $j=3,5, \ldots , \od (n)$. Each
  $\s \in \S^j$ is a unique product $ \s = \s'
\xi_{\od (n)} \cdots \xi_{j+2} \xi_j$ where $\s'\in \S^j$ and each
$ \xi_k$ is as in (\ref{1ksibi}) with
 $\xi_k(x_i) - x_i\in K\lfloor x_1, \ldots , \widehat{x_i}, \ldots ,
x_n\rfloor_k$, $i=1, \ldots , n$. Moreover, for all $i=1, \ldots ,
n$  the following conditions hold:
\begin{enumerate}
\item $\s^{-1} (x_i) \equiv - \xi_j (x_i) \mod ((x_i) + \L_{n,
\geq j+2}^{od})$, \item $\xi_k \cdots \xi_{j+2} \xi_j \s^{-1}
(x_i) \equiv - \xi_{k+2} (x_i) \mod ((x_i) + \L_{n, \geq
k+4}^{od})$, $k=j,j+2, \ldots , \od (n)-2$.
\end{enumerate}
\item For each odd natural number $j$ such that $3\leq j \leq n$,
the map
$$ \S'^j\times \CF_n''^j\ra \S^j, \;\; (\s' , \xi_{\od (n)},
\ldots , \xi_{j+2}, \xi_j) \mapsto \s'  \xi_{\od (n)} \cdots
\xi_{j+2} \xi_j, $$ is an isomorphism of the algebraic varieties
with the inverse map $\s \mapsto \s'  \xi_{\od (n)} \cdots
\xi_{j+2} \xi_j$ given by the decomposition of statement 3. \item
The map
$$ \S'\times \CF_n''\ra \G, \;\; (\s' , \xi_{\od (n)},
\ldots , \xi_5, \xi_3) \mapsto \s'  \xi_{\od (n)} \cdots \xi_5
\xi_3, $$ is an isomorphism of the algebraic varieties with the
inverse map $\s \mapsto \s'  \xi_{\od (n)} \cdots \xi_5 \xi_3$
given by the decomposition of statement 3.
\end{enumerate}
\end{corollary}

{\it Proof}. 1. By Theorem \ref{13Oct06}.(1), $\G = \Phi \S''$.
Note that $\S \subseteq \G$, $\S''\subseteq \S$, and $ \S'= \S
\cap \Phi$.  Now,
$$ \S = \S \cap \G = \S \cap \Phi \S''= (\S \cap \Phi ) \S'' = \S'
\S''.$$ This proves statement 1. The rest follows at once from
Theorem \ref{13Oct06}. $\Box$

By Corollary \ref{a8Oct06} and (\ref{dimEns}), in order to find
the dimension (and generators) for the Jacobian group $\S$ it
suffices to find the dimension (and generators) for $\S'$.

{\bf The largest normal subgroup of $\G$ in $\S$}.

Let $A\subseteq B$ be groups. Then 
\begin{equation}\label{NAB}
\CN (A, B):= \{ a\in A\, | \, bab^{-1} \in A, \;\; {\rm for \;
all} \;\; b\in B\},
\end{equation}
is a  normal subgroup of $B$ in $A$. If $N$ is a normal subgroup
of $B$ that is contained in $A$ then $N\subseteq \CN (A,B)$.
Therefore, $\CN (A,B)$ is the {\em largest normal} subgroup of $B$
that is contained in $A$. {\em The group $ A$ is a normal subgroup
of $B$ iff $A= \CN (A, B)$}.

\begin{theorem}\label{c8Oct06}
Let $K$ be  a commutative ring. Then
$$ \CN (\S , \G ) = \{ \tau \in \S \, | \, \tau (\l ) = \l \;\;
{\rm for \; all} \;\; \l \in {\rm im} (\CJ ) \}.$$
\end{theorem}

{\it Proof}. $\tau \in \CN (\S , \G )$ iff $\s \tau \s^{-1} \in
\S$ for all $\s \in \G$ iff
$$ 1= \CJ ( \s \tau \s^{-1} ) = \CJ (\s ) \s (\CJ ( \tau )) \s
\tau (\CJ (\s^{-1}))= \CJ (\s ) \s \tau  (\CJ (\s^{-1}))$$ iff
$\tau (\CJ ( \s^{-1} )) = \s^{-1} (\CJ ( \s )^{-1} )$ iff $\tau
(\CJ (\s^{-1} )) = \CJ ( \s^{-1})$ for all $\s \in \S$  (since
$\s^{-1} (\CJ ( \s )^{-1}) = \CJ ( \s^{-1})$) iff $\tau (\l ) =
\l$ for all $\l \in {\rm im } (\CJ )$.  $\Box $

\begin{corollary}\label{n8Oct06}
Let $K$ be  a commutative ring and $n\geq 4$. Then the group $\S$
is not a normal subgroup of $\G$ iff $n\geq 5$.
\end{corollary}

{\it Proof}. For $n=4$, the group $\G$ is abelian, and so $\S$ is
a normal subgroup of $\G$. Let $n\geq 5$. By Theorem
\ref{c8Oct06}, $\S$ is a normal subgroup of $\G$ iff $\S = \CN (
\S , \G )$ iff
$${\rm im} (\CJ ) \subseteq {E'}_n^{\S} : =\{ e\in
E_n'\, | \, \s (e) = e, \; \forall \, \s \in \S \}.$$
 For the
automorphism $\G \ni \s : x_1\mapsto x_1( 1+x_2x_3), x_i\mapsto
x_i, i\neq 1$, the Jacobian $\CJ (\s ) = 1+x_2x_3$ does not belong
to ${E'}_n^{\S}$ since $\tau (\CJ (\s ))\neq \CJ (\s )$ where $\S
\ni \tau : x_2\mapsto x_2+x_1x_4x_5, x_j\mapsto x_j, j \neq 2$.
Therefore, $\S$ is not a normal subgroup of $\G$. $\Box$


\section{The algebraic group $\S'$ and its dimension}\label{GSID}

In this section, the group $\S'$ is studied in detail over a
commutative ring $K$. It is proved that the group $\S'$ is a
unipotent affine  group over $K$ of dimension $(n-2)2^{n-2}
-n+\pi_n$ (Corollary \ref{d25Oct06}). Important subgroups $\{
\Phi'^{2s+1}\}$ are introduced and results are proved for these
groups  (Lemma \ref{j23Oct06} and Theorem \ref{g23Oct06}) that
play a crucial role in finding the dimension and coordinates of
the Jacobian group $\S$.

\begin{lemma}\label{s24Oct06}
Let $K$ be an arbitrary ring, $n\geq 4$, and $s=1,2, \ldots ,
[\frac{n-1}{2}]$. Each element $a\in \L_{n, 2s}$ is a unique sum
$a= a_{n-2s}+a_{n-2s+1} +\cdots + a_n$ where
\begin{eqnarray*}
 a_{n-2s}&:=& \l x_{n-2s+1}x_{n-2s+2} \cdots x_n,\\
 a_{n-2s+p}&:=& c_p x_{n-2s+p+1}x_{n-2s+p+2} \cdots x_n,\;
 c_p:= \sum_{1\leq i_1<\cdots <i_p\leq n-2s+p-1} \l_{i_1, \ldots , i_p}x_{i_1} \cdots x_{i_p},  \\
 a_n&:=&c_{2s}:= \sum_{1\leq i_1<\cdots <i_{2s}\leq n-1} \l_{i_1, \ldots , i_{2s}}x_{i_1} \cdots x_{i_{2s}},
\end{eqnarray*}
$1\leq p \leq 2s-1$, and  the lambdas are from $K$.
\end{lemma}

{\it Proof}. This follows directly from Theorem \ref{14Sep06}.(1).
$\Box $

Let $K$ be a commutative ring and $n\geq 4$. For each {\em fixed}
natural number $s$ such that  $1\leq s\leq  [\frac{n-1}{2}]$, we
define the $K$-module
$$ V:=V_{n,2s}:= \L_{n,2s}(1) \oplus \cdots \oplus  \L_{n,2s}(n), \;\;
\L_{n,2s}(i):= K\lfloor x_1, \ldots , \widehat{x_i}, \ldots ,
x_n\rfloor_{2s}, $$ the direct sum of free $K$-modules of finite
rank over $K$. Each element $v=(v_1, \ldots , v_n) \in V$ is a
unique sum $v =\sum_{i=1}^nv_ie_i$ where $v_i\in \L_{n, 2s}(i) $
and $e_1:= (1, 0 , \ldots , 0), \ldots , e_n:=(0, \ldots , 0,1)$.
By Lemma \ref{s24Oct06}, the $K$-module homomorphism
\begin{equation}\label{VJ1}
\bCJ :=\bCJ_{n,2s} : V\ra \L_{n,2s}, \; (v_1, \ldots , v_n)\mapsto
v_1+\cdots + v_n,
\end{equation}
is a surjection, and the $K$-homomorphism (where $a= a_{n-2s}
+\cdots +a_n$ as in Lemma \ref{s24Oct06})
$$f=f_{n,2s} : \L_{n,2s} \ra V, \;\; a= a_{n-2s} +\cdots + a_n\mapsto
\sum_{i=n-2s}^na_ie_i= (0, \ldots , 0, a_{n-s} , \ldots , a_n) $$
is a section of $\bCJ$, i.e. $\bCJ f= {\rm id}$. Hence,
\begin{equation}\label{VJf}
V= \ker (\bCJ ) \oplus f(\L_{n,2s}), \;\; f(\L_{n,2s})= A_{n-2s}
\oplus \cdots \oplus A_n,
\end{equation}
where $A_i:= A_{n,2s, i}:=f(\L_{n, 2s} ) \cap \L_{n, 2s}(i)$.
 By Lemma \ref{s24Oct06},
 \begin{eqnarray*}
 A_{n-2s}&=& Kx_{n-2s+1}x_{n-2s+2} \cdots x_n, \\
 A_{n-2s+p}&=& (\bigoplus_{1\leq i_1<\cdots <i_p\leq n-2s+p-1} Kx_{i_1}\cdots x_{i_p})x_{n-2s+p+1}x_{n-2s+p+2} \cdots x_n, \\
 A_n&=& \bigoplus_{1\leq i_1<\cdots <i_{2s}\leq n-1} Kx_{i_1}\cdots x_{i_{2s}}, \\
\end{eqnarray*}
where $1\leq p \leq 2s-1$. So, the  $K$-modules $A_{n-2s}, \ldots
, A_n$ are free finitely generated $K$-modules  generated by
monomials (as above) of degree $2s$.

We will see later that the map $\bCJ$ in (\ref{VJ1}) is, up to
isomorphism, the Jacobian map $\bCJ $ (\ref{VJ2}).  Our goal is to
find a special $K$-basis for the kernel $\ker (\bCJ )$ that has
connection with certain automorphisms of the group $\S'$. For, we
will define, so-called, avoidance  functions  which allows one to
produce the required basis and then explicit automorphisms of  the
group  $\S'$.

{\bf Avoidance functions}. For each monomial $u= x_{i_1} \cdots
x_{i_t}$ of $\L_n$ the set $\{ i_1, \ldots , i_t\}$ is called the
{\em support} of the monomial $u$. Let us fix a natural number $s$
such that $1\leq s \leq [\frac{n-1}{2}]$. Next, for each $i=1,
\ldots , n-1$, we are going to define a set $S_i'$ and a function
$j_i$ on it. We do this in two steps: first, for $i=n-2s, \ldots ,
n-1$; and then for $i=1, \ldots , n-2s-1$.

 For each $i=n-2s, n-2s+1,
\ldots , n$, let $S_i:=S_{i,s}:= \Supp (A_i)$ be the set of
supports of all the monomials from the module $A_i$ (the
$K$-module $A_i$ is generated by monomials), and let
$S_i':=S_{i,s}'$ be its  complement in the set $\Supp (\L_{n,
2s}(i))= \{ \alpha \subseteq \{ 1, \ldots , \widehat{i}, \ldots ,
n\} \, | \, |\alpha | =2s\}$ where $|\alpha |$ is the number of
elements in the set $\alpha $. In more detail,
\begin{eqnarray*}
 S_{n-2s, s}&=& \{ \{ n-2s+1, \ldots , n \} \} , \\
S_{n-2s+p, s}&=& \{ \alpha \subseteq \{ 1, \ldots
,\widehat{n-2s+p}, \ldots , n \}\, | \,
 \alpha \supseteq \{ n-2s+p+1, \ldots , n\} , |\alpha | =2s \} , \\
 S_{n, s}&=& \{ \alpha \subseteq \{ 1, \ldots , n-1 \}\, | \,
 |\alpha | =2s \} , \; {\rm and} \\
  S_{n-2s, s}'&=&\{ \alpha \subseteq \{ 1, \ldots ,\widehat{n-2s}, \ldots ,  n \}\, | \,
   \alpha \neq \{ n-2s+1, \ldots , n\} , |\alpha | =2s \} , \\
 S_{n-2s+p, s}'&=& \{ \alpha \subseteq \{ 1, \ldots ,\widehat{n-2s+p}, \ldots  ,n \}\, | \,
 \alpha \not\supseteq \{ n-2s+p+1, \ldots , n\} , |\alpha | =2s \} , \\
 S_{n, s}'&=&\emptyset ,
\end{eqnarray*}
where $1\leq p \leq 2s-1$. One can easily verify that the
following sets are {\em the only empty sets} among the sets $ \{
S_{n-2s+p', s}'\, | \, s=1, \ldots , [\frac{n-1}{2}]; \,
p'=0,1,\ldots , 2s-1\}$: if $n$ is an {\em odd} number,
$s=[\frac{n-1}{2}]$, and $p'=0,1,\ldots , 2s-1$, i.e.
\begin{equation}\label{Sem}
S_{1,[\frac{n-1}{2}]}'=S_{2,[\frac{n-1}{2}]}'=\cdots = S_{n-1,
[\frac{n-1}{2}]}'=\emptyset .
\end{equation}
In particular, for $n=5$ we have 
\begin{equation}\label{Sem1}
S_{1,2}'=S_{2,2}'=\cdots = S_{4,2}'=\emptyset .
\end{equation}

 Let us stress that {\em for each $i=n-2s,
\ldots , n-1$, the set $S_i':=S_{i,s}'$ is equal to the set of all
$\alpha \in \Supp (\L_{n, 2s}(i))$ such that} $\{ i+1, i+2, \ldots
, n\} \backslash \alpha \neq \emptyset$. So, we can fix a function
$$ j_i:=j_{i,s}: S_i'\ra \{ i+1, i+2, \ldots , n\} , \;\; \alpha \mapsto
j_i(\alpha )\in \{ i+1, i+2, \ldots , n\} \backslash \alpha .$$ If
the set $S_i'$ is an empty set then this definition is vacuous
since we have an `empty function' defined on the empty set. It is
convenient to have these  `empty functions' in order to save on
notation.

  For each $i=1, \ldots , n-2s-1$,   let $S_i':=S_{i,s}':=\Supp
(\L_{n,2s}(i))$, and  we can fix a function
$$ j_i:=j_{i,s}: S_i'\ra \{ i+1, i+2, \ldots , n\} , \;\; \alpha \mapsto
j_i(\alpha )\in \{ i+1, i+2, \ldots , n\} \backslash \alpha ,$$
(if not then $\{ i+1, \ldots , n\} \subseteq \alpha$ for some
$\alpha\in S_i'$, and so $\alpha \cup \{ 1, \ldots , i\}= \{ 1,
\ldots , n\}$ but then one has the contradiction: $n= |\alpha \cup
\{ 1, \ldots , i\}|\leq |\alpha | +i\leq 2s +n-2s-1= n-1<n$).

{\it Definition}. For the fixed number $s$ (as above),   the
functions $\{ j_i\, | \, 1\leq i \leq n-1\}$ are called {\bf
avoidance functions}.

Note that there are  many avoidance functions, in general.  The
importance of avoidance functions $\{ j_i\}$ is the fact that, for
any $\alpha \in S_i'$,  we can attach the 1-dimensional abelian
subgroup of $\S'$, $\{ \rho_{i, j_i(\alpha ) ; \l x^\alpha
}\}\simeq K$, $\rho_{i, j_i(\alpha ) ; \l x^\alpha }\mapsto \l$,
by the rule 
\begin{equation}\label{rjet}
\rho_{i, j_i(\alpha ) ; \l x^\alpha }:x_i\mapsto x_i(1+\l x^\alpha
) , \;\; x_{j_i(\alpha )} \mapsto x_{j_i(\alpha )}(1-\l x^\alpha
), \; \; x_k\mapsto x_k, \;\; k\neq i, j_i(\alpha ).
\end{equation}

\begin{lemma}\label{k24Oct06}
Let $K$ be a commutative ring, $n\geq 4$, and $\{ j_i=j_{i,s}\}$
be avoidance functions for a fixed number $s$ such that $1\leq
s\leq [\frac{n-1}{2}]$. Then the  set $\cup_{i=1}^{n-1} \{
x^\alpha (e_i-e_{j_i(\alpha )})\, | \, \alpha \in S_i'\}$ is a
$K$-basis for the kernel $\ker (\bCJ )$ of the map $\bCJ
=\bCJ_{n,2s}$, (\ref{VJ1}); and the rank of the free $K$-module
$\ker (\bCJ_{n,2s} )$ is equal to
$$ {\rm rk}_K(\ker (\bCJ_{n,2s} ))=\sum_{i=1}^{n-1}|S_{i,s}'|=
n{n-1\choose 2s}-{n\choose 2s}.$$
\end{lemma}

{\it Proof}. By the very definition, the elements from the union
are from the kernel $\ker (\bCJ )$ and   $K$-linear independent
(use the fact that $i<j_i(\alpha )$ for all $i$ and $\alpha \in
S_i'$; and the fact that, for each $i=1, \ldots , n-1$, the
monomials $x^\alpha$, $\alpha \in S_{i,s}'$ are $K$-linear
independent). Let $U$ be the $K$-submodule of $V$ that these
elements generate. It follows from the definition of the sets
$S_i'$ that $U+f(\L_{n, 2s})=V$, hence $U= \ker (\bCJ )$, by
(\ref{VJf}) and the inclusion $U\subseteq \ker (\bCJ )$.

By (\ref{VJf}), the rank of the free $K$-module $\ker
(\bCJ_{n,2s})$ is equal to
 \begin{eqnarray*}
{\rm rk}_K(\ker (\bCJ_{n,2s} ))&=&\sum_{i=1}^{n-1}|S_{i,s}'|={\rm
rk}_K(V_{n,2s})-{\rm rk}_K(f(\L_{n,2s}))\\
&=& n \, {\rm rk}_K(\L_{n-1,2s})- {\rm rk}_K(\L_{n,2s})=
n{n-1\choose 2s}-{n\choose 2s}.\;\; \Box
\end{eqnarray*}

{\bf The groups $\Phi'^{2s+1}$}. Let $K$ be a commutative ring,
and  $n\geq 4$. For each number $s=1,2, \ldots , [\frac{n-1}{2}]$,
let  $\Phi'^{2s+1}$ be the subset  of $\Phi$ that contains all the
elements  of the following type:
\begin{eqnarray*}
 \s(x_i)&=& x_i(1+\cdots ), \;\; 1\leq  i\leq  n-2s-1,  \\
\s(x_{n-2s})&=& x_{n-2s}(1+\l x_{n-2s+1} x_{n-2s+2}\cdots
x_n+\cdots ),\\
\s(x_{n-2s+1})&=& x_{n-2s+1}(1+(\sum_{1\leq i_1\leq n-2s} \l_{i_1}
x_{i_1})  x_{n-2s+2}\cdots
x_n+\cdots ),\\
\cdots &=& \cdots ,\\
 \s(x_{n-2s+p})&=& x_{n-2s+p}(1+(\sum_{1\leq
i_1<\cdots <i_p\leq n-2s+p-1} \l_{i_1, \ldots , i_p} x_{i_1}\cdots
x_{i_p}) x_{n-2s+p+1}\cdots
x_n+\cdots ),\\
\cdots &=& \cdots ,\\
\s(x_{n-1})&=& x_{n-1}(1+(\sum_{1\leq i_1<\cdots <i_{2s-1}\leq
n-2} \l_{i_1, \ldots , i_{2s-1}} x_{i_1}\cdots x_{i_{2s-1}})
x_n+\cdots ),\\
\s(x_n)&=& x_n(1+\sum_{1\leq i_1<\cdots <i_{2s}\leq n-1} \l_{i_1,
\ldots , i_{2s}} x_{i_1}\cdots x_{i_{2s}}+\cdots ),
\end{eqnarray*}
where the lambdas  are elements of $K$ and the three dots mean
higher terms. By Theorem \ref{11Sep06}.(2), $\Phi'^{2s+1}$ is a
subgroup of $\Phi$. In the notation of Lemma \ref{s24Oct06}, the
automorphism $\s =\s_a$ above (where $a=a_{n-2s}+\cdots +a_n$ as
in Lemma \ref{s24Oct06})  can be written as
\begin{eqnarray*}
 \s (x_i)&=&x_i(1+\cdots ), \;\; 1\leq i \leq n-2s-1,  \\
 \s (x_i)&=&x_i(1+a_i+\cdots ), \;\; n-2s\leq i \leq n.
\end{eqnarray*}
Clearly,
$$ \Phi^{2s+3}\subseteq \Phi'^{2s+1} \subseteq \Phi^{2s+1}, $$
and $\Phi^{2s+3}$ is a normal subgroup of $\Phi'^{2s+1}$ since
$\Phi^{2s+3}$ is a normal subgroup of $\Phi$.  For each $s=1,
\ldots , [\frac{n-1}{2}]$, we have the Jacobian map
$$ \CJ : \Phi^{2s+1} \ra E_{n,2s}':= 1+\sum_{i\geq  s} \L_{n,
2s}, \;\; \s \mapsto \CJ (\s ).$$ The factor group $ E_{n,2s}'/
E_{n,2s+2}'= \{ (1+a)E_{n,2s+2}'\, | \, a\in \L_{n, 2s}\}$ is
canonically isomorphic to the additive group $\L_{n, 2s}$ via the
isomorphism
$$ E_{n,2s}'/E_{n,2s+2}'\ra \L_{n,2s}, \;\;
(1+a)E_{n,2s+2}'\mapsto a.$$ The Jacobian map $\CJ : \Phi^{2s+1}
\ra E_{n, 2s}'$ yields the Jacobian maps 
\begin{equation}\label{VJ2}
\bCJ : \Phi^{2s+1} / \Phi^{2s+3} \ra E_{n, 2s}'/E_{n, 2s+2}', \;\;
\s\Phi^{2s+3} \mapsto \CJ (\s ) E_{n, 2s+2}',
\end{equation}
and
$$
\bCJ : \Phi'^{2s+1} / \Phi^{2s+3} \ra E_{n, 2s}'/E_{n, 2s+2}',
\;\; \s\Phi^{2s+3} \mapsto \CJ (\s ) E_{n, 2s+2}'.
$$
There is the  natural  isomorphism of the abelian groups:
\begin{eqnarray*}
 \Phi^{2s+1}/\Phi^{2s+3} &\ra &
\bigoplus_{i=1}^n \L_{n,2s}(i)=V=V_{n,2s},\\
\{ \s : x_i\mapsto x_i(1+a_i) \Phi^{2s+3} \} &\mapsto & (a_1,
\ldots , a_n),
\end{eqnarray*}
where $a_i\in \L_{n,2s}(i)$ for all $i=1, \ldots , n$.  When we
identify the groups $\Phi^{2s+1}/\Phi^{2s+3}$ and $V$ on the one
hand, and the groups $E_{n,2s}'/E_{n,2s+2}'$ and $\L_{n,2s}$ on
the other via the isomorphisms above, then {\em the Jacobian map}
$\bCJ : \Phi^{2s+1}/\Phi^{2s+3}\ra E_{n,2s}'/E_{n,2s+2}'$,
$\s\Phi^{2s+3}\mapsto \CJ (\s )E_{n,2s+2}'$,  {\em coincides with
the map} (\ref{VJ1}), $\bCJ : V\ra \L_{n,2s}$, $(a_1, \ldots ,
a_n) \mapsto a_1+\cdots +a_n$. Then  Lemma \ref{j23Oct06} follows,
which is one of the key results in finding generators for the
Jacobian group $\S$ and its subgroup $\S'$.

\begin{lemma}\label{j23Oct06}
Let $K$ be a commutative ring, $n\geq 4$, and $s=1, \ldots ,
[\frac{n-1}{2}]$. The Jacobian map $\bCJ :
\Phi'^{2s+1}/\Phi^{2s+3} \ra E_{n,2s}'/E_{n,2s+2}'$,
$\s\Phi^{2s+3}\mapsto \CJ (\s )E_{n,2s+2}'$,  is an isomorphism of
the abelian groups which is given by the rule
\begin{eqnarray*}
\bCJ (\s \Phi^{2s+3}) &=&(1+\l x_{n-2s+1} \cdots x_n
+\sum_{p=1}^{2s-1}(\sum_{1\leq i_1<\cdots <i_s\leq n-2s+p-1}
\l_{i_1, \ldots , i_p} x_{i_1}\cdots x_{i_p})x_{n-2s+p+1} \cdots
x_n\\
&+& \sum_{1\leq i_1<\cdots <i_{2s}\leq n-1} \l_{i_1, \ldots , i_p}
x_{i_1}\cdots x_{i_{2s}})E_{n,2s+2}'
\end{eqnarray*}
for the  element $\s\in \Phi'^{2s+1}$ as above (i.e. in the
definition of $\Phi'^{2s+1}$).
\end{lemma}

{\it Proof}. When one writes down the determinant $\CJ (\s )$ for
the  element $\s \in \Phi'^{2s+1}$ as above it is easy to see that
$\CJ (\s )E_{n,2s+2}'$ is the product of the diagonal elements in
the determinant $\CJ (\s )$ modulo $E_{n,2s+2}'$:
\begin{eqnarray*}
\bCJ (\s \Phi^{2s+3}) &=&(1+\l x_{n-2s+1}\cdots x_n)\cdots\\
&\cdots & (1+\sum_{1\leq i_1<\cdots <i_s\leq n-2s+p-1} \l_{i_1,
\ldots ,
i_p} x_{i_1}\cdots x_{i_p})x_{n-2s+p+1}\cdots x_n)\cdots \\
& \cdots &(1+\sum_{1\leq i_1<\cdots <i_{2s}\leq n-1} \l_{i_1,
\ldots , i_p} x_{i_1}\cdots x_{i_{2s}})E_{n,2s+2}'\\
&=&(1+\l x_{n-2s+1} \cdots x_n +\sum_{p=1}^{2s-1}(\sum_{1\leq
i_1<\cdots <i_s\leq n-2s+p-1} \l_{i_1, \ldots , i_p} x_{i_1}\cdots
x_{i_p})x_{n-2s+p+1} \cdots
x_n\\
&+& \sum_{1\leq i_1<\cdots <i_{2s}\leq n-1} \l_{i_1, \ldots , i_p}
x_{i_1}\cdots x_{i_{2s}})E_{n,2s+2}',
\end{eqnarray*}
and so we obtain  the formula for $\bCJ (\s \Phi^{2s+3})$ in Lemma
\ref{j23Oct06}. Now, it is obvious that the map $\bCJ $ is the
isomorphism of the abelian groups since each element of $\L_{n,
2s}$ can be uniquely written as a sum $s$ in the formula for $\bCJ
(\s \Phi^{2s+3})= (1+s)E_{n,2s+2}'$ above (see Lemma
\ref{s24Oct06}). $\Box $

\begin{theorem}\label{g23Oct06}
Let $K$ be a commutative ring, $n\geq 4$, and $s=1, \ldots ,
[\frac{n-1}{2}]$. Then
\begin{enumerate}
\item $\Phi^{2s+1}= \Phi'^{2s+1}\S'^{2s+1} = \S'^{2s+1}
\Phi'^{2s+1}$.\item $\Phi^{2s+1}/\Phi^{2s+3} \simeq
\S'^{2s+1}/\S'^{2s+3} \times \Phi'^{2s+1}/\Phi^{2s+3}$, the direct
product of abelian groups. \item Each automorphism $\s \in
\Phi^{2s+1}$: $x_i \mapsto x_i(1+b_i+\cdots )$, $b_i\in \L_{n, 2s}
(i)$, $i=1, \ldots , n$, is the unique product modulo
$\Phi^{2s+3}$ as follows, 
\begin{equation}\label{srp1}
\s\equiv \phi'\prod_{i=1}^{n-1} \prod_{\alpha \in S_i'}
\rho_{i,j_i(\alpha ); \l_\alpha x^\alpha} \mod \Phi^{2s+3},
\end{equation}
where
\begin{eqnarray*}
 \phi'(x_k)&:=&x_k, \;\; 1\leq k \leq n-2s-1, \\
 \phi'(x_i)&:=&x_i(1+a_i), \;\; n-2s\leq i \leq n,
\end{eqnarray*}
and
$$ (b_1, \ldots , b_n)= \sum_{i=1}^{n-1}\sum_{\alpha\in S_i'}
\l_\alpha x^\alpha (e_i-e_{j_i(\alpha )})+\sum_{i=n-2s}^na_ie_i$$
in $V= \ker (\bCJ )\oplus f(\L_{n, 2s})$ for unique $\l_\alpha \in
K$ and unique elements $a_i$ as in Lemma \ref{s24Oct06}.
\end{enumerate}
\end{theorem}

{\it Proof}. Recall that the Jacobian map $\bCJ
:\Phi^{2s+1}/\Phi^{2s+3} \ra E_{n, 2s}'/ E_{n, 2s+2}'$,
$\s\Phi^{2s+3}\mapsto \CJ (\s ) E_{n, 2s+2}'$, is naturally
identified with the map (\ref{VJ1}),
$$ \bCJ : V\ra \L_{n, 2s}, \;\; (a_1, \ldots , a_n)\mapsto
a_1+\cdots + a_n,$$ under the identifications
$\Phi^{2s+1}/\Phi^{2s+3}\equiv V$ and $ E_{n, 2s}'/ E_{n,
2s+2}'\equiv \L_{n, 2s}$. Consider the direct sum (\ref{VJf}): $V=
\ker (\bCJ ) \oplus f(\L_{n, 2s})$. By Lemma \ref{k24Oct06}, the
free $K$-module $\ker (\bCJ )$ has the $K$-basis $\cup_{i=1}^{n-1}
\{ x^\alpha (e_i-e_{j_i(\alpha )})\, | \, \alpha \in S_i'\}$, and
each element of $f(\L_{n, 2s})$ is a unique sum
$\sum_{i=n-2s}^na_ie_i$ where $a_i$ are as in Lemma
\ref{s24Oct06}. To each basis element $x^\alpha (e_i-e_{j_i(\alpha
)})$ we attach the 1-dimensional abelian subgroup of $\S'^{2s+1}$:
$ \{ \rho_{i,j_i(\alpha ) ; \l x^\alpha}\}_{\l \in K} \simeq K$,
by (\ref{rjet}). To each element $a= \sum_{i=n-2s}^na_ie_i$ it
corresponds (under the identification
$\Phi^{2s+1}/\Phi^{2s+3}\equiv V$) the automorphism $\phi_a\in
\Phi'^{2s+1}$:
\begin{eqnarray*}
 \phi_a(x_k)&:=&x_k, \;\; 1\leq k \leq n-2s-1, \\
 \phi_a(x_i)&:=&x_i(1+a_i), \;\; n-2s\leq i \leq n.
\end{eqnarray*}
Each element $v$ of $V$ is a unique sum
$$ v= \sum_{i=1}^{n-1}\sum_{\alpha \in S_i'} \l_\alpha x^\alpha
(e_i-e_{j_i(\alpha )})+\sum_{k=n-2s}^na_ke_k$$ for unique
$\l_\alpha \in K$ and unique $a_k$ as in Lemma \ref{s24Oct06}.
Under the identification $\Phi^{2s+1}/\Phi^{2s+3}\equiv V$, the
element $v$ can be identified with the automorphism $\s_v$ modulo
$\Phi^{2s+3}$ (i.e. $v\equiv \s_v\Phi^{2s+3}$) where
\begin{equation}\label{svas1}
\s_v = \phi_a\s', \;\; \s':= \prod_{i=1}^{n-1} \prod_{\alpha \in
S_i'}\rho_{i,j_i(\alpha ); \l_\alpha x^\alpha}.
\end{equation}
Conversely, any coset $\s\Phi^{2s+3}$ where $\s \in \Phi^{2s+1}$
can be identified with the element $v\in V$ (i.e.
$\s\Phi^{2s+3}\equiv v$) by the rule: let $\s(x_i)= x_i(1+
b_i+\cdots ) $ for some $b_i\in \L_{n,2s}(i)$, $i=1, \ldots , n$,
then $(b_1, \ldots , b_n)\in V= \oplus_{i=1}^n \L_{n,2s}(i)$ and
$\s\Phi^{2s+3} \equiv (b_1, \ldots , b_n)$. Now, statement 1
follows immediately from (\ref{svas1}). Under the identification
$\Phi^{2s+1}/\Phi^{2s+3}\equiv V$, the decomposition $V= \ker
(\bCJ ) \oplus f(\L_{n, 2s})$ corresponds to the decomposition
(the direct product of groups)
$$\Phi^{2s+1}/\Phi^{2s+3} \simeq
\S'^{2s+1}\Phi^{2s+3}/\Phi^{2s+3} \times
\Phi'^{2s+1}/\Phi^{2s+3}.$$ Since $
\S'^{2s+1}\Phi^{2s+3}/\Phi^{2s+3}\simeq \S'^{2s+1}/\S'^{2s+3}\cap
\Phi^{2s+3}\simeq \S'^{2s+1}/\S'^{2s+3}$, the statement 2 follows.
Statement 3 is just (\ref{svas1}).
 $\Box $

\begin{theorem}\label{23Oct06}
Let $K$ be a commutative ring, and  $n\geq 4$.
\begin{enumerate}
\item Then each automorphism $\s \in \Phi$ is a unique product
\begin{equation}\label{srp2}
\s = \prod_{i=1}^{[\frac{n-1}{2}]}\phi_{2s+1}\s_{2s+1}=
\phi_3\s_3\phi_5\s_5\cdots
\phi_{2[\frac{n-1}{2}]+1}\s_{2[\frac{n-1}{2}]+1}
\end{equation}
for unique elements $\phi_{2s+1}\in \Phi'^{2s+1}$ and
$\s_{2s+1}\in \S'^{2s+1} $ from   (\ref{srp1}).  Moreover,
\begin{eqnarray*}
 \s &\equiv &\phi_3\s_3\mod \Phi^5, \\
(\phi_3\s_3\cdots \phi_{2s-1}\s_{2s-1})^{-1} \s &\equiv
&\phi_{2s+1}\s_{2s+1}\mod \Phi^{2s+3}, \;\; 2\leq s\leq
[\frac{n-1}{2}].
\end{eqnarray*}
\item Each automorphism $\s \in \S'$ is a unique product
\begin{equation}\label{srp3}
\s = \prod_{i=1}^{[\frac{n-1}{2}]}\s_{2s+1}= \s_3\s_5\cdots
\s_{2[\frac{n-1}{2}]+1}
\end{equation}
for unique elements $\s_{2s+1}\in \S'^{2s+1} $ from (\ref{srp1}).
Moreover,
\begin{eqnarray*}
 \s &\equiv &\s_3\mod \Phi^5, \\
(\s_3\cdots \s_{2s-1})^{-1} \s &\equiv &\s_{2s+1}\mod \Phi^{2s+3},
\;\; 2\leq s\leq [\frac{n-1}{2}].
\end{eqnarray*}
\end{enumerate}
\end{theorem}

{\it Proof}. 1. This statement follows from Theorem
\ref{g23Oct06}.

2. We need only to show that, for $\s \in \S'$,  $\phi_3= \cdots =
\phi_{2[\frac{n-1}{2}]+1}=e$ in (\ref{srp2}). Suppose that
$\phi_{2s+1} \neq e$ for some $s$ and the $s$ is the least
possible with this property. We seek a contradiction. Without loss
of generality we may assume that $\s_3=\cdots = \s_{2s-1}=e$, i.e.
$\s = \phi_{2s+1}\s_{2s+1} \cdots$. Since $\CJ (\s ) =1$, we must
have $\bCJ (\s )=1$ in $E_{n, 2s}'/ E_{n, 2s+2}'$. On the other
hand, $\bCJ (\s ) = \bCJ (\phi_{2s+1})\neq 1$ in $E_{n, 2s}'/
E_{n, 2s+2}'$ (Lemma \ref{j23Oct06}), by the choice of the
$\phi_{2s+1}$, a contradiction. $\Box$

{\bf The dimension of the algebraic group $\S'$}. Recall that
$n\geq 4$ and for each number $s=1, \ldots , [\frac{n-1}{2}]$, we
defined the sets $S_i':= S_{i, s}'$, $1\leq i \leq n-1$,  and the
avoidance functions $\{ j_i:= j_{i,s}\}$. By Theorem
\ref{23Oct06}.(2) and Theorem \ref{g23Oct06}.(3), each element
$\s$ of $\S'$ is a unique {\em ordered} product 
\begin{equation}\label{scooS}
 \s = \prod_{s=1}^{[\frac{n-1}{2}]}\prod_{i=1}^{n-1}\prod_{\alpha
\in S_{i, s}'}\rho_{i,j_{i,s}(\alpha ); \l_\alpha x^\alpha }
\end{equation}
where $\alpha = \alpha_{i,s}$ (they depend on $i$ and $s$) and
$\l_{\alpha} = \l_{\alpha, i,s}\in K$. Therefore, $\{ \l_\alpha =
\l_{\alpha, i,s} \}$ are affine coordinates for the algebraic
group $\S'$ over the ring $K$, and the algebra of (regular)
functions $\OO (\S')$ on the algebraic group $\S'$ is a {\em
polynomial algebra} in 
\begin{equation}\label{dimS1}
\dim (\S') =\sum_{s=1}^{[\frac{n-1}{2}]}\sum_{i=1}^{n-1}|S_{i,s}'|
\end{equation}
variables. Consider the function $$\pi_n:=\begin{cases}
1& \text{if $n$ is odd},\\
2& \text{if $n$ is even}.\\
\end{cases}$$

\begin{corollary}\label{d25Oct06}
Let $K$ be a commutative ring and $n\geq 4$. The group $\S'$ is a
unipotent affine  group over $K$ of dimension
\begin{eqnarray*}
 \dim (\S') &=&\sum_{s=1}^{[\frac{n-1}{2}]}(n{n-1\choose 2s} -
{n\choose 2s})= \sum_{s=1}^{[\frac{n-1}{2}]}(n-2s-1){n\choose
2s}\\
&=& \begin{cases}
(n-2)2^{n-2}-n+2& \text{if $n$ is even},\\
(n-2)2^{n-2}-n+1& \text{if $n$ is odd}.\\
\end{cases}
\end{eqnarray*}
over $K$, i.e.  the algebra of regular functions on $\S'$ is a
polynomial algebra over the ring $K$ in $\dim (\S')$ variables $\{
\l_\alpha \}$.
\end{corollary}

{\it Proof}. The only statement which is needed to be proven is
the formula for the dimension. For each $s$, by Lemma
\ref{k24Oct06},
\begin{eqnarray*}
\sum_{i=1}^{n-1}|S_{i,s}'| &=& n{n-1\choose 2s} - {n\choose 2s}=
(n-2s-1){n\choose 2s}.
\end{eqnarray*}
The first part of the formula for $\dim (\S')$ then follows from
(\ref{dimS1}). Note that
$$ n\sum_{s=1}^{[\frac{n-1}{2}]}{n-1\choose 2s}=n(\sum_{s=0}^{[\frac{n-1}{2}]}{n-1\choose
2s}-1)=n(2^{n-2}-1).$$ If $n$ is  even then
$$\sum_{s=1}^{[\frac{n-1}{2}]}{n\choose 2s}=\sum_{s=0}^{[\frac{n-1}{2}]}{n\choose 2s}-1=
\sum_{s=0}^{[\frac{n}{2}]}{n\choose 2s}-{n\choose
2[\frac{n}{2}]}-1=2^{n-1}-2.$$ If $n$ is odd then
$$\sum_{s=1}^{[\frac{n-1}{2}]}{n\choose 2s}=\sum_{s=1}^{[\frac{n}{2}]}{n\choose 2s}=
\sum_{s=0}^{[\frac{n}{2}]}{n\choose 2s}-1=2^{n-1}-1.$$ By the
first part of the formula for $\dim (\S' )$ and the calculations
above, we have
$$ \dim (\S' ) = n(2^{n-2}-1) -(2^{n-1}-\pi_n) = (n-2)2^{n-2}
-n+\pi_n.  \;\;\;\; \Box $$ The subgroup $\S'$ of $\G$ is `twice
smaller' than $\G$ in the following sense (see (\ref{dimG}))
\begin{equation}\label{dSGh}
\lim_{n\rightarrow \infty} \frac{\dim (\S')}{\dim (\G )}=
\lim_{n\rightarrow \infty}
\frac{(n-2)2^{n-2}-n+\pi_n}{n(2^{n-1}-n)}=\frac{1}{2}.
\end{equation}

{\bf The group $\S'$ is not a normal subgroup of $\S$ if $n\geq 6$
and $2\neq 0$ in $K$}. Clearly, $\S\ni \s : x_1\mapsto
x_1+x_2x_3x_4, \; x_i\mapsto x_i, \; i\neq 1$; and $\S'\ni \tau :
x_1\mapsto x_1(1+x_5x_6),  x_2\mapsto x_2(1-x_5x_6), \; x_i\mapsto
x_i, \; i\neq 1,2$. Then $\s \tau \s^{-1} \tau^{-1}(x_1) =
x_1+2x_2x_3x_4x_5x_6$, hence $\s \tau \s^{-1} \tau^{-1} \not\in
\S'$. This means that the subgroup $\S'$ of $\S$ is not normal if
$n\geq 6$ and $2\neq 0$ in $K$.

We will see later that the group $\S''$ is a {\em closed normal}
subgroup of the Jacobian group $\S$ (Theorem \ref{20Nov06}.(2)).


\section{A (minimal) set of generators for the Jacobian group $\S$ and its
dimension }\label{MSGJG}

Let $K$ be a {\em commutative} ring. In this section, a minimal
set of generators for the Jacobian group $\S$ is found explicitly
(Theorem \ref{14Nov06}). The dimensions and coordinates of   the
following algebraic groups are found explicitly: $\S $ (Theorem
\ref{19Nov06}), $\S'\cap \S''$ (Lemma \ref{i14Nov06}), $\S''$
(Theorem \ref{20Nov06}). It is proved that the sets of cosets $\S'
\diagdown \S$ and $ \S'\cap \S'' \diagdown \S''$ have natural
structure of an affine variety of dimension $n(2^{n-2}-n+1)$ over
$K$ (Corollary \ref{s19Nov06}).

{\bf (Minimal) set of generators for the Jacobian group $\S $}. We
keep the notations of Section \ref{GSID}. For $s=1$, $1\leq p \leq
2s-1=1$ , i.e. $p=1$. Consider the sets $S_{i,1}'$, $1\leq i \leq
n-1$, defined in Section \ref{GSID}:
\begin{eqnarray*}
 S_{i,1}'&=&\Supp (\L_{n,2}(i)),\;\;  1\leq i\leq n-3,  \\
 S_{n-2,1}'&=& \{ \alpha \subseteq \{ 1, \ldots , \widehat{n-2},
 n-1, n\, | \, \alpha \neq \{ n-1, n\}, |\alpha |=2\}, \\
 S_{n-1,1}'&=& \{ \alpha \subseteq \{ 1, \ldots , \widehat{n-1},
  n\, | \, \alpha \not\ni n, |\alpha |=2\}.
\end{eqnarray*}
Fix avoidance functions
$$ j_i: S_{i,1}'\mapsto \{ i+1, \ldots , n\} , \;\; 1\leq i \leq
n-1. $$

The next theorem provides a (minimal) set of generators for the
Jacobian group $\S$ for $n\geq 7$.
\begin{theorem}\label{14Nov06}
Let $K$ be a commutative ring, $n\geq 7$, and for $s=1$ let $\{
j_i:=j_{i,1}\}$ be avoidance functions. If either $n$ is odd; or
$n$ is even and $\frac{1}{2}\in K$; then
$$ \S = \langle \rho_{i,j_i(\alpha ) ; \l x^\alpha}, \s_{i', \l
x_jx_kx_l}\, | \, \l \in K; 1\leq i\leq n-1; \alpha \in S_{i,1}';
1\leq i'\leq n; j<k<l; i'\not\in \{ j,k,l\}\rangle$$ and the
1-dimensional abelian subgroups $\{\rho_{i,j_i(\alpha ) ; \l
x^\alpha}\} \simeq K$ and $\{ \s_{i', \l x_jx_kx_l} \}\simeq K$ of
$\S$ form a minimal set of generators for $\S$ in the sense that
no subgroup can  be dropped.
\end{theorem}

{\it Proof}. Recall that $\S = \S'\S''$ (Corollary
\ref{a8Oct06}.(1)); and, by (\ref{scooS}),
$$ \S'=\langle  \rho_{i,j_{i, s}(\alpha_s ) ; \l x^{\alpha_s}}\, | \,
\l \in K; 1\leq i\leq n-1; 1\leq s \leq [\frac{n-1}{2}]; \alpha_s
\in S_{i,s}' \rangle$$  where $j_{i, s}$ are avoidance functions;
$$ \S''=\langle  \s_{i', \l x_jx_kx_l}\, | \, \l \in K;
1\leq i'\leq n; j<k<l; i'\not\in \{ j,k,l\}\rangle$$ (Theorem
\ref{s15Oct06}.(1)), and, by the definition, the group $\S''$ is
 generated by all the automorphisms
 $$ \xi_{i,b_i}:x_i\mapsto x_i+b_i, \;\; x_j\mapsto x_j, \;\;
 j\neq i,$$
 where $b_i\in K\lfloor x_1, \ldots , \widehat{x_i}, \ldots ,
 x_n\rfloor^{od}_{\geq 3}$.

 If $s=1$ then the elements $\{ \rho_{i,j_{i,1}(\alpha_1); \l
 x^{\alpha_1}}\}$ are precisely the $\rho$-part of the generators
 in the theorem. If $s\geq 2$ then, by (\ref{com1}), each element
 $\rho_{i,j_{i,s}(\alpha_s); \l x^{\alpha_s}}$  belongs to the group
 $\S''$. Therefore,
$$ \S = \langle \rho_{i,j_i(\alpha ) ; \l x^\alpha}, \s_{i', \l
x_jx_kx_l}\, | \, \l \in K; 1\leq i\leq n-1; \alpha \in S_{i,1}';
1\leq i'\leq n; j<k<l; i'\not\in \{ j,k,l\}\rangle ,$$ i.e. the
first part of the theorem is proved. To prove the second part of
the theorem (about minimality), note that, by Theorem
\ref{11Sep06}.(2),  the map
$$(\L_{n,3})^n\ra U^3/U^5, \;\; a=(a_1, \ldots , a_n)\mapsto
\s_aU^5,$$ is a group isomorphism where $\s_a\in U^3: x_i\mapsto
x_i+a_i$. We identify these two groups via the isomorphism above,
then the elements $\rho_{i,j_i(\alpha ); \l x^\alpha }U^5$ and
$\s_{i'; \l x_jx_kx_l}U^5$ are identified correspondingly with the
elements $\l x^\alpha (x_ie_i-x_{j_i(\alpha )}e_{j_i(\alpha )})$
and $ \l x_jx_kx_le_{i'}$ of the $K$-module $W:= (\L_{n,3})^n=
\oplus_{i=1}^n \L_{n,3}e_i$ where $e_1:=(1,0, \ldots , 0), \ldots
, e_n:=(0, \ldots, 0,1)$. To prove the minimality it suffices to
show that the elements $ \{ x^\alpha (x_ie_i-x_{j_i(\alpha
)}e_{j_i(\alpha )}), \;   x_jx_kx_le_{i'} \}$ are $K$-linearly
independent. To prove this let us consider the descending
filtration $\{ W_i:= \oplus_{j=i}^n \L_{n,3}e_j\}$ on $W$.
Clearly, $W_i/W_{i+1} = (\L_{n,3}e_i\oplus W_{i+1} )/ W_{i+1}
\simeq \L_{n,3}e_i\simeq \L_{n,3}$, $i\geq 1$. Suppose that $r:=
\sum \l_{i\alpha} x^\alpha (x_ie_i-x_{j_i(\alpha )}e_{j_i(\alpha
)}) +\sum \mu_{i', j,k,l} x_jx_kx_le_{i'}=0$ is a nontrivial
relation. Let $i$ be minimal index such that either some
$\l_{i\alpha}\neq 0$ or some $\mu_{i,j,k,l}\neq 0$. Then $r\in
W_i$. Taking the relation $r$ modulo $W_{i+1}$ we have
$$ r\equiv (\sum
\l_{i\alpha} x^\alpha x_i +\sum \mu_{i, j,k,l} x_jx_kx_l)e_i\equiv
0 \mod W_{i+1},$$ (we used the fact that $j_i(\alpha) >i$), i.e.
$\sum \l_{i\alpha} x^\alpha x_i +\sum \mu_{i, j,k,l} x_jx_kx_l=0$
in $\L_{n,3}$, hence all $\l_{i\alpha}=0$ and all
$\mu_{i,j,k,l}=0$ since all the monomials are distinct, a
contradiction. This finishes the proof of the theorem.
 $\Box $

{\bf The dimension of $\S'\cap \S''$}. Recall that
$$\pi_n:=\begin{cases}
1& \text{if $n$ is odd},\\
2& \text{if $n$ is even}.\\
\end{cases}$$
 Let $\S / \S'':= \{ \s \S''\, | \, \s \in \S\} = \{ \s
\S''\, | \, \s \in \S'\}$ since $\S = \S'\S''$ (Corollary
\ref{a8Oct06}.(1)).

The next result shows that the subgroup $\S''$ of $\S$ is quite
large and that the intersection $\S'\cap \S''$ is a large subgroup
of $\S'$.

\begin{lemma}\label{i14Nov06}
Let $K$ be a commutative ring and $n\geq 4$. Then
\begin{enumerate} \item $\S'\cap \S''= \S'^5$ where $\S'^5:= \{ \s \in \S'\, |
\, (\s -1) (\gm )\subseteq \gm^5\}$, and so $\S'\cap \S''$ is a
closed subgroup of $\S$.\item The group $\S'\cap \S''$ is a
unipotent affine group over $K$ of dimension $$\dim (\S'\cap
\S'')= (n-2)2^{n-2}-n+\pi_n-(n-3){n\choose 2}$$ over $K$. \item
There is the  natural bijection $\S/\S'' \ra \S'/\S'^5\simeq
\prod_{i=1}^{n-1} \prod_{\alpha \in S_{i,1}'}\{ \rho_{i,j_i(\alpha
) ; \l x^\alpha }\}_{\l \in K},$ $ \s \S''\mapsto \s \S'$, where
$\s \in \S'$ (see Corollary \ref{a8Oct06}.(1)). The set
$\S'/\S'^5$ is an affine variety of dimension
$$ \dim (\S'/\S'^5)= n{n-1\choose 2}- {n\choose 2}= (n-3)
{n\choose 2}.$$
\end{enumerate}
\end{lemma}

{\it Proof}. 1. For $n=4$, the first statement is obvious as $
\S'^5= \S''^5= \{ e\}$ and $\S'\cap \S'' = \{ e\}$, by the very
definitions of the groups $\S'$ and $\S''$. We assume that $n\geq
5$.

 By (\ref{srp3}), any element $\s$ of $\S'$ is a
product $\s = \s_3\s_5\cdots \s_{2[\frac{n-1}{2}]+1}$ where $\s_i$
is a product of elements of the type $\rho_{i,j_{i,s}(\alpha_s );
\l x^{\alpha_s}}$, $\alpha_s\in S_{i,s}'$, $1\leq s\leq
[\frac{n-1}{2}]$, by (\ref{scooS}). Any element $\s$ of $\S'^5$ is
a product $\s = \s_5\s_7\cdots \s_{2[\frac{n-1}{2}]+1}$ (with
$\s_3=e$) where each $\s_i$ is product of elements of the type
$\rho_{i,j_{i,s}(\alpha_s ); \l x^{\alpha_s}}$, $\alpha_s\in
S_{i,s}'$, $2\leq s\leq [\frac{n-1}{2}]$. For $n\geq 6$, by
(\ref{com1}), if $s\geq 2$ then all $\rho_{i,j_{i,s}(\alpha_s );
\l x^{\alpha_s}}\in \S''$. Therefore, an element $\s =
\s_3\s_5\cdots \s_{2[\frac{n-1}{2}]+1}\in \S'$ belongs to the
group $\S''$ (i.e. $\s \in \S'\cap \S''$) iff $\s_3\in \S''$
(since $\s_5\cdots \s_{2[\frac{n-1}{2}]+1}\in \S''$)  iff
$\s_3=e$.

An element $\s = \s_3\s_5\cdots \s_{2[\frac{n-1}{2}]+1}\in \S'$
belongs to the group $\S'^5$  iff $\s_3\in \S'^5$ (since
$\s_5\cdots \s_{2[\frac{n-1}{2}]+1}\in \S'^5$)  iff $\s_3=e$.
Therefore, $\S'\cap \S''=\S'^5$, if $n\geq 6$.

If $n=5$ then $\S'^5= \{e\}$ by (\ref{Sem1}), and so $\S'\cap \S''
= \{ e \}$, by the very definitions of the groups $\S'$ and
$\S''$.

3.  By Corollary \ref{a8Oct06}.(1), $\S = \S'\S''$. Using
statement 1, we se that
$$ \S/\S''= \S'\S''/ \S''\simeq \S'/\S'\cap \S''\simeq
\S'/\S'^5\simeq \S'^3/\S'^5\simeq  \prod_{i=1}^{n-1} \prod_{\alpha
\in S_{i,1}'}\{ \rho_{i,j_i(\alpha ) ; \l x^\alpha }\}_{\l \in K}.
$$
So, $\S'/\S'^5$ is an affine variety.  Now, using the
 identifications as in the proof of Theorem \ref{g23Oct06} in the
 case $s=1$ there it follows at once that
 $$ \dim (\S'/\S'^5)= {\rm rk}_K(V)-{\rm rk}_K(\L_{n,2})=n{n-1\choose 2}-
{n\choose 2}= (n-3){n\choose 2}. $$ One can prove this fact
directly. Note that $|S'_{i,1}|= {n-1\choose 2}$, $1\leq i \leq
n-3$; $|S'_{n-2, 1}|= {n-1\choose 2}-1$ and $ |S'_{n-1, 1}|= {n-2
\choose 2}$. Then
$$ \dim (\S'/\S'^5) = \sum_{i=1}^{n-1}|S'_{i,1}|= (n-3) {
n-1\choose 2} +{n-1\choose 2} -1+{n-2\choose 2} =(n-3){n\choose
2}.$$
 2. By statement 1, $\S'\cap
\S'' = \S'^5$. Hence,
$$\dim (\S'\cap \S'' ) = \dim (\S'
)-\dim (\S'/\S'^5)=(n-2)2^{n-2}-n+\pi_n-(n-3){n\choose 2},$$ by
Corollary \ref{d25Oct06}. $\Box$

{\bf The dimension of $\S $}. The next theorem gives the dimension
 of the Jacobian group $\S$.

\begin{theorem}\label{19Nov06}
Let $K$ be a commutative ring and $n\geq 4$. The Jacobian group
$\S$ is a unipotent affine group over $K$ of dimension
$$\dim (\S )=\begin{cases}
(n-1)2^{n-1} -n^2+2& \text{if $n$ is even},\\
(n-1)2^{n-1} -n^2+1& \text{if $n$ is odd},\\
\end{cases}
 $$
 over $K$, i.e. the algebra of regular functions on $\S$ is a
 polynomial algebra in $\dim (\S )$ variables over $K$.
\end{theorem}

{\it Proof}. Recall that the algebraic group $\S'$ is affine and
$\dim (\S') = (n-2)2^{n-2}-n+\pi_n$ (Corollary \ref{d25Oct06});
$\S \simeq \S'\times \CF_n''$ (Corollary \ref{a8Oct06}.(5)); $\dim
( \CF_n'')= n(2^{n-2} -n+1)$, see (\ref{dimEns}). Therefore, the
algebraic group $\S$ is affine and
$$ \dim (\S ) = \dim (\S' )+\dim (\CF_n'')= (n-2)2^{n-2}
-n+\pi_n+n(2^{n-2} -n+1) = (n-1) 2^{n-1} -n^2+\pi_n.\Box $$

The Jacobian group $\S$ is a large subgroup  of $\G$ since
\begin{equation}\label{lmSG}
\lim_{n\rightarrow\infty} \frac{\dim (\S )}{\dim (\G
)}=\lim_{n\rightarrow\infty} \frac{(n-1) 2^{n-1}
-n^2+\pi_n}{n(2^{n-1} -n)}=1.
\end{equation}

{\bf The coordinates of $\S$}. The isomorphism (\ref{scooS}) and
the isomorphism in Theorem \ref{13Oct06}.(5) provide the explicit
coordinates for the Jacobian group $\S$ if $n\geq 4$.

{\bf The dimension of $\S''$}. The following theorem gives the
dimension of the group $\S''$ and proves that the group $\S''$ is
a {\em closed normal} subgroup of $\S$ (which is not obvious from
the outset).

\begin{theorem}\label{20Nov06}
Let $K$ be a commutative ring and $n\geq 4$. Then
\begin{enumerate}
\item $\S''= (\S'\cap \S'')\CF_n''= \S'^5\CF_n''$.\item $\S''$ is
 the  closed normal algebraic subgroup of $\S$. Moreover, $\S''$ is an affine group of dimension
$$\dim (\S'') = \begin{cases}
(n-1)2^{n-1}-n^2+2-(n-3){n\choose 2}& \text{if $n$ is even},\\
(n-1)2^{n-1}-n^2+1-(n-3){n\choose 2}& \text{if $n$ is odd},\\
\end{cases}
$$
and the factor group $\S/\S''\simeq \S'/\S'^5$ is an abelian
affine group of dimension $\dim (\S/\S'') = n{n-1\choose 2}-
{n\choose 2}= (n-3){n\choose 2}$. \item The map $ \S'\cap
\S''\times \CF_n''\ra \S''$, $(\s', \xi_{\od (n)}, \ldots , \xi_5,
\xi_3)\mapsto \s'\xi_{\od (n)} \cdots  \xi_5 \xi_3$, is an
isomorphism of algebraic varieties over $K$ with the inverse $\s
\mapsto \s'\xi_{\od (n)} \cdots  \xi_5 \xi_3$ given in Corollary
\ref{a8Oct06}.(5).
\end{enumerate}
\end{theorem}

{\it Proof}. $1$. The first equality follows from statement 3,
then the second equality follows from $\S'\cap \S''= \S'^5$ (Lemma
\ref{i14Nov06}.(1)).

 3. Since $ \CF''_n\subseteq \S''\subseteq \S$, statement 3
 follows from Corollary \ref{a8Oct06}.(5).

2. By Lemma \ref{i14Nov06}.(2) and statement 3, the group $\S''$
is affine and
\begin{eqnarray*}
\dim (\S'') &=& \dim (\S'\cap \S'') +\dim ( \CF_n'')\\
& =& (n-2) 2^{n-2} -n+\pi_n-(n-3) { n \choose 2} + n( 2^{n-2}
-n+1)\\
& =&  (n-1) 2^{n-1} - n^2+ \pi_n-(n-3) {n\choose 2},
\end{eqnarray*}
by Lemma \ref{i14Nov06}.(2) and (\ref{dimEns}).  Recall that $\S'$
is a closed subgroup of $\S$ and $\S'^5$ is a closed subgroup of
$\S'$, hence $\S''$ is a closed subgroup of $\S$ since
$$ \S'' \simeq (\S'\cap \S'')\times \CF_n''= \S'^5\times
\CF_n''\subseteq \S'\times \CF_n''\simeq \S.$$
Let us prove that
the group $\S''$ is a normal subgroup of the Jacobian group $\S$.
First, note that 
\begin{equation}\label{S5Ss}
\S^5\subseteq \S'',
\end{equation}
since $\S^5= \S'^5\S''^5$ (Corollary \ref{a8Oct06}.(2)) and
$\S'^5= \S'\cap \S''\subseteq \S''$ (Lemma \ref{i14Nov06}.(1)).
The subgroup $\S''$ is normal in $\S$ iff $\s \S''= \S''\s$ for
all $\s \in \S$. Note that
$$[\S , \S ]\subseteq [\G , \G ] \cap \S = [\G^3, \G^3] \cap \S \subseteq \G^5\cap \S=
 \S^5.$$ For any
$\tau \in \S''$,
$$ \s \tau = \s\tau \s^{-1} \tau^{-1} \tau \s = [\s , \tau ] \tau
\s \in [\S , \S ] \S''\s \subseteq \S^5\S'' \s =\S''\s , $$ by
(\ref{S5Ss}). This means that $ \s\S'' \subseteq  \S''\s$.
Similarly,
$$ \tau \s = \s \tau [ \tau^{-1} , \s^{-1}]\in \s \S''[\S, \S ]
\subseteq \s \S'' \S^5= \s \S'', $$
 by (\ref{S5Ss}). This means that $  \S''\s \subseteq \s \S''$.
 Therefore, $\s\S''= \S''\s$ for all $\s \in \S$, i.e. $\S''$ is
 a normal subgroup of $\S$. Finally, the factor group
 $$ \S / \S'' = \S'\S''/\S''\simeq \S'/\S'\cap \S'' = \S'/\S'^5$$
 is an abelian affine  group of dimension $(n-3){n\choose 2}$ (Lemma \ref{i14Nov06}.(3)),
 and its coordinates are given
 explicitly by Lemma \ref{i14Nov06}.(3). $\Box $

 {\bf The coordinates on $\S''$}. By Theorem \ref{20Nov06}.(3), each
automorphism $\s \in \S''$ is a unique product $\s = \s'\xi_{\od
(n)}\cdots  \xi_5 \xi_3$ where, by (\ref{scooS}), the $\s'$ is a
unique product 
\begin{equation}\label{1scooS}
 \s' = \prod_{s=2}^{[\frac{n-1}{2}]}\prod_{i=1}^{n-1}\prod_{\alpha
\in S_{i, s}'}\rho_{i,j_{i,s}(\alpha ); \l_\alpha x^\alpha }
\end{equation}
where $\alpha = \alpha_{i,s}$ (they depend on $i$ and $s$) and
$\l_{\alpha} = \l_{\alpha, i,s}\in K$. Therefore, $\{ \l_\alpha =
\l_{\alpha, i,s} \}$ and the coefficients of the elements that
define the automorphisms $\xi_i$ are  affine coordinates for the
algebraic group $\S''$ over the ring $K$. The group $\S''$
 is a large subgroup  of the Jacobian group $\S$ and the group
 $\S'$ is of `half size' of $\S''$ since
\begin{equation}\label{2scooS}
\lim_{n\rightarrow\infty} \frac{\dim (\S'' )}{\dim (\S )}=1,
\;\;\;\;\;  \lim_{n\rightarrow\infty} \frac{\dim (\S' )}{\dim
(\S'' )}=\frac{1}{2}.
\end{equation}

{\bf The dimension of $ \S'\backslash \S$}. For groups $A\subseteq
B$, let $A\diagdown B:= \{ Ab\, | \, b\in B\}$ and $B/A:= \{ bA\,
| \, b\in B\}$. If a group $A$ acts on sets $X$ and $Y$ then a map
$f: X\ra Y$ that respects the actions of the group $G$ on the sets
$X$ and $Y$ is called a $G$-{\em map}, i.e. $f(ax) = af(x)$ for
all $a\in A$ and $x\in X$. The isomorphism in Corollary
\ref{a8Oct06}.(5) is a $\S'$-isomorphism where the group $\S'$
acts by left multiplication on $\S$ and $\S'$ (in $\S'\times
\CF_n''$). Therefore, the set $\S' \diagdown \S$ is naturally
isomorphic to the set $\S' \diagdown\S'\times \CF_n''\simeq
\CF_n''$. The set $\CF_n''$ is an affine variety over $K$ of
dimension $n(2^{n-2} -n+1)$ (by (\ref{dimEns})), hence so is
$\S'\diagdown \S$. Note that
$$ \S'\diagdown\S = \S'\diagdown\S'\S''\simeq \S'\cap
\S''\diagdown\S'' \simeq \S'^5\diagdown\S''$$ since $\S= \S'\S''$
(Corollary \ref{a8Oct06}.(1)) and $\S'\cap \S'' = \S'^5$ (Lemma
\ref{i14Nov06}.(1)). So, we have proved the next corollary.

\begin{corollary}\label{s19Nov06}
Let $K$ be a commutative ring and $n\geq 4$. There are natural
isomorphisms of affine varieties over $K$: $ \CF_n''\simeq
\S'\diagdown\S \simeq \S'\cap \S''\diagdown\S'' \simeq
\S'^5\diagdown\S''$, each of them has dimension $n(2^{n-2}-n+1)$
over $K$.
\end{corollary}


\section{The image of the Jacobian map,
 the dimensions of the Jacobian ascents and of $\G /\S $}\label{IJMJA}

In this section, it is proved that all the Jacobian ascents
$\G_{2s}$ are distinct groups with a single exception (Corollary
\ref{a20Nov06}) and that their structure is completely determined
by the Jacobian group, $\G_{2s} = \G^{2s+1}\S$ (Theorem
\ref{G20Nov06}); each quotient space $\G_{2s}/\G_{2t}$ is an
affine variety (Corollary \ref{t27Nov06}) which is, via the
Jacobian map, canonically isomorphic to the affine variety
$E_{n,2s}'/E_{n,2t}'$ (Theorem \ref{x27Nov06}). In particular, the
quotient space $\G / \S$ is an affine variety of dimension
$2^{n-1}-\pi_n$ (Corollary \ref{y27Nov06}). The Jacobian map is a
surjective map if $n$ is odd and is not if $n$ is even. (Theorem
\ref{8Oct06}).

{\bf The equalities $\G_{2s} = \G^{2s+1} \S$}. It follows directly
from (\ref{char1}) and (\ref{JtriB}) that $\G^{2s+1}\S\subseteq
\G_{2s}$ for all $s=1,2, \ldots , [\frac{n-1}{2}]$. The next
theorem states that, in fact, the equalities hold, i.e.
$\G^{2s+1}\S= \G_{2s}$. The groups $\{ \G^{2s+1}\}$ have clear
structure and are given explicitly, therefore studying the
Jacobian ascents $\{ \G_{2s}\}$ is immediately reduced to studying
the Jacobian group $\S$.

\begin{theorem}\label{G20Nov06}
Let $K$ be a commutative ring and $n\geq 4$. Then
\begin{enumerate}
\item $\G_{2s}= \G^{2s+1}\S = \Phi^{2s+1}\S = \Phi'^{2s+1} \S$ for
 each   $s=1,2, \ldots , [\frac{n-1}{2}]$.  \item If $n$ is an even
number then $\G_n= \S $, i.e. $\G_n = \G_{n+2} =\S $.
\end{enumerate}
\end{theorem}

{\it Proof}. 1. The equalities $\G^{2s+1}\S = \Phi^{2s+1}\S $ are
obvious due to Theorem \ref{13Oct06}.(2). The equalities
$\Phi^{2s+1}\S = \Phi'^{2s+1} \S$ are obvious due to Theorem
\ref{g23Oct06}.(1).  The inclusions  $\G^{2s+1}\S\subseteq
\G_{2s}$ are obvious for $s=1,2, \ldots , [\frac{n-1}{2}]$. To
prove the reverse inclusions let $\s \in \G_{2s}$ (i.e. $\CJ (\s )
\in E_{2s}'$) where $1\leq s\leq  [\frac{n-1}{2}]$. We have to
show that $\s \in \G^{2s+1}\S$. If $\CJ (\s ) =1$, i.e. $\s \in
\S$, there is nothing to prove. So, we assume that $\CJ (\s ) \neq
1$, i.e. $\s \not\in \S$.  By Theorem \ref{13Oct06}.(1), $\s =
\phi \xi$ for some  elements $\phi\in \Phi^{2m+1}\backslash
\Phi^{2m+3}$ and $\xi \in \S''$. Now, we fix a presentation, say
$\s = \phi \xi$, with the least $m$. Using Lemma \ref{j23Oct06}
and Theorem \ref{g23Oct06}.(2,3), we see that (by the minimality
of $m$)
$$ \CJ (\s ) = \CJ (\phi \xi ) = \CJ (\phi ) \phi (\CJ (\xi ))= \CJ
(\phi ) \in E_{2m}'\backslash E_{2m+2}',$$ and so $\s \in
\G_{2m}\backslash \G_{2m+2}$, hence $s\leq m$ by the choice of $m$
and  since $\s \in \G_{2s}$. This proves that $\s \in \Phi^{2m+1}
\S \subseteq \G^{2m+1} \S \subseteq \G^{2s+1} \S$, as required.

2. Note that $\G_n \subseteq \G_{n-2} = \Phi'^{n-1}\S $ (by
statement 1) and $\Phi^{n+1} = \{ e \}$, and so $\Phi'^{n-1} /
\Phi^{n+1} = \Phi'^{n-1}$. If $\s \in \G_n$ then $\s = \phi \tau $
for some automorphisms $\phi \in \Phi'^{n-1}$ and $ \tau \in \S$,
and $E_{n,n}'= 1+ Kx_1x_2\cdots x_n \ni \CJ (\s ) = \CJ (\phi \tau
) = \CJ (\phi )$, hence $\phi \in \Phi^{n+1} = \{ e\}$, by Theorem
\ref{g23Oct06}.(3) and Lemma \ref{j23Oct06}. Therefore, $\s = \tau
\in \S$. This proves that $\G_n = \S$.  $\Box $

 Note that the number
$t:=2[\frac{n-1}{2}]+2$ is equal to $n+1$ if $n$ is odd; and to
$n$ if $n$ is even. Correspondingly, 
\begin{equation}\label{G2n2}
\G_{2[\frac{n-1}{2}]+2} = \begin{cases}
\G_{n+1} =  \S & \text{if $n$ is odd},\\
\G_{n} =  \S & \text{if $n$ is even (by Theorem \ref{G20Nov06}.(2))}.\\
\end{cases}
\end{equation}

Combining two results together, namely Theorem \ref{G20Nov06}.(1)
and Theorem \ref{g23Oct06}.(2), for each $s=1,2, \ldots ,
[\frac{n-1}{2}]$, there is a natural isomorphism of the abelian
groups: 
\begin{equation}\label{bG2s}
\G_{2s}/\G_{2s+2} \simeq \Phi'^{2s+1}/\Phi^{2s+3}.
\end{equation}
In more detail,
\begin{eqnarray*}
\G_{2s}/\G_{2s+2}  & = & \Phi^{2s+1} \S / \Phi^{2s+3} \S \simeq
\Phi^{2s+1} / \Phi^{2s+3}(\S \cap \Phi^{2s+1})= \Phi^{2s+1} /
\Phi^{2s+3} \S'^{2s+1}\\
&\simeq & (\Phi^{2s+1}/\Phi^{2s+3}) / (\Phi^{2s+3}
\S'^{2s+1}/\Phi^{2s+3}) \\
&\simeq &(\S'^{2s+1}/ \S'^{2s+3}\times \Phi'^{2s+1} /
\Phi^{2s+3})/ (\S'^{2s+1}/ \S'^{2s+3})\simeq
\Phi'^{2s+1}/\Phi^{2s+3}.
\end{eqnarray*}
By Lemma \ref{j23Oct06} and (\ref{bG2s}), for each  $s=1,2, \ldots
, [\frac{n-1}{2}]$, there is the natural isomorphism of the
abelian groups: 
\begin{equation}\label{bG3s}
 \G_{2s}/\G_{2s+2} \simeq \Phi'^{2s+1}/\Phi^{2s+3}\ra
E_{n,2s}'/E_{n,2s+2}', \;\; \s \G_{2s+2}\mapsto \CJ (\s )
E_{n,2s+2}'.
\end{equation}
This isomorphism and its inverse, (\ref{ibG3s}),  are some  of the
key results in finding the image of the Jacobian map (Theorem
\ref{8Oct06}). Recall that the map $ \L_{n,2s}\ra
E_{n,2s}'/E_{n,2s+2}'$, $a\mapsto (1+a)E_{n,2s+2}'$, is an
isomorphism of the abelian groups. By Theorem \ref{g23Oct06}.(3),
(see also Lemma \ref{j23Oct06}), the map 
\begin{equation}\label{ibG3s}
\L_{n,2s}\simeq E_{n,2s}'/E_{n,2s+2}'\ra \Phi'^{2s+1}/
\Phi^{2s+3}\simeq \G_{2s}/\G_{2s+2}, \;\; (1+a)E_{n,2s+2}'\mapsto
\phi'_a\Phi^{2s+3}(\simeq \phi'_a\G_{2s+2}),
\end{equation}
is the {\em inverse} map to the isomorphism (\ref{bG3s}) where
$\L_{n,2s}\ni a= a_{n-2s}+\cdots +a_n$ is the unique sum as in
Lemma \ref{s24Oct06} and the automorphism $\phi'_a$ is defined in
(\ref{srp1}), namely,  \begin{eqnarray*}
 \phi'_a(x_k)&:=&x_k, \;\; 1\leq k \leq n-2s-1, \\
 \phi'_a(x_i)&:=&x_i(1+a_i), \;\; n-2s\leq i \leq n.
\end{eqnarray*}
The fact that the map (\ref{ibG3s}) is the inverse of the map
(\ref{bG3s}) means that, for all $a\in \L_{n,2s}$,
\begin{equation}\label{jbG3s}
\CJ (\phi'_a) \equiv 1+a\mod E_{n,2s+2}',
\end{equation}
or, equivalently, for all $\s \in \G_{2s}$, 
\begin{equation}\label{kbG3s}
\s  \equiv \phi_a'\mod \G_{2s+2},
\end{equation}
where $\CJ (\s ) \equiv 1+a \mod E_{n,2s+2}'$ for a unique element
$a\in \L_{n,2s}$.

 Recall that the group $\G$ is equipped with the
Jacobian filtration $\{ \G_{2s}\}$,  and the group $E_n'$ is
equipped with the filtration $\{ E_{n, 2s}'\}$. Both filtrations
are descending. The Jacobian map $\CJ : \G \ra E_n'$, $\s \mapsto
\CJ (\s )$, is a filtered map, i.e. $\CJ (\G_{2s}) \subseteq
E_{n,2s}'$ for all $s$.

\begin{theorem}\label{27Nov06}
Let $K$ be a commutative ring, $n\geq 4$, and $s=1, \ldots ,
[\frac{n-1}{2}]$. Then each automorphism $\s \in \G$ is a unique
product $\s = \phi'_{a(2)}\phi'_{a(4)}\cdots
\phi'_{a(2[\frac{n-1}{2}])}\g$ for unique elements $a(2s)\in
\L_{n,2s}$ and $\g \in \G_{2[\frac{n-1}{2}]+2}=\S$ (by
(\ref{G2n2})). Moreover,
\begin{eqnarray*}
 a(2)&\equiv & \CJ (\s )-1 \mod E_{n,4}', \\
  a(2t)&\equiv & \CJ ({\phi'}_{a(2t-2)}^{-1} \cdots {\phi'}_{a(2)}^{-1}\s )-1 \mod E_{n,2t+2}', \; t=2, \ldots , [\frac{n-1}{2}], \\
  \g &= & (\phi'_{a(2)}\phi'_{a(4)}\cdots
\phi'_{a(2[\frac{n-1}{2}])})^{-1} \s.
\end{eqnarray*}
\end{theorem}

{\it Proof}. In brief, the theorem is a direct consequence of
repeated application of (\ref{kbG3s}). For $s=1$, by
(\ref{kbG3s}), $\s \equiv \phi_{a(2)}'\mod \G_4$ for a unique
element $a(2) \in \L_{n,2}$ such that $\CJ (\s ) \equiv 1+a(2)
\mod E_{n,4}'$. Now, $\s = \phi'_{a(2)}\s_4$ where $\s_4:=
\phi'^{-1}_{a(2)}\s \in \G_4$.  Repeating the same argument for
the automorphism $\s_4\in \G_4$ (i.e. for $s=2$), we have  $\s_4
\equiv \phi_{a(4)}'\mod \G_6$ for a unique element $a(4) \in
\L_{n,4}$ such that $\CJ (\s_4 ) \equiv 1+a(4) \mod E_{n,6}'$.
Then, $\s = \phi'_{a(2)}\phi'_{a(4)}\s_6$ where $\s_6:=
\phi'^{-1}_{a(4)}\phi'^{-1}_{a(2)}\s \in \G_6$. Continue in this
way we prove the theorem. $\Box$

We know already that the group $\G$ is an affine variety over $K$
where the coefficients $\{ \l_{\s , i, \alpha }\}$ of the
monomials $x^\alpha$ in the decomposition $\s (x_i) = x_i+
\sum_{|\alpha | \geq 2} \l_{\s , i, \alpha }x^\alpha$ (where $\s
\in \G$) are the coordinate functions on $\G$. Theorem
\ref{27Nov06} introduces the isomorphic  affine structure on $\G$
where the coefficients of the monomials $x^\alpha$ in $a( 2s)$ and
the coordinate functions on the Jacobian group $\S$ are new
coordinate functions on $\G$. We will see that this affine
structure on $\G$ is very useful in studying the spaces
$\G_{2s}/\G_{2t}$.

The next corollary is a direct consequence of Theorem
\ref{27Nov06}.

\begin{corollary}\label{s27Nov06}
Let $K$ be a commutative ring, $n\geq 4$, and $s=1, \ldots ,
[\frac{n-1}{2}]$. Then each automorphism $\s \in \G_{2s}$ is a
unique product $\s = \phi'_{a(2s)}\phi'_{a(2s+2)}\cdots
\phi'_{a(2[\frac{n-1}{2}])}\g$ for unique elements $a(2t)\in
\L_{n,2t}$, $s\leq t \leq [\frac{n-1}{2}]$,  and $\g \in
\G_{2[\frac{n-1}{2}]+2}=\S$ (by (\ref{G2n2})). Moreover,
\begin{eqnarray*}
 a(2s)&\equiv & \CJ (\s )-1 \mod E_{n,2s+2}', \\
  a(2t)&\equiv & \CJ ({\phi'}_{a(2t-2)}^{-1} \cdots {\phi'}_{a(2s)}^{-1}\s )-1 \mod E_{n,2t+2}', \; s<t \leq [\frac{n-1}{2}], \\
  \g &= & (\phi'_{a(2s)}\phi'_{a(4)}\cdots
\phi'_{a(2[\frac{n-1}{2}])})^{-1} \s.
\end{eqnarray*}
\end{corollary}

{\bf The dimension of $\G_{2s}/\G_{2t}$}. The next corollary shows
that the sets $\G_{2s}/\G_{2t}$ are affine varieties over $K$.

\begin{corollary}\label{t27Nov06}
Let $K$ be a commutative ring, $n\geq 4$, and $1\leq s<t \leq
[\frac{n-1}{2}]+1$. Then the set
$$\G_{2s}/\G_{2t}=\{ \phi'_{a(2s)}\cdots \phi'_{a(2t-2)}\G_{2t}\, | \, a(2s)\in \L_{n,2s}, \ldots , a(2t-2)\in
\L_{n,2t-2}\}$$ is an affine variety over $K$ of dimension $\dim
(\G_{2s}/\G_{2t})=\sum_{k=s}^{t-1} {n\choose 2k}$.
\end{corollary}

{\it Proof}. The first part of the corollary follows from
Corollary \ref{s27Nov06}, where the coefficients of the elements
$a(2s), \ldots , a(2t-2)$ are coordinate functions of the affine
variety $\G_{2s}/\G_{2t}$. Clearly, $\dim
(\G_{2s}/\G_{2t})=\sum_{k=s}^{t-1}{\rm rk}_K(\L_{n,2k})
=\sum_{k=s}^{t-1} {n\choose 2k}$. $\Box$

{\bf The dimension of the Jacobian ascents}. By Corollary
\ref{s27Nov06}, for each $n\geq 4$ and $s=1, \ldots
,[\frac{n-1}{2}]$, the Jacobian group 
\begin{equation}\label{pG2s}
\G_{2s} = \{  \phi'_{a(2s)}\cdots \phi'_{a(2[\frac{n-1}{2}])}\g \,
| \, a(2i) \in \L_{n,2i}, \g \in \S\}
\end{equation}
is an affine variety of dimension 
\begin{equation}\label{dimG2s}
\dim (\G_{2s})= \dim (\S ) + \sum_{i=s}^{[\frac{n-1}{2}]}{n\choose
2i}.
\end{equation}
 The coordinate functions for the affine group $\G_{2s}$ are
the coefficients of all the elements $a(2i)$ and the coordinate
functions on the Jacobian  group $\S$. In particular case when
$s=1$, one has 
\begin{equation}\label{pG22}
\G = \G_{2} = \{  \phi'_{a(2)}\cdots \phi'_{a(2[\frac{n-1}{2}])}\g
\, | \, a(2i) \in \L_{n,2i}, \g \in \S\}.
\end{equation}
It follows from (\ref{pG2s}) and (\ref{pG22}) that each Jacobian
ascent $\G_{2s}$, $s=1, \ldots ,[\frac{n-1}{2}]$, is a {\em
closed} subgroup of $\G$ that satisfies {\em exactly} $\dim (\G )
- \dim (\G_{2s})= \sum_{i=1}^{s-1}{n\choose 2i}$ defining
equations, namely, all coefficients of the elements $a(2), a(4),
\ldots , a(2s-2)$ are equal to zero.

Note that, for an even number $n$, $\G_{2[\frac{n}{2}]}=\G_n= \S $
(Theorem \ref{G20Nov06}.(2)). This means that, for each $n\geq 4$
(not necessarily even), the groups $\S$ and $\G_{2s}$, $s=1,
\ldots ,[\frac{n-1}{2}]$, are {\em all} the Jacobian ascents.

\begin{corollary}\label{a28Nov06}
Let $K$ be a commutative ring, $n\geq 4$. Then all the Jacobian
ascents are affine groups over $K$ and closed subgroups of $\G$,
and $\dim (\G_{2s})= \dim (\S ) +
\sum_{i=s}^{[\frac{n-1}{2}]}{n\choose 2i}$, $s=1, \ldots
,[\frac{n-1}{2}]$.
\end{corollary}

{\bf The isomorphisms $\bCJ_{s,t}$}.  For $1\leq s<t\leq
[\frac{n-1}{2}]+1$, the abelian group $E_{n,2s}'/E_{n,2t}'$ is an
affine variety over $K$ of the same dimension as the affine
variety $\G_{2s}/\G_{2t}$. The next result shows that the Jacobian
map 
\begin{equation}\label{Jbst}
\bCJ_{s,t}: \G_{2s}/\G_{2t}\ra E_{n,2s}'/E_{n,2t}', \;\; \s
\G_{2t}\mapsto \CJ (\s ) E_{n,2t}',
\end{equation}
is an isomorphism of affine varieties.

\begin{theorem}\label{x27Nov06}
Let  $K$ be a commutative ring, $n\geq 4$, and $1\leq s<t\leq
[\frac{n-1}{2}]+1$. Then the Jacobian map $\bCJ_{s,t}$,
(\ref{Jbst}), is an isomorphism of affine varieties.
\end{theorem}

{\it Proof}. By the definition, the map $\bCJ_{s,t}$ is a
polynomial map. In order to finish the proof of the theorem it
suffices to show that the map $ \bCJ_{s,t}$ is a bijection and its
inverse is also a polynomial map. For a given $t$, to prove these
two statements, we will use downward induction on  $s$ starting at
$s=t-1$ where the result is known due to (\ref{bG3s}),
(\ref{ibG3s}), (\ref{jbG3s}), and (\ref{kbG3s}). If $t=2$ then
$s=1$, and we are done. So, let $t\geq 3$ and $s<t-1$, and, by the
inductive hypothesis, we assume that the map $\bCJ_{s+1, t}$ is an
isomorphism of affine varieties. We are going to present the
inverse map for $\bCJ_{s,t}$ which is, by construction,  a
polynomial map.

By Corollary \ref{t27Nov06}, each element of $\G_{2s}/\G_{2t}$ can
be written uniquely in the form $\phi_a'\tau \G_{2t}$ where $a\in
\L_{n,2s}$ and $\tau \G_{2t}\in \G_{2s+2} / \G_{2t}$. Similarly,
each element of $E_{n,2s}'/E_{n,2t}'$ can be written uniquely in
the form $(1+a) bE_{n,2t}'$ where $a\in \L_{n,2s}$ and $
bE_{n,2t}'\in  E_{n, 2s+2}'/E_{n,2t}'$. To finish the proof we
have to show that, for a given element $(1+a) bE_{n,2t}'\in E_{n,
2s}'/E_{n,2t}'$, and an unknown $\phi'_{a'}\tau \G_{2t}\in
\G_{2s}/\G_{2t}$, the equation
$$ \bCJ_{s,t} (\phi'_{a'}\tau \G_{2t})=(1+a) bE_{n,2t}'$$
has a unique solution $\phi'_{a'}\tau \G_{2t}$ that depends
polynomially on the RHS. By taking the equation modulo
$E_{n,2s+2}'$, we obtain the equality $a=a'$, by (\ref{jbG3s}):
$$1+a\equiv \CJ ( \phi'_{a'})\equiv 1+a'\mod E_{n,2s+2}'.$$
Now, we can solve the equation explicitly which can be written as
follows
$$ \CJ ( \phi'_a) \phi'_a(\bCJ_{s+1, t}(\tau \G_{2t}))= \bCJ_{s,t}
(\phi'_a\tau \G_{2t}) = (1+a) b E_{n, 2t}'.$$ Namely,
\begin{equation}\label{st2t}
\tau \G_{2t} = (\bCJ_{s+1, t})^{-1} \phi'^{-1}_a(\CJ
(\phi'_a)^{-1} (1+a) bE_{n,2t}')
\end{equation}
is the unique solution that depends polynomially on the RHS (note
that $\CJ (\phi'_a)^{-1} (1+a)\in E_{n,2s+2}'$, by (\ref{jbG3s})),
as required.  $\Box $

{\bf The Jacobian ascents are distinct groups except one case}.
Now, we are ready to give an answer to the question of whether the
Jacobian ascents are distinct groups or not.

\begin{corollary}\label{a20Nov06}
Let $K$ be a commutative ring and $n\geq 4$.
\begin{enumerate}
\item If $n$ is an odd number then the Jacobian ascents $$ \G =
\G_2\supset \G_4\supset \cdots \supset \G_{2s} \supset \cdots
\supset \G_{2[\frac{n}{2}]}\supset \G_{2[\frac{n}{2}]+2}=\S$$ are
distinct groups. \item If $n$ is an even  number then the Jacobian
ascents
$$ \G = \G_2\supset \G_4\supset \cdots \supset \G_{2s} \supset
\cdots \supset\G_{2[\frac{n}{2}]-2}\supset\G_{2[\frac{n}{2}]}
=\G_{2[\frac{n}{2}]+2}=\S$$ are distinct groups except the last
two groups, i.e. $\G_{2[\frac{n}{2}]} =\G_{2[\frac{n}{2}]+2}$.
\end{enumerate}
\end{corollary}

{\it Proof}. 1. If $n$ is odd then $[\frac{n-1}{2}]=[\frac{n}{2}]$
and the result follows from Theorem \ref{x27Nov06} since the
groups $\{ E_{n,2s}'\}$ are distinct for $s=1,2,\ldots ,
[\frac{n-1}{2}]$.

2.  If $n$ is even then $[\frac{n-1}{2}]=[\frac{n}{2}]-1$ and
$2[\frac{n}{2}]=n$. By Theorem \ref{x27Nov06}, the following
groups are distinct: $\G = \G_2\supset \G_4\supset \cdots \supset
\G_{2s} \supset \cdots \supset\G_{2[\frac{n}{2}]-2} \supset \S$.
By Theorem \ref{G20Nov06}.(2), $\G_{2[\frac{n}{2}]}
=\G_{2[\frac{n}{2}]+2}=\S$, and so the result. $\Box $

{\bf The dimension of $\G / \S$}. By taking the extreme values for
$s$ and $t$ in Corollary  \ref{t27Nov06}, namely, $s=1$ and
$t=[\frac{n-1}{2}]+1$,  we see that $\G / \S$ is an affine variety
due to (\ref{G2n2}) and $\G = \G_2$.  The next corollary gives the
dimension of the variety  $\G/ \S$.

\begin{corollary}\label{y27Nov06}
Let $K$ be a commutative ring and $n\geq 4$. Then $\G / \S$ is an
affine variety.
\begin{enumerate}
\item If $n$ is odd then the Jacobian map $\G / \S \ra E_n'$,
$\s\S\mapsto \CJ (\s )$, is an isomorphism of the affine varieties
over $K$, and $\dim (\G / \S ) = 2^{n-1} -1$. \item If $n$ is even
then the Jacobian map $\G / \S \ra E_n'/E_{n,n}'$, $\s\S\mapsto
\CJ (\s )E_{n,n}'$, is an isomorphism of the affine varieties over
$K$ (where $E_{n,n}'= 1+Kx_1\cdots x_n$), and $\dim (\G / \S ) =
2^{n-1} -2$.
\end{enumerate}
\end{corollary}

{\it Proof}. 1. Take $s=1$ and $t= [\frac{n-1}{2}]+1$ in Theorem
\ref{x27Nov06}. Since $n$ is an odd number, $2t+2= n+1$, and so
$E_{n,2t+2}'= E_{n,n+1}'=\{1\}$ and $\G_{2t+2}= \G_{n+1} =\S$
(Corollary \ref{a20Nov06}.(1)). By Theorem \ref{x27Nov06}, the
Jacobian map $\G / \S \ra E_n'$, $\s\S\mapsto \CJ (\s )$, is an
isomorphism of the affine varieties over $K$. Now,
$$ \dim (\G / \S ) = \dim (E_n') =
\sum_{s=1}^{[\frac{n}{2}]}{n\choose 2s}=
\sum_{s=0}^{[\frac{n}{2}]}{n\choose 2s}-1=2^{n-1} -1.$$

2. Similarly, take $s=1$ and $t= [\frac{n-1}{2}]+1$ in Theorem
\ref{x27Nov06}. Since $n$ is an even number, $2t+2= n$, and so
$E_{n,2t+2}'= E_{n,n}'=1+Kx_1\cdots x_n$ and $\G_{2t+2}= \G_{n}=
\S$ (Corollary \ref{a20Nov06}.(2)). By Theorem \ref{x27Nov06}, the
Jacobian map $\G / \S \ra E_n'/E_{n,n}'$, $\s\S\mapsto \CJ (\s
)E_{n,n}'$, is an isomorphism of the affine varieties over $K$,
and
$$ \dim (\G / \S ) = \dim (E_n'/E_{n,n}') =
\sum_{s=1}^{[\frac{n}{2}]}{n\choose 2s}-1=
\sum_{s=0}^{[\frac{n}{2}]}{n\choose 2s}-2=2^{n-1} -2.\;\; \Box $$

Theorem \ref{8Oct06} gives an answer to the natural question of
whether the Jacobian map $\CJ : \G \ra E_n'$, $\s \mapsto \CJ (\s
)$, is surjective? The answer is `yes' for odd numbers $n$, and,
surprisingly,  `no' for even numbers $n$. In the second case, the
image of the Jacobian map is large. More precisely, it is a {\em
closed} subvariety of $E_n'$ of codimension 1 which is defined by
a single equation. Moreover, it is canonically isomorphic to the
affine variety $E_n'/E_{n,n}'$.

\begin{theorem}\label{8Oct06}
Let $K$ be a commutative ring, $n\geq 4$,  $\CJ :\G \ra E_n'$, $\s
\mapsto \CJ (\s )$, be the Jacobian map, and $s=1,2, \ldots ,
[\frac{n-1}{2}]$. Then,
\begin{enumerate}
\item for an odd number $n$, the Jacobian map $\CJ$ is surjective.
Moreover, for each $s$, the map $ \CJ : \G_{2s}\ra E_{n,2s}'$, $\s
\ra \CJ (\s )$, is surjective; and
 \item for an even  number $n$, the Jacobian map $\CJ$ is not surjective but very close to be a surjective
 map. In more detail,
\begin{enumerate}
\item the image ${\rm im } (\CJ )$ is a closed algebraic variety
of $E_n'$ of codimension 1 (i.e. $\dim ({\rm im } (\CJ
))=2^{n-1}-2$) which  is defined by a single equation (see the
proof), \item ${\rm im } (\CJ )\cap E_{n,n}'=\{ 1\}$ where
$E_{n,n}'= 1+Kx_1\cdots x_n, $\item the image ${\rm im }(\CJ )$ is
canonically isomorphic to the algebraic group $E_n'/E_{n,n}'$ via
the map ${\rm im}(\CJ ) \ra E_n'/E_{n,n}'$, $\alpha  \mapsto
\alpha   E_{n,n}'$.
\end{enumerate}
\end{enumerate}
\end{theorem}

{\it Proof}. 1. The fact that the Jacobian map $\CJ : \G \ra
E_n'$, $ \s \mapsto \CJ (\s )$, is surjective follows from
Corollary \ref{y27Nov06}.(1) and (\ref{GSim}).

 2. By Corollary \ref{y27Nov06}.(2) and (\ref{GSim}), the Jacobian
 map $\CJ$ is not surjective, though there is a bijection between
 the image ${\rm im} (\CJ )$ and $E_n'/E_{n,n}'$. The set
 $E_n'/E_{n,n}'$ may be identified with the closed affine
 subvariety of the affine variety $E_n'$ that is given by a single
 equation: the coefficient of the element $x_1\cdots x_n$ is
 equal to zero. In more detail, $E_n'=
 1+\oplus_{s=1}^{[\frac{n}{2}]}\L_{n,2s}$ and $E_n'/E_{n,n}'$ is
 identified with $ 1+\oplus_{s=1}^{[\frac{n}{2}]-1}\L_{n,2s}$. Then
 the bijection between ${\rm im} (\CJ )$ and $E_n'/E_{n,n}'$ means
 that the last coordinate, say $\l = \l (\s )$, of each element
 $\CJ (\s )= 1+\cdots + \l x_1\cdots x_n$ is a polynomial function
 of the previous coordinates (in the three dots expression). This
 is the defining equation of the image ${\rm im } (\CJ )$ in
 $E_n'$.  So, the statements (a) and (c) follow. The statement (b)
 follows at once from the equality $\G_n= \S $ (Theorem
 \ref{G20Nov06}.(2)): $\s \in {\rm im} (\CJ )\cap E_{n,n}'$
  iff $\s \in \G_n= \S$.  $\Box$


\section{Analogues of the Poincar\'{e} Lemma}\label{APL}
In this section, two results (Theorems \ref{s14Sep06} and
\ref{P14Sep06}) are proved  that have flavour of the Poincar\'{e}
Lemma. Theorem \ref{s14Sep06} is used in the proof of Theorem
\ref{M30Sep06}.

\begin{theorem}\label{14Sep06}
Let $K$ be an arbitrary (not necessarily commutative)  ring. Then
\begin{enumerate}
\item the Grassmann ring $\L_n(K)$ is a direct sum of right
$K$-modules
\begin{eqnarray*}
\L_n(K)&=& x_1\cdots x_nK \oplus x_1\cdots x_{n-1}K \oplus
x_1\cdots x_{n-2}K\lfloor x_n\rfloor \oplus\cdots \\
&\cdots &\oplus x_1\cdots x_iK\lfloor x_{i+2}\ldots ,  x_n\rfloor
\oplus\cdots \oplus  x_1 K\lfloor x_3\ldots ,  x_n\rfloor\oplus
K\lfloor x_2\ldots , x_n\rfloor .
\end{eqnarray*}
\item So, each element $a\in \L_n(K)$ is a unique sum
$$ a= x_1\cdots x_na_n+ x_1\cdots x_{n-1}b_n+\sum_{i=1}^{n-2}
x_1\cdots x_ib_{i+1} + b_1$$ where $a_n, b_n\in K$, $b_i\in
 K\lfloor x_{i+1}\ldots ,  x_n\rfloor$, $1\leq i\leq n-1$.
 Moreover,
 \begin{eqnarray*}
 a_n&=& \der_n\der_{n-1}\cdots \der_1(a), \\
b_{i+1}&=&\der_i\der_{i-1}\cdots \der_1(1-x_{i+1}\der_{i+1})(a), \; 1\leq i\leq n-1, \\
 b_1&=&(1-x_1\der_1)(a).
\end{eqnarray*}
So, $$ a= x_1\cdots x_n\der_n\der_{n-1}\cdots
\der_1(a)+\sum_{i=1}^{n-1} x_1\cdots x_i\der_i\cdots
\der_1(1-x_{i+1}\der_{i+1})(a)+(1-x_1\der_1)(a).$$
\end{enumerate}
\end{theorem}

 {\it Proof}. For each $i=1, \ldots , n$, let $K_i:=
K\lfloor x_i, \ldots , x_n\rfloor$ and $K_{n+1}:= K$.

 1. Existence of the decomposition is a consequence of a repeated
 use of the fact that $K_i= x_iK_{i+1}\oplus K_{i+1}$. Namely,
 \begin{eqnarray*}
 K_n&=& x_1K_2\oplus K_2= x_1(x_2K_3\oplus K_3)\oplus K_2= x_1x_2K_3\oplus x_1K_3\oplus K_2 \\
 &=&x_1x_2(x_3K_4\oplus K_4)\oplus x_1K_3\oplus
 K_2\\
 &=& x_1x_2x_3(x_4K_5\oplus K_5)\oplus x_1x_2K_4\oplus x_1K_3\oplus
 K_2=\cdots
\end{eqnarray*}
when this process stops after $n$ steps we get the required
decomposition.

$2$. The crucial steps in finding the coefficients for the element
$a$ are $(i)$ $\der_i^2=\cdots =\der_n^2=0$, and $(ii)$ for each
$i=1, \ldots , n$, the map $\phi_i:= 1-x_i\der_i:\L_n\ra \L_n$ is
the projection onto the Grassmann subring $ K\lfloor x_1, \ldots ,
\widehat{x_i}, \ldots , x_n \rfloor$ in the decomposition $\L_n
=K\lfloor x_1, \ldots , \widehat{x_i}, \ldots , x_n \rfloor\oplus
x_i K\lfloor x_1, \ldots , \widehat{x_i}, \ldots , x_n \rfloor$
(Lemma \ref{p9Sep06}). The tail $t$ in the sum $a= x_1\cdots x_n
a_n +t$ has (total) degree in the variables $x_1, \ldots , x_n$
strictly less than $n$, hence $t$ is killed by the map
$\der_n\cdots \der_1$. Therefore,
$$\der_n\cdots \der_1(a)=\der_n\cdots \der_1(x_1\cdots x_na_n)=\der_n\cdots \der_2(x_2\cdots
x_na_n)=\cdots =a_n.$$ To find the elements $b_i$ we use induction
on $i$. Since the map $\phi_1= (1-x_1\der_1): \L_n\ra \L_n$ is a
projection onto $K\lfloor x_2, \ldots , x_n \rfloor$ and all
summands of $a$ but the last belong to the ideal $(x_1)$ (which is
annihilated by $\phi_1$), it follows at once that $\phi_1(a)=
\phi_1(b_1)=b_1$. Similarly, applying $\phi_2$ to $a$ we see that
$\phi_2(a) = x_1b_2+\phi_2(b_1)$. Since $\phi_2( b_1) \in K\lfloor
x_3, \ldots ,  , x_n \rfloor$, we have $\der_1\phi_2(b_1) =0$, and
so $\der_1\phi_2(a) =\der_1(x_1b_2) = b_2$. Suppose that the
formula for the $b_k$ in the theorem is true for all $k=1, \ldots
, i$, we have to prove it for $i+1$. The cases $i=1,2$ have been
established already. So, let $i\geq 3$. Now,
$$ \phi_{i+1} (a) = x_1\cdots x_i b_{i+1} +\phi_{i+1}
(\sum_{k=1}^{i-1} x_1\cdots x_k\der_k\cdots \der_1\phi_{k+1}
(a))+\phi_{i+1} (b_1).$$ Note that the skew derivations $\der_1,
\ldots , \der_i$ {\em commute} with $\phi_{i+1}$. $\phi_{i+1}
(b_1) \in K\lfloor x_2, \ldots , \widehat{x_{i+1}}, \ldots , x_n
\rfloor$ implies $\der_1\phi_{i+1} (b_1) =0$, and so $\der_i\cdots
\der_1\phi_{i+1} (b_1)=0$. For each $k=1, \ldots  , i-1$, let
$c_k= \phi_{i+1} (x_1\cdots x_k\der_k\cdots \der_1\phi_{k+1}(a))$.
Using the commutation relations for the Grassmann $K$-algebra
$\L_k=\oplus_{\alpha , \beta \in \CB_k} \der^\alpha x^\beta K$ one
can write (in $\L_k$)
$$ x_1\cdots x_k\der_k\cdots \der_1 = 1+d_k \;\; {\rm where}\;\;
d_k\in \oplus_{0\neq \alpha \in \CB_k , \beta \in \CB_k}
\der^\alpha x^\beta K.$$ Since $\der_k\cdots \der_1 d_k=0$ (as
$\der_1^2=\cdots = \der_k^2=0$) and $i>k$, we have
\begin{eqnarray*}
 \der_i\cdots \der_1c_k&=& \phi_{i+1}\der_i\cdots \der_1(1+d_k) \phi_{k+1} (a)=
 \phi_{i+1} \der_i\cdots \der_{k+1} \cdots \der_1\phi_{k+1} (a)  \\
 &=& (-1)^k \phi_{i+1} \der_i\cdots \widehat{\der_{k+1}} \cdots \der_1\der_{k+1}\phi_{k+1}
 (a)=0
\end{eqnarray*}
since $\der_{k+1}\phi_{k+1}=0$. Now, we see that
\begin{eqnarray*}
 \der_i\cdots \der_1\phi_{i+1}(a)&=& \der_i\cdots \der_1
 (x_1\cdots x_ib_{i+1})\\
 &=& \der_i\cdots \der_1
 (x_1\cdots x_i)b_{i+1}\;\;\; ({\rm as}\;\; b_{i+1} \in K\lfloor x_{i+2}, \ldots ,
 x_n\rfloor\subseteq \cap_{k=1}^i\ker (\der_k))\\
 &=& b_{i+1},
\end{eqnarray*}
as required. $\Box $

By Theorem \ref{14Sep06}, the identity map ${\rm id}_{\L_n}
:\L_n\ra \L_n$ is equal to 
\begin{equation}\label{idLn}
{\rm id}_{\L_n}=x_1\cdots x_n\der_n\der_{n-1}\cdots \der_1+
\sum_{i=1}^{n-1} x_1\cdots x_i\der_i\cdots
\der_1(1-x_{i+1}\der_{i+1})+(1-x_1\der_1).
\end{equation}
If $n'\geq n$ then the RHS of (\ref{idLn}) is a map from $\L_{n'}$
to itself. Therefore, 
\begin{equation}\label{1idLn}
{\rm id}_{\L_n'}=x_1\cdots x_n\der_n\der_{n-1}\cdots \der_1+
\sum_{i=1}^{n-1} x_1\cdots x_i\der_i\cdots
\der_1(1-x_{i+1}\der_{i+1})+(1-x_1\der_1).
\end{equation}

\begin{theorem}\label{s14Sep06}
Let $K$ be an arbitrary ring, $u_1, \ldots , u_n\in \L_n(K)$, and
$ a\in \L_n(K)$ be an unknown. Then the system of equations
$$\begin{cases}
x_1a=u_1 \\
x_2a=u_2 \\
\;\;\;\; \;\;\;\vdots \\
x_na=u_n
\end{cases}
$$
has a solution in $\L_n$ iff the following two conditions hold
\begin{enumerate}
\item $u_1\in (x_1), \ldots , u_n\in (x_n)$, and \item
$x_iu_j=-x_ju_i$ for all $i\neq j$.
\end{enumerate}
In this case, 
\begin{equation}\label{alsol}
a= x_1\cdots x_na_n+\sum_{i=1}^{n-1} x_1\cdots x_i\der_i\cdots
\der_1\der_{i+1}(u_{i+1})+\der_1(u_1), \;\; a_n \in K,
\end{equation}
are all the solutions.
\end{theorem}

{\it Remark}. An analogue of the Poincar\'{e} Lemma for $\L_n$ is
given later (Theorem \ref{P14Sep06}). Theorem \ref{s14Sep06} is a
sort of Poincar\'{e} Lemma for the Grassmann algebra since the map
$l_{x_i}:\L_n\ra \L_n$, $ u \mapsto x_iu$,  the left
multiplication by $x_i$,  is a sort of skew partial derivatives on
$x_i$ as follows from the following two properties:

1. Each element $a\in \L_n$ is a unique sum $ a= x_i\alpha +\beta$
with $\alpha , \beta \in K\lfloor x_1, \ldots , \widehat{x_i},
\ldots x_n\rfloor$ and $l_{x_i}(a) = x_i \beta$; and

$2$. for any two elements $a_s\in \L_{n, s}$ and $a_t\in \L_{n,
t}$ where $s,t \in \Z_2$:
$$ x_i( a_sa_t) = \frac{1}{2}x_ia_sa_t+ \frac{1}{2}x_ia_sa_t=
(\frac{1}{2}x_ia_s)a_t+ (-1)^sa_s(\frac{1}{2}x_ia_t),$$
 provided $\frac{1}{2}\in K$.

{\it Proof}. Suppose that $a\in \L_n$ is a solution then $u_i=
x_ia\in (x_i)$ for all $i$; and, for all $i\neq j$,
$$ x_iu_j+x_ju_i= x_ix_ja+ x_jx_ia= x_ix_ja -x_ix_ja=0.$$
So, conditions 1 and 2 hold. Evaluating the skew derivation
$\der_i$ at the equality $u_i= x_ia$ one sees that
\begin{equation}\label{aid1i}
\der_i(u_i) = \der_i(x_ia) = (1-x_i\der_i) (a).
\end{equation}
Let us write the element $a$ as the sum in Theorem \ref{14Sep06}.
Note that if $a$ is a solution to the system then $a+x_1\cdots
x_na_n$ is also a solution for an arbitrary choice of $a_n \in K$,
and vice versa. By (\ref{aid1i}) and Theorem \ref{14Sep06}.(2),
$$ a= x_1\cdots x_n a_n+\sum_{i=1}^{n-1} x_1\cdots x_i
\der_i\cdots \der_1\der_{i+1} (u_{i+1})+\der_1(u_1).$$ This proves
(\ref{alsol}).

It remains to show that if conditions 1 and 2 hold then
(\ref{alsol}) are solutions to the system. We prove directly that
$x_ja=u_j$ for all $j$. An idea of the proof is to use the
identity (\ref{1idLn}). For $j=1$, note that $x_1\der_1(u_1) =
u_1$ since $u_1\in (x_1)$, and so $x_1a= x_1\der_1 (u_1) = u_1$.
Suppose that $2\leq j \leq n$. Then
$$ x_ja= x_1\cdots x_{j-1} \der_{j-1} \cdots \der_1 x_j\der_j(u_j)
+\sum_{i=1}^{j-2}x_1\cdots x_i \der_i \cdots \der_1
x_j\der_{i+1}(u_{i+1})+x_j\der_1(u_1).$$ Note that
$x_j\der_{i+1}(u_{i+1}) = -\der_{i+1}(x_ju_{i+1})=
-\der_{i+1}(-x_{i+1} u_j) = (1-x_{i+1}\der_{i+1}) (u_j)$;
$x_j\der_j(u_j) = u_j$ since $u_j\in (x_j)$; and $x_j\der_1(u_1) =
-\der_1(x_ju_1) = -\der_1(-x_1u_j)= (1-x_1\der_1)(u_j)$. Using
these equalities, we see that
$$ x_ja=(x_1\cdots x_{j-1}\der_{j-1} \cdots
\der_1+\sum_{i=1}^{j-2}x_1\cdots x_i\der_i\cdots
\der_1(1-x_{i+1}\der_{i+1})+(1-x_1\der_1))(u_j)=u_j,$$ by
(\ref{1idLn}). $\Box$

\begin{theorem}\label{P14Sep06}
Let $K$ be an arbitrary ring, $u_1, \ldots , u_n\in \L_n(K)$, and
$ a\in \L_n(K)$ be an unknown. Then the system of equations
$$\begin{cases}
\der_1(a)=u_1 \\
\der_2(a)=u_2 \\
\;\;\;\; \;\;\;\vdots \\
\der_n(a)=u_n
\end{cases}
$$
has a solution in $\L_n$ iff the following two conditions hold
\begin{enumerate}
\item for each $i=1, \ldots , n$, $u_i\in K\lfloor x_1, \ldots
,\widehat{x_i} ,\ldots , x_n \rfloor$, and \item
$\der_i(u_j)=-\der_j(u_i)$ for all $i\neq j$.
\end{enumerate}
In this case, 
\begin{equation}\label{aPlsol}
a=\l +\sum_{0\neq \alpha \in \CB_n} \phi (u_\alpha ) x^\alpha,
\;\; \l \in K,
\end{equation}
are all the solutions where $\phi$ is defined in Lemma
 \ref{p9Sep06}.(3) and, for $\alpha = \{ i_1<\cdots < i_k\}$,
$u_\alpha := \der_{i_k}\der_{i_{k-1}}\cdots \der_{i_2}(u_{i_1})$.
\end{theorem}

{\it Proof}. Suppose that $a\in \L_n$ is a solution then $u_i=
\der_i(a)\in {\rm im} (\der_i) = K\lfloor x_1, \ldots
,\widehat{x_i} ,\ldots , x_n \rfloor$, and so the first condition
holds. For all $i\neq j$,
$$ \der_i(u_j)= \der_i\der_j (a) = -\der_j\der_i(a) =
-\der_j(u_i), $$ and so the second condition holds. Note that if
$a$ is a solution then $a+\l$, $\l\in K$, are all the solutions
since $K= \cap_{i=1}^n \ker (\der_i)$. By Theorem
\ref{s9Sep06}.(1),
$$a= \sum_{\alpha \in \CB_n} \phi (\der^\alpha (a)) x^\alpha = \l
+\sum_{0\neq \alpha \in \CB_n} \phi (\der^\alpha (a)) x^\alpha =
\l +\sum_{0\neq \alpha \in \CB_n} \phi (u_\alpha ) x^\alpha,$$ so
(\ref{aPlsol}) holds.

It remains to show that if conditions 1 and 2 hold then
(\ref{aPlsol}) are solutions to the system. We prove directly that
$\der_i(a) = u_i$ for all $i$. An idea of the proof is to use the
equality of Theorem \ref{s9Sep06}.(1) together with conditions 1
and 2.
$$ \der_i(a) = \sum_{i\in \alpha \in \CB_n} \phi (u_\alpha )
(-1)^{\alpha_1+\cdots + \alpha_{i-1}}x^{\alpha\backslash \{ i \}}=
\sum_{i\in \alpha \in \CB_n} \phi (\der^{\alpha\backslash \{ i
\}}(u_i))x^{\alpha\backslash \{ i \}}=u_i.$$ The second equality
above is due to the fact that $(-1)^{\alpha_1+\cdots +
\alpha_{i-1}}u_\alpha =\der^{\alpha\backslash \{ i \}}(u_i)$, by
condition 2. The last equality follows from Theorem
\ref{s9Sep06}.(1) and condition 1. $\Box$

\section{The unique  presentation $\s = \o_{1+a} \g_b\s_A $ for
$\s\in \Aut_K(\L_n)$ }\label{PREXS}

In this section, $K$ is a {\em reduced commutative} ring with
$\frac{1}{2}\in K$. By Theorem \ref{29Sep06}.(3), $G= \O \G
\GL_n(K)^{op}$. So, each element $\s \in G$ has the unique
presentation as the product $\s = \o_{1+a} \g_b\s_A $ where
$\o_{1+a} \in \O$ ($a\in \L_n'^{od}$), $\g_b \in \G$, $\s_A\in
\GL_n(K)^{op}$ where $\L_n'^{od} := \oplus_i \L_{n,i}$ and  $i$
runs through all {\em odd} natural numbers such that $1\leq i\leq
n-1$.

\begin{theorem}\label{M30Sep06}
Let  $K$ be a  reduced commutative ring with $\frac{1}{2}\in K$.
Then each element $\s \in G$ is a unique product $\s = \o_{1+a}
\g_b\s_A$ (Theorem \ref{29Sep06}.(3)) where $a\in \L_n'^{od}$ and
\begin{enumerate}
\item $\s (x) = Ax +\cdots $ (i.e. $\s (x) \equiv Ax\mod \gm$) for
some $A\in \GL_n(K)$, \item $b= A^{-1} \s(x)^{od}-x$, and \item
$a= -\frac{1}{2}\g_b (\sum_{i=1}^{n-1} x_1\cdots x_i\der_i\cdots
\der_1\der_{i+1} (a_{i+1}') +\der_1(a_1'))$ where $a_i':= (A^{-1}
\g_b^{-1} (\s (x)^{ev}))_i$, the $i$'th component of the
column-vector $A^{-1} \g_b^{-1} (\s (x)^{ev})$.
\end{enumerate}
\end{theorem}

{\it Remark}.  Recall that  $x=\begin{pmatrix}
  x_1\\
 \vdots  \\
x_n \\
\end{pmatrix}$, $b=\begin{pmatrix}
  b_1\\
 \vdots  \\
b_n \\
\end{pmatrix}$, $\s (x)=\begin{pmatrix}
  \s (x_1)\\
 \vdots  \\
\s (x_n) \\
\end{pmatrix}$, $\s (x)^{ev}=\begin{pmatrix}
  \s (x_1)^{ev}\\
 \vdots  \\
\s (x_n)^{ev} \\
\end{pmatrix}$, $\s (x)^{od}=\begin{pmatrix}
  \s (x_1)^{od}\\
 \vdots  \\
\s (x_n)^{od} \\
\end{pmatrix}$, and any element $u\in \L_n$ is a unique sum $u=
u^{ev}+u^{od}$ of its even and odd components.

{\it Proof}. Statement 1 is obvious. Note that 
\begin{equation}\label{sact}
\s (x) = \o_{1+ a} \g_b(Ax) = \o_{1+a} (A(x+b))= A(x+b) +
2aA(x+b).
\end{equation}
 Then $\s(x)^{od}
= A(x+b)$ and $\s (x)^{ev}=2aA(x+b)$. The first equality is
equivalent to statement 2, and the second equality can be
rewritten as follows,
$$ - \frac{1}{2}A^{-1} \s (x)^{ev} = (x+b) a = \g_b(x) a = \g_b(
x\g_b^{-1} (a)),$$ or, equivalently, $x\g_b^{-1} (a) = -
\frac{1}{2}A^{-1}\g_b^{-1} ( \s (x)^{ev})$. This is the system of
equations
$$\begin{cases}
x_1\g_b^{-1} (a)=-\frac{1}{2}a_1', \\
x_2\g_b^{-1} (a)=-\frac{1}{2}a_2', \\
\;\;\;\; \;\;\;\vdots \\
x_n\g_b^{-1} (a)=-\frac{1}{2}a_n'.
\end{cases}
$$
Its solutions are given by Theorem \ref{s14Sep06},
\begin{equation*}
a=-\frac{1}{2} \g_b(\sum_{i=1}^{n-1} x_1\cdots x_i\der_i\cdots
\der_1\der_{i+1}(a_{i+1}')+ \der_1(a_1')) +a_nx_1\cdots x_n, \; \;
a_n\in K,
\end{equation*}
where we have used the fact that $\g_b(a_nx_1\cdots x_n) =
a_nx_1\cdots x_n$. Since $ a_nx_1\cdots x_n = a+\frac{1}{2}
\g_b(\ldots )\in \L_n'^{od}$ we must have $a_n=0$, hence statement
3 holds (where $(\ldots )$ are the elements in the bracket above).
$\Box $

Department of Pure Mathematics

University of Sheffield

Hicks Building

Sheffield S3 7RH

UK

email: v.bavula@sheffield.ac.uk

\end{document}